\documentclass[english]{amsart}
\usepackage{amssymb,amsmath,amsfonts,amsthm}
\usepackage{stmaryrd}
\usepackage[all,cmtip,poly]{xy}
\usepackage{svn}
\usepackage{enumerate}
\usepackage{hyperref}

\SVN $LastChangedDate: 2013-05-15 11:39:46 -0400 (Wed, 15 May 2013) $
\SVN $LastChangedRevision: 435 $

\newcommand{\WBDJ}{W^{\textrm{BDJ}}}
\newcommand{\Wcris}{W^{\cris}}
\newcommand{\Wexplicit}{W^{\operatorname{explicit}}}
\newcommand{\WBT}{W^{\operatorname{BT}}}

\newcommand{\Ghat}{\hat G}
\newcommand{\hodgetateembedding}{{\kappa}}

\swapnumbers

\newcommand{\para}{\subsubsection{}}

\newtheorem{thm}[subsubsection]{Theorem}
\newtheorem{lemma}[subsubsection]{Lemma}
\newtheorem{lem}[subsubsection]{Lemma}
\newtheorem{cor}[subsubsection]{Corollary}

\newtheorem{prop}[subsubsection]{Proposition}

\newtheorem{ithm}{Theorem}

\theoremstyle{definition}

\newtheorem{defn}[subsubsection]{Definition}

\theoremstyle{remark}
\newtheorem{remark}[subsubsection]{Remark}
\newtheorem{rem}[subsubsection]{Remark}

\newtheorem{example}[subsubsection]{Example}

\def\numequation{\addtocounter{subsubsection}{1}\begin{equation}}
\def\nummultline{\addtocounter{subsubsection}{1}\begin{multline}}
\def\anumequation{\addtocounter{subsection}{1}\begin{equation}}

\newcommand{\F}{\FF}

\newcommand{\Q}{\QQ}

\newcommand{\Z}{\ZZ}

\newcommand{\D}{\cD}

\renewcommand{\O}{\cO}

\newcommand{\m}{\frakm}

\newcommand{\rbar}{\bar{r}}

\newcommand{\varepsilonbar  }{\overline{\varepsilon}}

\newcommand{\kappabar     }{\overline{\kappa}}

 \newcommand{\rhobar   }{{\overline{\rho}}}

 \newcommand{\chibar   }{\overline{\chi}}

\newcommand{\betat         }{\widetilde{\beta}}

\newcommand{\FF}{{\mathbb F}}

\newcommand{\QQ}{{\mathbb Q}}

\newcommand{\ZZ}{{\mathbb Z}}

\newcommand{\cD}{{\mathcal D}}

\newcommand{\cI}{{\mathcal I}}

\newcommand{\cO}{{\mathcal O}}

\newcommand{\cS}{{\mathcal S}}

\newcommand{\frakm}{\mathfrak{m}}

\newcommand{\Fbar}{\overline{\F}}
\newcommand{\Qbar}{\overline{\Q}}

\newcommand{\Fp}{\F_p}

\newcommand{\Fpbar}{\Fbar_p}

\newcommand{\Zp}{\Z_p}

\newcommand{\Ql}{\Q_{\ell}}
\newcommand{\Qp}{\Q_p}

\newcommand{\Qpbar}{\Qbar_p}

\newcommand{\psBT}{\psi\mathrm{BT}}
\newcommand{\calE}{\mathcal{E}}
\newcommand{\EpsBT}{\calE_{\psBT}}
\newcommand{\barP}{\overline{\mathfrak{P}}}

\DeclareMathOperator{\Ext}{Ext}
\DeclareMathOperator{\Fil}{Fil}

\DeclareMathOperator{\GL}{GL}

\DeclareMathOperator{\Hom}{Hom}

\DeclareMathOperator{\Sym}{Sym}

\newcommand{\ab}{\mathrm{ab}}
\newcommand{\cris}{\mathrm{cris}}

\newcommand{\HT}{\mathrm{HT}}

\newcommand{\llb}{\llbracket}
\newcommand{\rrb}{\rrbracket}
 
                               \newcommand{\into}{\hookrightarrow}
\newcommand{\onto}{\twoheadrightarrow}
\newcommand{\toisom}{\buildrel\sim\over\to}

\newcommand{\longnto}[1]{\buildrel#1\over\longrightarrow}
\newcommand{\lto}{\longrightarrow}

\newcommand{\ue}{\underline \epsilon}

\newcommand{\col}{\colon}

\newcommand{\fM}{{\mathfrak M}}
\newcommand{\fS}{{\mathfrak S}}

\newcommand{\M}{{\mathfrak M}}

\newcommand{\e}{{\mathfrak e}}

\newcommand{\hR}{\widehat{\mathcal R}}

\newcommand{\mf}{\mathfrak}

\newcommand{\Art}{{\operatorname{Art}}}
\newcommand{\barM}{\overline{\M}}
\newcommand{\barN}{\overline{\mathfrak{N}}}
\newcommand{\barK}{\overline{{K}}}
\newcommand{\barQQ}{\overline{{\Q}}}
\newcommand{\ve}{\varepsilon}

\newcommand{\barMp}{{\overline{\M'}}}

\newcommand{\gr}{{\rm gr}}

\newcommand{\fijd}{{\Fil ^{\{m_{ij }\}}\D}}

\newcommand{\f}{{\mf{f}}}
\newcommand{\U}{\mathrm{U}}

\begin{document}
\title[The weight part of Serre's conjecture for $\GL(2)$]{The weight part of Serre's conjecture for $\GL(2)$}
\author{Toby Gee} \email{toby.gee@imperial.ac.uk} \address{Department of
  Mathematics, Imperial College London}\author{Tong
  Liu}\email{tongliu@math.purdue.edu}\address{Department of
  Mathematics, Purdue University}\author{David Savitt} \email{savitt@math.arizona.edu}
\address{Department of Mathematics, University of Arizona}
\thanks{The first author is
  supported in part by a Marie Curie Career Integration Grant, and by an
  ERC Starting Grant. The second author was partially
  supported by NSF grant  DMS-0901360; the third author was partially
  supported by NSF  grant DMS-0901049 and NSF CAREER grant
  DMS-1054032.}  

\begin{abstract}Let $p>2$ be prime. We use purely local methods to
  determine the possible reductions of certain two-dimensional
  crystalline representations, which we call ``pseudo-Barsotti--Tate
  representations'', over arbitrary finite extensions of $\Qp$. As a
  consequence, we establish (under the usual Taylor--Wiles hypothesis)
  the weight part of Serre's conjecture for $\GL(2)$ over arbitrary
  totally real fields.
\end{abstract}
\maketitle

\setcounter{tocdepth}{1}
\tableofcontents

\section*{Overview}

Let $p$ be a prime number. Given an irreducible modular
representation $\rhobar:G_\Q\to\GL_2(\Fpbar)$, the weight part of
Serre's conjecture predicts the set of weights $k$ such that $\rhobar$
is isomorphic to the mod $p$ Galois representation $\rhobar_{f,p}$
associated to some eigenform of weight $k$ and level prime to $p$. The
conjectural set of weights is determined by the local representation
$\rhobar|_{G_{\Qp}}$. In recent years, beginning with the work
of~\cite{bdj}, generalisations of the weight part of Serre's
conjecture have become increasingly important, in particular because
of their importance in formulating a $p$-adic Langlands correspondence
(\emph{cf.}~\cite{MR2931521}).

The weight part of Serre's original conjecture was settled in the
early 1990s (at least if $p>2$;
see~\cite{MR1185584,MR1176206,GrossTameness}). The paper~\cite{bdj}
explored the generalisation of the weight part of Serre's original
conjecture~\cite{MR885783} to the setting of Hilbert modular forms
over a totally real field $F$ in which $p$ is unramified. Already in
this case even formulating the conjecture is far more difficult; there
are many more weights, and the conjectural description of them
involves subtle questions in integral $p$-adic Hodge theory. The
conjecture was subsequently extended to arbitrary totally real fields
in~\cite{ScheinRamified,GeePrescribed,blggu2}; see~\cite[\S 4]{blggu2}
for an extensive discussion. The present paper completes the proof of
this general Serre weight conjecture for $\GL(2)$ (under the assumption
that $p>2$ and that a standard Taylor--Wiles hypothesis holds).

Work of \cite{blggu2} and \cite{GeeKisin,Newton} has reduced the weight part of
Serre's conjecture for $\GL(2)$ to a certain statement about local Galois
representations:\ namely, two sets of Serre weights associated to a mod $p$
local Galois representation, each defined in $p$-adic Hodge-theoretic terms,
must be seen to be equal.  It is this problem that we resolve in this paper. In
our earlier paper~\cite{GLS-BDJ} we proved the same result for totally real
fields in which $p$ is unramified, by establishing a structure theorem for Kisin
modules which slightly extended Fontaine--Laffaille theory in this case. In
contrast, in this paper we must work over an arbitrarily ramified base field,
and the proof requires a delicate analysis of the Kisin modules associated to a
certain class of crystalline representations, which we call
pseudo-Barsotti--Tate representations. To the best of our knowledge, these are
the first general results about the reductions of crystalline representations in
any situation where arbitrary ramification is permitted.

\section{Introduction}

We begin by
tracing the history of the work on the weight part of Serre's conjecture over the last decade.
The first breakthrough towards proving the conjecture of~\cite{bdj} was the
paper~\cite{geebdj}, which proved the conjecture under a mild global
hypothesis as well as a genericity hypothesis on the local mod $p$
representations. The global hypothesis comes from the use of the
Taylor--Wiles--Kisin method, which is used to prove modularity lifting
theorems; the modularity lifting theorems used in~\cite{geebdj} are
those of~\cite{KisinModularity,MR2280776} for potentially
Barsotti--Tate representations. The genericity hypothesis was needed
for a complicated combinatorial argument relating  Serre weights to the reduction modulo
$p$ of the types associated to certain potentially Barsotti--Tate
representations. The natural output of the argument
was a description of the Serre weights (in this generic setting) in
terms of potentially Barsotti--Tate representations, while the
conjecture of~\cite{bdj} is in terms of crystalline representations,
and the comparison between the two descriptions involved a delicate
calculation in integral $p$-adic Hodge theory. It was clear that the
combinatorial arguments would not extend to cover the non-generic case, or to
settings in which $p$ is allowed to ramify.

All subsequent results on the problem have followed~\cite{geebdj} in
making use of modularity lifting theorems, and have assumed that
$p>2$, which we do for the remainder of this introduction. The next progress (other than special
cases such as~\cite{geesavitttotallyramified}) was due
to the work of~\cite{BLGGT}, which proved automorphy lifting theorems
for unitary groups of arbitrary rank, in which the weight of the
automorphic forms is allowed to vary. As the conjectural description
of the set of weights is purely local, there is a natural analogue of
the weight part of Serre's conjecture for inner forms of $\U(2)$, and the paper~\cite{blggu2} used
the results of~\cite{BLGGT} to prove a general result on these
conjectures. Specifically, given a global modular representation
$\rhobar$ associated to some form of $\U(2)$, the main results
of~\cite{blggu2} show that $\rhobar$ is modular of all the weights
predicted by the generalisations of the weight part of Serre's
conjecture. The problem is then to prove that $\rhobar$ cannot be
modular of any other weight.

This remaining problem is easily reduced to a local problem. To
describe this, let $K/\Qp$ be a finite extension with residue field $k$, and
$\rbar:G_K\to\GL_2(\Fpbar)$ a continuous representation. In this local
context a Serre weight is an irreducible $\Fpbar$-representation of
$\GL_2(k)$ (\textit{cf.}\ Definition~\ref{defn:serre weight}).
Associated to $\rbar$ there are two sets of Serre weights
$\Wexplicit(\rbar)\subseteq\Wcris(\rbar)$, defined in
terms of $p$-adic Hodge theory.  The set $\Wcris(\rbar)$ is defined in terms of
the existence of crystalline lifts with certain Hodge--Tate weights,
and $\Wexplicit(\rbar)$ is defined as a subset of $\Wcris(\rbar)$ by
explicitly writing down examples of these crystalline lifts (as
inductions and extensions of crystalline characters); see
Definitions~\ref{defn: Wcris} and \ref{defn: Wexplicit} below for precise definitions. These
definitions are then extended to global representations $\rhobar$ by
taking the tensor products of the sets of weights for the restrictions of
$\rhobar$ to decomposition groups at places dividing $p$. In the case
that $p$ is unramified or $\rbar$ is semisimple, $\Wexplicit(\rbar)$
is the set of weights predicted by the conjectures
of~\cite{bdj,ScheinRamified}.

If $\rhobar$ is a global modular representation for $\U(2)$, then it is
almost immediate from the definition that the set of weights in which
$\rhobar$ is modular is contained in $\Wcris(\rhobar)$. The main
result of~\cite{blggu2} shows that the set of weights contains
$\Wexplicit(\rhobar)$. Therefore, to complete the proof of the weight
part of Serre's conjecture for $\U(2)$ (i.e.\ to prove that
$\Wexplicit(\rhobar)$ is the
set of weights in which
$\rhobar$ is modular) it is only necessary to show in
the local setting
that $\Wexplicit(\rbar)=\Wcris(\rbar)$.

In the case that $K/\Qp$ is unramified (that is, the setting of the
original conjecture of~\cite{bdj}), this was proved by purely local
means in our earlier paper~\cite{GLS-BDJ}. That work uses the second
author's theory of $(\varphi,\Ghat)$-modules to prove a structure
theorem for the Kisin modules associated to crystalline
representations (over an unramified base) with Hodge--Tate weights
just beyond the Fontaine--Laffaille range. A careful analysis of this
structure theorem, and of the ways in which these Kisin modules can be
extended to $(\varphi,\Ghat)$-modules, allowed us to compute the possible
reductions of the crystalline representations, and explicitly check
that they were all of the required form.

Until the present paper, the equality
$\Wexplicit(\rbar)=\Wcris(\rbar)$ was not known in any greater
generality. However, in the paper~\cite{GLS11} we were able to show
that if $p$ is totally ramified, then the set of modular weights for $\U(2)$
is exactly $\Wexplicit(\rhobar)$, by making a rather baroque global
argument using the results of~\cite{blggu2} and a comparison to
certain potentially Barsotti--Tate representations, motivated by the
approach of~\cite{geebdj}. This approach did not show that
$\Wexplicit(\rbar)=\Wcris(\rbar)$; as in~\cite{geebdj}, this approach
relies on combinatorial results that do not extend to the general
case, although generalisations subject to a genericity hypothesis  (and other
related results)  were proved using these techniques
in~\cite{DiamondSavitt}.

The results of~\cite{blggu2,GLS-BDJ,GLS11} only concerned the
analogues for $\U(2)$ of the conjecture of~\cite{bdj} and its
generalisations, which were formulated using modular forms on
quaternion algebras over totally real fields. The quaternion algebra
setting of~\cite{bdj} is a more
natural generalisation of the original conjectures of~\cite{MR885783},
but it is harder, because of a parity obstruction coming from the
units of the totally real field: algebraic Hilbert modular forms
necessarily have paritious weight. This means, for example, that there
are mod $p$ Hilbert modular forms of level prime to $p$ and some
weight which cannot be lifted to characteristic zero forms of the same
weight, and is one of the reasons for studying potentially
Barsotti--Tate lifts instead (which correspond to Hilbert modular
forms of parallel weight two).

The papers~\cite{GeeKisin,Newton} independently succeeded (by rather different
means) in transferring the results for $\U(2)$ described above to the
setting of quaternion algebras over totally real fields. In
particular, \cite{GeeKisin} defines a set of weights $\WBT(\rbar)$,
and shows (by global methods, using the results of~\cite{blggu2}) that
there are
inclusions \[\Wexplicit(\rbar)\subseteq\WBT(\rbar)\subseteq\Wcris(\rbar).\]
The definition of $\WBT(\rbar)$ can again be extended to global representations
in exactly the same manner as before,
and (in either the $\U(2)$ or quaternion algebra settings) the set
of weights in which $\rhobar$ is modular is always the set
$\WBT(\rhobar)$. In particular, this shows that the set of weights is
determined purely locally, and in order to complete the proof of the
weight part of Serre's conjecture, it would be enough to solve the
purely local problem of showing that $\WBT(\rbar)=\Wexplicit(\rbar)$.

Unfortunately, the definition of $\WBT(\rbar)$ is rather indirect,
being defined as a linear combination of the Hilbert--Samuel
multiplicities of certain potentially Barsotti--Tate deformation
rings. These rings have only been computed when $K=\Qp$
(\cite{SavittDuke}) or when $K/\Qp$ is unramified and $\rbar$ is generic
(\cite{breuillatticconj,BreuilMezardRaffinee}), and they appear to be extremely difficult to
determine in any generality.

The main local result of this paper (see Theorem~\ref{thm:cyclotomic})
is that
$\Wexplicit(\rbar)=\Wcris(\rbar)$ in complete generality (provided
that $p>2$). The proof is
purely local, and will be described below. As a consequence, we deduce
that $\WBT(\rbar)=\Wexplicit(\rbar)$ when $p>2$, and thus
we obtain the following theorem (see Theorem~\ref{thm:our main global
  result}).

\begin{ithm}\label{thm:our main global result - intro version}Let $p>2$ be
  prime, let $F$ be a totally real field, and let
  $\rhobar:G_F\to\GL_2(\Fpbar)$ be a continuous representation. Assume
  that $\rhobar$ is modular, that $\rhobar|_{G_{F(\zeta_p)}}$ is
  irreducible, and if $p=5$ assume further that the projective image
  of $\rhobar|_{G_{F(\zeta_p)}}$ is not isomorphic to $A_5$.

  For each place $v|p$ of $F$ with residue field $k_v$, let $\sigma_v$
  be a Serre weight of $\GL_2(k_v)$. Then $\rhobar$ is modular of
  weight $\otimes_{v|p}\sigma_v$ if and only if
  $\sigma_v\in\Wexplicit(\rhobar|_{G_{F_v}})$ for all $v$.
\end{ithm}
In particular, this immediately implies the generalisations of the
weight conjecture of~\cite{bdj} proposed
in~\cite{ScheinRamified,GeePrescribed}.

We remark that Newton and Yoshida (\cite{NewtonYoshida}, forthcoming) have
also established Theorem~\ref{thm:our main global result -
  intro version} for many weights $\otimes_{v|p} \sigma_v$ by methods that are rather different from ours,
namely via novel arguments involving the study of the special fibres of
integral models for Shimura curves.

We now turn to an outline of our approach. To prove that
$\Wexplicit(\rbar)=\Wcris(\rbar)$, it is necessary to study the
possible reductions of certain 2-dimensional crystalline
representations. The Hodge--Tate weights of these representations
(which were first considered in~\cite{GeePrescribed}) have a
particular form, and we call these representations
\emph{pseudo-Barsotti--Tate} representations, because most of their
Hodge--Tate weights (in a sense that we make precise in Definition~\ref{defn:pseudo-bt}) agree with those
of Barsotti--Tate representations. In the case that $K/\Qp$ is
unramified, these are exactly the representations considered
in~\cite{GLS-BDJ}, and our techniques are a wide-ranging extension of
those of that paper to a setting where $p$ may be arbitrarily
ramified.

The equality $\Wexplicit(\rbar)=\Wcris(\rbar)$ is by definition equivalent to
the statement that the possible reductions of pseudo-Barsotti--Tate
representations with given Hodge--Tate weights are in an explicit
list, namely the list of representations arising as the reductions of
pseudo-Barsotti-Tate representations (of the same Hodge--Tate weights)
which are extensions or inductions of crystalline characters. The most
obvious way to try to establish this would be to classify (the lattices in)
pseudo-Barsotti--Tate representations (perhaps in terms of their
associated $(\varphi,\Ghat)$-modules), and then to compute all of
their reductions modulo $p$. Experience suggests that this is likely
to be a very difficult problem, and this is not the approach that we
take. Instead, we proceed more indirectly, guided by the particular
form of $\Wexplicit(\rbar)$.

Our first step is to make a detailed study of the filtrations on the
various objects in integral $p$-adic Hodge theory which are attached
to lattices in pseudo-Barsotti--Tate representations. This study,
which is a considerable generalisation of much of the work carried out
in~\cite[\S4]{GLS-BDJ} in the unramified case, culminates in
Theorem~\ref{thm: structure}, which classifies the possible structures
of their underlying Kisin modules in terms of their Hodge--Tate
weights. It is then relatively straightforward to determine the
possible characters that can occur in the reduction modulo $p$ of such
a representation, and comparing this to the form of
$\Wexplicit(\rbar)$, we prove the conjecture in the case that $\rbar$
is a direct sum of two characters (see
Theorem~\ref{thm:reducible-case}).

We next turn to the case that $\rbar$ is irreducible. In this case, we
know that $\rbar$ becomes a direct sum of two characters after
restriction to the absolute Galois group of the unramified quadratic
extension, and applying the previous result in this case gives a
constraint on the form of $\rbar$. It is not at all obvious that this
constraint is sufficient to prove the conjecture in this case, but we
are able to establish this by a somewhat involved combinatorial
argument (see Theorem~\ref{thm:irreducible-case}).

At this point we have established the conjectures
of~\cite{ScheinRamified}, which treat the case that $\rbar$ is
semisimple. There is, however, still a considerable amount of work to
be done to deal with the case that $\rbar$ is an extension of
characters. It is not surprising that more work should be needed in
the non-semisimple case. For example, if the ramification degree of $K$ is at least
$p$ and $\rbar$ is semisimple, then the inclusion $\Wcris(\rbar) \subseteq \Wexplicit(\rbar)$ is
essentially trivial,
because the set $\Wexplicit(\rbar)$ consists of all the weights
satisfying a simple determinant condition. On the other hand
when~$\rbar$ is not semisimple,
the definition of $\Wexplicit(\rbar)$,
which is given in terms of certain crystalline $\Ext^1$ groups, only
becomes more complicated as the ramification index increases.

In the remaining case that $\rbar$ is an extension of characters, the
result that we proved in the semisimple case shows that the two
characters are of the predicted form, and it remains to show that the
extension class is one of the predicted extension classes. We again proceed indirectly.
From the definition of $\Wexplicit(\rbar)$, we need
to show that $\rbar$ has a pseudo-Barsotti--Tate lift of the given
weight and which is an extension of crystalline characters. An approach to
this problem naturally suggests itself: we could compute the dimension
of the space of extensions in characteristic $p$ that arise from the
reductions of extensions of crystalline characters, and try to use our
structure theorem to prove that the set of extensions that can arise
as the reductions of possibly irreducible pseudo-Barsotti--Tate
representations is contained in a space of this same
dimension.

This is in effect what we do, but there are a number of serious
complications that arise when we try to
compute our upper bound on the set of extension classes coming from
pseudo-Barsotti--Tate representations. It is natural to return to our
structural result Theorem~\ref{thm: structure}, and this gives us
non-trivial information on the possible Kisin modules, which we
assemble with some effort.  We would then like to compare these Kisin modules with those arising
from extensions of crystalline characters. It is at this point that
two difficulties arise. One is that the functor from lattices
in Galois representations to Kisin modules is not exact (see
Example~\ref{ex:structure-of-argumemt} and the discussion that follows
it), which means that the Kisin
module corresponding to the reduction of such an extension need not
correspond to an extension of the corresponding rank-one Kisin
modules. The other related difficulty is that when the ramification
degree is large, there are many different crystalline characters with
different Hodge--Tate weights to consider, and it is necessary to
relate their reductions. We are able to overcome these difficulties by showing that there
is a ``maximal'' pair of rank-one Kisin modules for the representation
$\rbar$ and the Hodge--Tate
weights under consideration, and reducing to the
problem of studying their extensions.  This is done in
Sections~\ref{sec:comparison} to \ref{sec:weight-part-serres}.

The second complication is that Kisin modules do not completely determine the corresponding
$G_K$-representations, but rather their restrictions to a certain
subgroup $G_{K_\infty}$. If $\rbar$ is an extension of $\chibar_2$ by
$\chibar_1$ with $\chibar_1\chibar_2^{-1}$ not equal to the mod $p$
cyclotomic character, it turns out that the natural restriction map
from extensions of $G_K$-representations to extensions of
$G_{K_\infty}$-representations is injective (Lemma~\ref{lem:
  injectivity of restriction on H^1}), and we have done enough to
complete the proof.  However, in the remaining case that $\chibar_1\chibar_2^{-1}$ is the
mod $p$ cyclotomic character, we still have more work to do; we need
to study the uniqueness or otherwise of the extensions of the Kisin
modules to $(\varphi,\Ghat)$-modules. We are able to do this with the
aid of~\cite[Cor.~4.10]{GLS-BDJ}, which constrains the possible
$\Ghat$-actions coming from the reductions of crystalline
representations.

\subsection{Acknowledgements} We would like to thank Florian Herzig and
the anonymous referees for
their helpful comments.

\subsection{Notation}

\subsubsection{Galois theory}
\label{sec:galois-theory}
If $M$ is a field, we let $G_M$ denote its
absolute Galois group.
If~$M$ is a global field and $v$ is a place of $M$, let $M_v$ denote
the completion of $M$ at $v$.   If
~$M$ is a finite extension of $\Ql$ for some $\ell$, we let $M_0$
denote the maximal unramified extension of $\Ql$ contained in $M$, and
we write $I_M$
for the inertia subgroup of~$G_M$. If $R$ is a local ring we write
$\mf{m}_{R}$ for the maximal ideal of $R$.

 Let $p$ be a prime number.   Let $K$ be a finite extension of $\Qp$, with ring of integers $\cO_K$
 and residue field~$k$.   Fix a uniformiser $\pi$ of $K$, let
 $E(u)$ denote the minimal polynomial of $\pi$ over $K_0$, and set $e
 = \deg E(u)$.  We also fix an algebraic closure $\barK$ of $K$.  The
 ring of Witt vectors $W(k)$ is the ring of integers in $K_0$.

Our representations of $G_K$ will have coefficients in subfields of $\Qpbar$,
another fixed algebraic closure of $\Qp$, whose residue field we denote $\Fpbar$.   Let $E$ be a finite
extension of $\Qp$ contained in $\Qpbar$ and containing the image of every
embedding into $\Qpbar$ of the unramified quadratic extension of $K$; let $\O_E$ be the ring of integers in
$E$, with uniformiser $\varpi$ and residue field $k_E \subset \Fpbar$.

We write $\Art_K \col K^\times\to W_K^{\ab}$ for
the isomorphism of local class field theory, normalised so that
uniformisers correspond to geometric Frobenius elements.  For each $\lambda\in \Hom(k,\Fpbar)$ we
define the fundamental character $\omega_{\lambda}$ corresponding
to~$\lambda$ to be the composite 
$$ I_K \lto W_K^{\ab} \longnto{\Art_K^{-1}} \O_K^{\times} \lto k^{\times} \longnto{\lambda}
\Fpbar^{\times}.$$ 

We fix a compatible system of $p^n$th roots of $\pi$: that is, we set
$\pi_0 = \pi$ and for all $n  > 0$ we fix a choice of $\pi_n$
satisfying $\pi_{n}^p = \pi_{n-1}$.  Define $K_\infty = \cup_{n=0}^\infty K(\pi_n).$

\subsubsection{Hodge--Tate weights}
\label{sec:p-adic-hodge}

If $W$ is a de Rham representation of $G_K$ over
$\barQQ_p$ and $\hodgetateembedding$ is an embedding $K \into \barQQ_p$ then the multiset
$\HT_\hodgetateembedding(W)$ of Hodge--Tate weights of $W$ with respect to $\hodgetateembedding$ is
defined to contain the integer $i$ with multiplicity $$\dim_{\barQQ_p} (W
\otimes_{\hodgetateembedding,K} \widehat{\barK}(-i))^{G_K},$$ with the usual notation
for Tate twists. (Here $\widehat{\barK}$ is the completion of $\barK$.)  Thus for example
$\HT_\hodgetateembedding(\varepsilon)=\{ 1\}$, where $\varepsilon$ is
the cyclotomic character. We will refer to the
elements of $\HT_\hodgetateembedding(W)$  as the
``$\hodgetateembedding$-labeled Hodge--Tate weights of $W$'', or simply as the
``$\hodgetateembedding$-Hodge--Tate weights of $W$''.

\subsubsection{$p$-adic period rings}
\label{sec:p-adic-period-rings}
Define $\fS = W(k)\llbracket u\rrbracket$.  The ring $\fS$
is equipped with a Frobenius endomorphism $\varphi$ via $u\mapsto u^p$
along with the natural Frobenius on $W(k)$.

We denote by $S$  the $p$-adic completion of the divided power
envelope of $W(k)[u]$ with respect to the ideal generated by $E(u)$.
Let $\Fil^r S$ be the closure in $S$ of the ideal generated by
$E(u)^i/i!$ for $i \ge r$.   Write $S_{K_0} = S[1/p]$ and
$\Fil^r S_{K_0} = (\Fil^r S)[1/p]$.
There is a unique Frobenius map $\varphi\col S \to S$ which extends
the Frobenius on $\fS$. We write $N_S$ for the $K_0$-linear derivation on $S_{K_{0}}$ such that
$N_S(u)= -u$.

\section{Kisin modules attached to pseudo-Barsotti--Tate
  representations}

\subsection{Finer filtrations on Breuil modules}

Let $V$ be a $d$-dimensional $E$-vector space with continuous
$E$-linear $G_K$-action which makes $V$ into a crystalline
representation of $G_K$.  Let $D := D_{\cris}(V)$ be the filtered
$\varphi$-module associated to~$V$ by Fontaine~\cite{AST223}.  Recall that we have assumed that
the coefficient field $E$ contains the image of every embedding of
$K$ into $\Qpbar$. Since $V$ has an $E$-linear structure, it turns out
that the filtration on the Breuil module $\D$ attached to $D$ in \cite[\S6]{BreuilGriffith} can be
endowed with
finer layers that encode the information of the
$\hodgetateembedding$-labeled Hodge--Tate weights of $D$. As we will
see in  later subsections, these finer filtrations play a crucial
role in understanding the structure of Kisin modules (\cite{KisinCrys}) associated to
pseudo-Barsotti--Tate representations (Theorem \ref{thm: structure}).  In this subsection, we construct such finer filtrations and study their basic properties.

\para Let $f = [K_0 : \Q_p]$ and recall that $e = [K: K_0]$.  Fix an element $\kappa_0
\in \Hom_{\Q_p} (K_0,E)$, and recursively define $\kappa_i \in
\Hom_{\Qp}(K_0,E)$ for $i \in \Z$ so that $\hodgetateembedding_{i+1}^p \equiv \hodgetateembedding_i
\pmod{p}$.  Then the distinct elements of
$\Hom_{\Qp}(K_0,E)$ are $\kappa_0,\ldots,\kappa_{f-1}$, while $\kappa_f = \kappa_0$.
Now we label the elements of $\Hom_{\Qp}(K,E)$ as
$$\Hom_{\Q_p}(K, E)= \{\kappa_{ij} \, : \,  i = 0, \dots , f-1 ,\, j=
0, \dots, e-1 \}$$ in any manner so that $\kappa_{ij}|_{K_0} =
\kappa_i$.

Let $\ve_i \in W(k) \otimes_{\Zp} \O_E$ be the unique idempotent
element such that $(x \otimes 1)\ve_i = (1 \otimes \hodgetateembedding_i(x))\ve_i$
for all $x \in W(k)$.  Then we have $\ve_i(W(k) \otimes_{\Zp} \O_E)
\simeq W(k) \otimes_{ W(k), \kappa_i} \O_E$.

Write $K_{0, E} = K_0 \otimes_{\Q_p} E$ and ${K}_E = K \otimes_{\Q_p}
E$. We have a natural decomposition of rings  $K_{0,E} = \prod
_{i=0}^{f-1} \ve_i ({K_{0, E }}) \simeq \prod _{i=0}^{f-1} E_{i}$
where $E_{i}= \ve_i ({K_{0, E}}) \simeq  K_0 \otimes_{K_0, \kappa_i}
E$. Similarly we have $K_E  = \prod_{i,j} E_{ij}$ with $E_{ij} = K
\otimes_{K , \kappa _{ij}} E$.  We will sometimes identify an element
$x \in E$ with an element of
$E_i$ via the map $x \mapsto 1 \otimes x$, and similarly for $E_{ij}$.

\para Recall that $D$ is a finite free $K_{0, E}$-module of rank $d$ (see
\cite[\S3.1]{BreuilMezard}, for example), so that $D_K:= K \otimes_{K_0} D$ is a finite free $K_E$-module of rank $d$.
 In particular we have natural decompositions $D = \oplus_{i =0}^{f-1} D_{i}$ with $D_i = \ve_i (D) \simeq D \otimes_{K_{0,E} }E _{i}$, and $D_K = \oplus_{i,j} D_{K, ij}$ with $D_{K,ij} = D \otimes_{K_E} E_{ij}$.
Note that $\Fil^m D_K$ and ${\rm gr} ^m D_K$ have similar
decompositions, though of course they need not be free. Given a
multi-set of integers $\{m_{ij}\}$ with $0 \le i \le f-1$ and $0 \le j
\le e-1$, we define
\[ \Fil ^{\{m_{ij}\}} D_K := \bigoplus_{i=0}^{f-1}
\bigoplus_{j=0}^{e-1} \Fil ^{m_{ij}}D_{K, ij}\subset D_K.
\]
If $m_{ij} = m$ for all $i,j$ then of course $\Fil ^{\{m_{ij}\}}D_K = \Fil ^m D_K$.

\para Let $A$ be an $\fS$-algebra. We write $A_{\O_E} : = A \otimes_{\Z_p}
\cO_E$.  This is a $W(k) \otimes_{\Z_p}\O_E$-algebra, so that
$A_{\O_E} \simeq \prod_{i =0}^{f-1} A_{ \O_{E},i}$ with $A_{\O_{E},i} =
\ve _i (A_{\O_E})  \simeq A  \otimes_{W(k), \kappa_i} \O_E$.
Similarly, we  write $A_E := A \otimes_{\Z_p} E$, so that we have $A_E
\simeq \prod_{i =0}^{f-1} A_{E,i}$ with $A_{E,i} :=  A \otimes_{
  W(k),\kappa_i} E$.

Let us write $\iota$ for the isomorphism $\fS_{\O_E}
\simeq   \prod_{i =0}^{f-1} \fS_{ \O_{E},i}$ and $\iota_i$ for the
projection  $\fS_{\O_E} \twoheadrightarrow \fS_{\O_{E},i} $. Note that
$\fS_{\O_{E},i} = \fS \otimes_{W(k),\kappa_i} \O_E  \simeq
\O_{E_i}\llb u \rrb$, where $\O_{E_i}$ denotes the ring of integers in
$E_i$.

Recall that  $\pi$ is a fixed uniformiser of $K$, with minimal
polynomial $E(u)$ over $K_0$.  Write $\pi_{ij}= \kappa_{ij}(\pi) \in
E$. For each $\kappa_i$, we define $E^{\kappa_i} (u) =
\prod_{j=0}^{e-1} (u - \pi_{ij})$ in  $E[u]$,  so that
$E^{\kappa_i}(u)$ is just the polynomial obtained by acting by $\kappa_i$
on the coefficients of $E(u)$; note that identifying $E_i$ with $E$ will
identify $\iota_i(E(u))$ with  $E^{\kappa_i} (u)$.

Let $f_\pi$ be the map $S \to \O_K$ induced by $u \mapsto \pi$.  We
will also write $f_{\pi}$ for the map  $f_\pi \otimes_{\Z_p} E : S_ E
\to (\O_K)_E$.  We have surjections $\iota_{ij}:  (\O_K)_E \to K_E \to
E_{ij}$ for all $i,j$, and composing with $f_{\pi}$ gives $E$-linear
maps $f_{ij}:= \iota_{ij} \circ f_\pi: S_E \to E_{ij}$. Restricting
$f_{ij}$ to $\fS_{\O_E}$ gives an $\O_E$-linear surjection $\fS_{\O_E}
\to \O_{E_{ij}}$ that we also denote $f_{ij}$.  (Here $\O_{E_{ij}}$
denotes the ring of integers in $E_{ij}$.)  Set $\Fil ^1_{ij} S_E :=
\ker(f_{ij})$ and $\Fil ^1 \fS_{\O_E} := \fS_{\O_E} \cap \Fil ^1_{ij}
S_E$. Let $E_{ij }(u)$ be the unique element in $\fS_{\O_E}$ such that
$\iota_{\ell} (E_{ij}(u))=  u - \pi _{ij}$ if $\ell= i$ and $\iota_{\ell}
(E_{ij}(u)) =1$ if $\ell \not =i$.  We see that $E(u) \otimes 1 = \prod_{i,j}
E_{ij}(u)$ in $\fS_{\O_E}$.

\begin{lemma} \label{lem: fil1ij}\hfill
\begin{enumerate} \item $E_{ij}(u) \in \Fil ^1_{ij}S_E$.
\item $\bigcap_{ij} \Fil ^1 _{ij}S_E = \Fil ^1 S \otimes_{\Z_p} E. $
\item $\Fil ^1 _{ij}\fS_{\O_E} = E_{ij} (u) \fS_{\O_E}.$
\end{enumerate}
\end{lemma}
\begin{proof}  Note that $S_{E,i} \subset K_0 \llb u \rrb \otimes_{K_0,
    \kappa_i} E \simeq E_i \llb u \rrb$, so that elements in $S_{E,i}$ can be
  regarded as power series with $E_i$-coefficients.

Unwinding the definitions, we see that $f_{ij} : S_E \to
  E_{ij}$ is the $E$-linear map sending $u$ to $\pi_{ij}$.  In fact
  the map $f_{ij}$ can be factored as
$$ f_{ij} : S_E \simeq \prod_{i'=0}^{f-1} S_{E,i'} \to \prod_{i'=0}^{f-1} (\cO_K)_{E,i'} \to
E_{ij}$$
where the second map is the product of the maps sending $u$ to $\pi
\otimes 1$, while the third map is the map $\cO_K
\otimes_{W(k),\kappa_i} E \onto \cO_K \otimes_{\cO_K,\kappa_{ij}} E$
on the $i$th factor and is zero on the remaining factors.  Now (1) is
clear, and moreover an element $h \in S_{E,i}$ lies in $\Fil^1_{ij}
S_E$ if and only if $h = (u - \pi_{ij})h'$ for some $h' \in E_i\llb u
\rrb$.  (We caution the reader that $h'$ need not be in $S_{E,i}$.)
Then (2) and (3) follow easily.
\end{proof}

For $0 \le \ell \le f-1$, define $\Fil ^1_{ij} S_{E,\ell}:= \Fil
^1_{ij} S_{E} \cap S_{E,\ell}$ and $\Fil ^1 _{ij} \fS_{E,\ell} := \Fil
^1 _{ij} \fS_{\O_E} \cap  \fS_{\O_{E,\ell}}$. Note that unless $\ell =
i$ we have $\Fil^1_{ij} S_{E,\ell} = S_{E,\ell}$, and similarly for $\fS_{E,\ell}$.

\para Let $N$ be the $W(k)$-linear differential operator on $\fS$ such that
 $N(u) = -u$, and extend $N$ to $\fS_{\O_E}$ in the unique $E$-linear way. It is easy to check that $N$ is compatible with $\iota$ in the following sense: if $N$ is the  $\O_{E}$-linear differential operator on $\O_{E,i} \llb u\rrb  $ such that $N(u)= -u$,  then $\iota_i (N(x))  = N(\iota_i (x))$ for any $x \in \fS_{\O_E}$.

Let $\D := S \otimes_{W(k)} D$ be the Breuil module attached to
$D$  (see e.g.~\cite[\S6]{BreuilGriffith}).  We have a natural isomorphism  $ D_K \simeq
\D \otimes_{S,f_{\pi}} \O_K$; therefore, we also have a natural
isomorphism  $D_{K, ij} \simeq \D \otimes_{S_E, f_{ij}} E_{ij}$. We
again denote  the projection  $\D \twoheadrightarrow D_K$ by $f_\pi$,
and the projection $\D \twoheadrightarrow D_K \twoheadrightarrow D_{K,
  ij}$ by $f_{ij}$.

\begin{remark}\label{rem: relation fij and fpi} Since $D_K =
  \oplus_{i,j} D_{K, ij}$, it is easy to check  that $f_\pi =
  \oplus_{i,j} f_{ij}$. This is a useful fact in the $\Qp$-rational
  theory. Unfortunately, this  fails in general in the  integral
  theory (unless $e=1$), essentially because the
  idempotents in $K\otimes_{\Qp} E$ may not be contained in $\cO_K
  \otimes_{\Zp} \cO_E$.  While it is not hard to see that
  $f_{ij}(\fS_{\O_E}) = \O_{E_{ij}}$, the map  $f_\pi : \fS_{\O_E} \to
  \oplus_{i,j} \O_{E_{ij}}$  in general is not surjective: indeed $f_\pi (\fS_{\O_E}) = \O_K \otimes_{\Z_p} \O_E$, and
  this need not equal $\oplus_{i,j} \O_{E_{ij}}$ unless
  $e=1$.
\end{remark}

\para \label{para:filmij} Let  $\{ m_{ij}\} $ be a collection of integers indexed by $i = 0, \dots
  f-1$ and $j =0, \dots ,e-1$. We recursively define a filtration $\Fil ^{\{m_{ij}\}} \D\subseteq \D$. We first set $\Fil ^{\{m_{ij}\}} \D = \D$ if $m_{ij} \le 0$ for all $i,j$. Then define
\begin{equation*}\label{defn: filmij}
\Fil ^{\{m_{ij}\}} \D =\{x \in \D \, : \,f_{ij} (x) \in \Fil
^{m_{ij}}D_{K, ij}  {\rm\ for\ all\ } i,j, \text{ and }  N(x) \in \Fil
^{\{m_{ij}-1\}}\D    \}.
\end{equation*}
This is a direct generalization of the usual
filtration $\Fil ^m \D$ on $\D$ defined in~\cite[\S6]{BreuilGriffith}, and evidently $\Fil ^m \D= \Fil
^{\{m_{ij}\}}\D$ when  $m_{ij} = m$ for all $i,j$.

\para \label{para:filmijl} Let us discuss a slight variation of the filtration $\fijd$
defined in  \ref{para:filmij}. Recall that  $ S _{E}\simeq \prod_{\ell
  =0}^{f-1}  S_{E,\ell}$. Then $\D = \oplus_{\ell = 0}^{f-1} \D_\ell$
with $\D_\ell = \D \otimes_{S_{E}} S_{E,\ell}$. We also have $\fijd =
\oplus_{\ell  =0} ^ {f-1}\Fil ^{\{m_{ij}\}} \D_{\ell  } $ with $\Fil
^{\{m_{ij}\}} \D_{\ell  } = \fijd \otimes_{S_E} {S_{E,\ell}}$. Note
that  $ N(\D_\ell) \subset \D_\ell$ (because $N$ is $E$-linear on
$D$). We see easily that $ \Fil ^{\{m_{i  j}\}}\D_\ell$ depends only on the $m_{\ell
  j}$ for $0 \le j \le e-1$, and  we can define
$ \Fil^{\{m_{\ell,0},\ldots,m_{\ell,e-1}\}} \D_{\ell} := \Fil ^{\{m_{i
    j}\}}\D_\ell$.  Note that $
\Fil^{\{m_{i,0},\ldots,m_{i,e-1}\}} \D_{i} $ also has the
recursive description
$$ \{x \in \D_i \, : \,
f_{i j} (x) \in \Fil ^{m_{i j}}D_{K, i j} \text{ for all $j$,
  and } N(x) \in \Fil^{\{m_{i,0}-1,\ldots,m_{i,e-1}-1\}}\D_i \}.$$

The following
proposition summarizes some useful properties of $\Fil
^{\{m_{ij}\}}\D$ and $\Fil^{\{m_{i,0},\ldots,m_{i,e-1}\}} \D_{i} $.

\begin{prop}\label{prop: properties of filmij} With notation as
 in \ref{para:filmij}, the filtration $\Fil^{\{m_{ij}\}} \D$ has the following properties.
\begin{enumerate}

\item If $m_{ij} \geq m'_{ij}$ for all $i,j$ then $\Fil ^{\{m_{ij}\}}\D \subseteq \Fil ^{\{m'_{ij}\}} \D$.

\item $N(\Fil ^{\{m_{ij}\}}\D) \subseteq \Fil ^{\{m_{ij}-1\}}\D$.

\item $\Fil^{\{m_{ij}\}}\D$ is an $S_E$-submodule of $\D$.

\item  $f_\pi( \fijd)  =  \oplus _{i,j } \Fil ^{m_{ij}} D_{K, ij}$.

\item $(\Fil^r S) \Fil^{\{m_{ij}-r\}} \D \subseteq \Fil^{\{m_{ij}\}} \D$ for
  all $r \ge 0$.

\item We have $(\Fil ^1 _{i' j'}S_E )\Fil ^{\{m_{ij}\}} \D \subseteq \Fil
  ^{\{m'_{ij }\}}\D $ with $m'_{i' j'}= m_{i' j'} +1$ and $m'_{ij} =
  m_{ij}$ if $(i, j) \neq (i', j')$.

\item Suppose that $E(u)x \in \Fil^{\{m_{ij}\}} \D$ with $x \in \D$.  Then $N^{\ell}(x) \in
  \Fil^{\{m_{ij}-1-\ell\}} \D$ for all $\ell \ge 0$.

\item Suppose that $E_{i' j'} (u) x \in \fijd $ with $x \in \D$. Then
  $x \in  \Fil ^{\{m'_{ij }\}}\D $ with $m'_{i' j' }= m_{i' j'} -1$
  and $m'_{ij} = m_{ij}$ if $(i, j) \not = (i', j')$.

\end{enumerate}
The analogous statements hold for the filtration
$\Fil^{\{m_{i,0},\ldots,m_{i,e-1}\}} \D_i$ of \ref{para:filmijl},
replacing $\D$ by $\D_i$ and $S_E$ by $S_{E,i}$.  (Note that we have defined $\Fil ^1_{i'j'} S_{E,i}$ below Lemma~\ref{lem: fil1ij}.)
\end{prop}

\begin{proof}
The proofs for $\Fil^{\{m_{i,0},\ldots,m_{i,e-1}\}} \D_i$ are essentially the
same as those for $\Fil^{m_{ij}} \D$, so we will concern ourselves
exclusively with the latter case.

(1) This is a straightforward induction on $m' =
\max_{i,j}\{m'_{ij}\}$, with the case $m' \leq  0$ serving as the base case. Let us suppose  that the
statement is true for $m'-1$ and consider the situation for $m'$.
If $x \in \Fil ^{\{m_{ij}\}}\D $, we certainly have $f_{ij} (x) \in
\Fil ^{m_{ij}}D_{K, ij} \subseteq \Fil ^{m'_{ij}} D_{K, ij}$,
so it suffices to show that $N(x) \in \Fil ^{\{m'_{ij}-1
  \}}\D$. But
$N(x) \in \Fil ^{\{m_{ij}-1\}}\D$, and $\Fil ^{\{m_{ij}-1\}}\D \subseteq
\Fil ^{\{m'_{ij}-1\}}\D   $ by induction,  which completes the proof.

(2) This is immediate from the definition.

For each of items (3)--(6), we proceed by induction on $m =
\max_{i,j}\{m_{ij}\}$.  These statements are all trivial if $m \leq
0$, except for (6) where we take $m < 0$ as the base case. For each of
these items in turn, let us suppose that the statement is true for $m-1$ and consider the situation for $m$.

(3) Pick $s \in S_E$ and $x \in \Fil ^{\{m_{ij}\}}\D$. It is clear
from the definitions that $f_{ij} (sx)= f_{ij} (s) f_{ij } (x) \in
\Fil ^{m_{ij}} D_{K, ij}$ for all $i,j$, so it remains to show that
$N(s x)  = N(s) x + s N(x)$ is inside $\Fil ^{\{m _{ij}-1
  \}}\D$. Since $N(x) \in \Fil ^{\{m_{ij}-1 \}}\D$, and $\Fil
^{\{m_{ij}-1\}}\D$ is an $S_E$-module by induction, we see that $ s
N(x) \in \Fil ^{\{m_{ij}-1\}}\D$. For the other term, we know from (1) that $\fijd \subseteq \Fil ^{\{m_{ij}-1 \}}\D$.
Therefore  $x$, and hence $N(s)x$, is in $\Fil ^{\{m_{ij}-1\}}\D $
(again using the induction hypothesis).

We now prove (6), (5), and (4) in that order.

(6) Let  $x = sy$ with $s \in \Fil ^1_{i'j'} S_E $ and $y \in \fijd$. We
note that $f_{ij} ( s y ) \in \Fil ^{m_{ij}}D_{K, ij}$ if $(i, j )
\not = (i', j' )$, and $f_{i'j'} (sy )= 0 \in \Fil ^{m_{i' j'} +1}
D_{K, i'j'}$. It remains to show that $N(x) = N(s) y + s N(y)$ lies in
$\Fil ^{\{m'_{ij} -1 \}}\D$.  We have
$N(s) y  \in \Fil ^{\{m_{ij}\}}\D  \subseteq  \Fil ^{\{m'_{ij}-1\}}\D$
by (3) and (1).  On the other hand $N(y) $ is in $\Fil ^{\{m_{ij}-1\}}\D$ by
definition, and so by the induction hypothesis we have $sN(y) \in \Fil
^{\{m'_{ij}-1 \}}\D $. This completes the induction.

(5) This follows by an argument essentially identical
to the proof of (6), using $f_{ij}(s)=0$ for all $i,j$ if $s \in \Fil^r S$, and that
$N(\Fil^r S) \subseteq \Fil^{r-1} S$ for $r \ge 1$.

(4) Pick  $x = (x_{ij}) \in \oplus_{i,j}\Fil ^{m_{ij}} D_{K, ij}$.
Since  $  x_{ij }\in \Fil ^{m_{ij}} D_{K, ij}\subseteq \Fil ^{m_{ij}-1} D_{K, ij}$, by induction there exists $\hat x \in \Fil ^{\{m_{ij}-1 \}}\D$ such that $f _{ij} (\hat x) = x_{ij} $ for all $i, j$.
Note that $N(E(u))$ is relatively prime to $E(u)$, so there exists
$R(u) \in K_0[u]$ such that $ E(u)$ divides $1 + R(u) N(E(u))$. Set $H(u): = R(u) E(u)$ and write $N(H(u)) +1 = E(u)Q(u)$.
Define \begin{equation*} \label{eq: lift to fild}
\hat y : = \sum_{\ell  =0 }^{m-1}\frac {H(u )^\ell  N^\ell  (\hat x)}{\ell !}.
\end{equation*}
It is easy to check that $f_{ij}(\hat y) = f_{ij }(\hat x ) = x_{ij}$ for all $i, j$. Now we have
\begin{eqnarray*}
N(\hat y ) & = &  N(\hat x ) + \sum_{\ell  = 1}^{m-1} \frac {1}{\ell!}\left (\ell H(u)^{\ell  -1} N(H(u)) N^\ell (\hat x)  + H(u)^\ell N^{\ell  + 1} (\hat x)\right ) \\
& = & \frac{H(u)^{m-1}}{(m-1)!}N^{m}(\hat x) + \sum_{\ell  =1}^{m-1} \frac {(1 + N(H(u))) H(u)^{\ell -1} N^{\ell }(\hat x)}{(\ell -1)! } \\
& = &  \frac{H(u)^{m-1}}{(m-1)!}N^{m }(\hat x) + \sum_{\ell  =1}^{m-1} \frac{Q(u) R(u)^{\ell -1} E(u)^{\ell } N^{\ell }(\hat x)}{(\ell-1)! }
\end{eqnarray*}
Since $\hat x \in \Fil^{\{m_{ij}-1\}}\D$, repeated applications of
(2) and (5) show that $E(u)^{\ell}  N^{\ell}  (\hat x) $ is in $\Fil
^{\{m_{ij}-1\}}\D$ for all $\ell$. Since $E(u)$ divides $H(u)$ we also
have
$H(u)^{m-1} \D \subseteq \Fil^{\{m-1\}} \D \subseteq
\Fil^{\{m_{ij}-1\}} \D$. By (5) we conclude that $\hat y  \in \fijd$ as required,
and the induction is complete.

(7) First we show that $f_{ij} (N^{\ell}(x)) \in \Fil
^{\{m_{ij}-1-\ell\}}D_{K, ij}$ for all $i,j$, proceeding by induction on $\ell$. Note that $N(E(u)x) = N(E(u))x
+ E(u)N(x)$ lies in $\Fil^{\{m_{ij}-1\}}\D$ and so $f_{ij} (N(E(u))x)
\in \Fil^{m_{ij}-1} D_{K, ij}$ for all $i,j$. Since $f_{ij}
(N(E(u)))\not = 0$ we see that $f_{ij} (x) \in \Fil^{m_{ij}-1} D_{K, ij}$ for all $i,j$. This proves the case $\ell =0$. In general, we have
\begin{equation*} N^{\ell+1}( E(u)x ) =  \sum_{n = 0} ^{\ell+1}
  \binom{\ell+1}{n} N^{\ell+1- n}(E(u)) N^{ n}(x) \in
  \Fil^{\{m_{ij}-1-\ell\}} \D,\end{equation*}
using (2) to deduce the containment.  After applying $f_{ij}$ to both sides of
the above equation, each term on the right-hand side with $n \neq \ell$
will be known to belong to $\Fil ^{\{m_{ij}-1-\ell\}}D_{K, ij}$. For
$n < \ell$ this follows from the induction hypothesis, while the
term with $n=\ell+1$ vanishes because $f_{ij}(E(u))=0$.  We conclude
that the same is true for the term with $n=\ell$, and so (again using
that $f_{ij}
(N(E(u)))\not = 0$) we deduce that $f_{ij} (N^{\ell}(x)) \in \Fil
^{\{m_{ij}-1-\ell\}}D_{K, ij}$.

Now, the claim $N^\ell(x) \in \Fil^{\{m_{ij}-1-\ell\}}\D$ is automatic for
  $\ell \ge \max_{i,j} \{m_{ij}\}$.  Suppose that we have $N^\ell(x) \in
  \Fil^{\{m_{ij}-1-\ell\}}\D$, or in other words $N(N^{\ell-1}(x)) \in
    \Fil^{\{m_{ij}-1-\ell\}}\D$.  Since we also have $f_{ij} (N^{\ell-1}(x)) \in \Fil
^{\{m_{ij}-1-(\ell-1)\}}D_{K, ij}$ from the previous paragraph, we deduce
that $N^{\ell-1}(x) \in \Fil^{\{m_{ij}-1-(\ell-1)\}}\D$.  Now (7)
  follows by (reverse) induction on $\ell$.

(8)  Again we proceed by induction on $m = \max_{ij}\{m_{ij}\}$.
Since $E_{i'j'}(u) x \in \fijd$, we have $f_{ij} ( E_{i'j'}(u) x) \in
\Fil ^{m_{ij}} D_{K, ij}$.  Then $f_{ij} (x) \in \Fil ^{m_{ij}} D_{K,
  ij}$ for any $(i,j) \neq (i',j')$ because $f_{ij} (E_{i'j'}(u)) \not
= 0$ in that case.
Note that \numequation\label{eq:ijeq} N(E_{i'j'}(u) x) = N(E_{i'j'}(u)) x + E_{i'j'}(u)N(x) \in
\Fil ^{\{m_{ij}-1\}}\D.\end{equation} Applying $f_{i'j'}$ to \eqref{eq:ijeq} and noting that $f_{i'j'}(N(E_{i'j'}(u) )) \not = 0$, we see that $ f_{i'j'} (x) \in \Fil ^{m_{i'j'}-1}D_{K, i'j'}$. Thus $f_{ij}(x)\in \Fil ^{m'_{ij}} D_{K, ij}$ for all $i, j$. It remains to show that
$N(x) \in \Fil ^{\{m'_{ij}-1\}}\D$. Note that $E_{i'j'}(u) x \in
\fijd$ implies that $E(u)x \in \fijd$, and (7) with $\ell=0$ shows
that that $ x \in \Fil^{\{m_{ij}-1\}}\D$. Hence $N( E_{i'j'}(u))x
\in\Fil^{\{m_{ij}-1\}}\D $ and then  \eqref{eq:ijeq} gives  $E_{i'j'} (u) N(x) \in
\Fil^{\{m_{ij}-1\}}\D$. By induction, we see that $N(x) \in \Fil
^{\{m'_{ij}-1\}}\D$, and we are done.
\end{proof}

\begin{rem}
  \label{rem:belowzero}
An argument analogous to the proof of Proposition~\ref{prop:
  properties of filmij}(1) shows that $\Fil^{\{m_{ij}\}}\D = \Fil^{\{\max(m_{ij},0)\}}\D$.
\end{rem}

The next proposition shows that the filtration $\Fil^{\{m_{ij}\}}\D$ is
essentially characterized by certain of the properties listed in Proposition~\ref{prop: properties of filmij}.

\begin{prop}
  \label{prop:uniqueness-of-fil}
Fix elements $n_{ij} \in \Z \cup \{\infty\}$ indexed by $i =
0,\ldots,f-1$ and $j = 0,\ldots,e-1$, and let $\cS$ be the set of collections
of integers $\{m_{ij}\}$  with $m_{ij} \le n_{ij}$ for all $i,j$.
Suppose that we have another
filtration of $\D$ by additive subsets $\widetilde \Fil
  ^{\{m_{ij}\}}\D$, defined for all
$\{m_{ij}\} \in \cS$ and
satisfying the following properties.
\begin{enumerate}
\item $\widetilde \Fil
  ^{\{m_{ij}\}}\D  = \D$ if $m_{ij} \le 0$ for all $i,j$.

\item For each $\{m_{ij}\} \in \cS$ there exists $r > 0$ such that
  $(\Fil^r S)\D \subseteq \widetilde \Fil^{\{m_{ij}\}} \D$.

\item The filtration $ \widetilde \Fil^{\{m_{ij}\}} \D$  satisfies
  properties (2) and (4) of
  Proposition~\ref{prop: properties of filmij}, as well as
  $E(u) \cdot \widetilde \Fil^{\{m_{ij}-1\}} \D \subseteq  \widetilde \Fil^{\{m_{ij}\}} \D$, for all $\{m_{ij}\} \in
  \cS$.
\end{enumerate}
Then for all $\{m_{ij}\} \in \cS$ we
have  $\widetilde \Fil ^{\{m_{ij}\}}\D =\Fil ^{\{m_{ij}\}}\D $.

The analogous statement holds for the filtration
$\Fil^{\{m_{i,0},\ldots,m_{i,e-1}\}} \D_i$ of \ref{para:filmijl},
 replacing $\D$ by $\D_i$.
\end{prop}

\begin{proof}  The proof for $\D_i$ is the same as the proof for
  $\D$, so again we concentrate on the latter case.  As usual we proceed by induction
  on $m=\max_{ij} \{m_{ij}\}$, with the base case $m = 0$ coming from (1).  Suppose
  that the claim holds for $m-1$, and consider the situation for $m$.

   Pick $x \in \widetilde \Fil ^{\{m_{ij}\}} \D$. By hypotheses \ref{prop: properties of filmij}(4) and
 \ref{prop: properties of filmij}(2) respectively, we have $f_{ij}(x)\in \Fil ^{m_{ij}} D_{K, ij}$ and $N(x) \in
 \widetilde \Fil ^{\{m_{ij}-1\}} \D$. By induction $\widetilde \Fil
 ^{\{m_{ij}-1\}} \D = \Fil ^{\{m_{ij}-1\}} \D $, so $x \in \fijd$ by definition, and
 we have shown that $\widetilde \Fil ^{\{m_{ij}\}} \D \subseteq \Fil
 ^{\{m_{ij}\}} \D$.

Conversely pick $x \in \fijd$.  By hypothesis \ref{prop: properties of filmij}(4) there exists $y \in \widetilde
\Fil^{\{m_{ij}\}}\D$ such that $f_{ij} (y) = f_{ij}(x)$ for all $i ,
j$.  We can (and do) modify $y$ so that $x= y + E(u)z$
with $z \in \D$ and $y$ still in $\widetilde \Fil^{\{m_{ij}\}}\D$.    This is possible because any $s \in \Fil^1 S$ can be
written as $E(u) s' + s''$ with $s' \in S[1/p]$ and $s'' \in \Fil^r S$
for $r \gg 0$; now use hypothesis (2).

Now $E(u)z=x-y \in \Fil^{\{m_{ij}\}}\D$ (from the second paragraph of the
proof), so $z \in  \Fil^{\{m_{ij}-1\}}\D$ by Proposition~\ref{prop:
  properties of filmij}(7).   Then $z \in
\widetilde\Fil^{\{m_{ij}-1\}}\D$ by induction, and hypothesis (3)
with $r=1$ shows that $E(u)z \in
\widetilde\Fil^{\{m_{ij}\}}\D$.  We conclude that $x  =  y + E(u)z \in
\widetilde\Fil^{\{m_{ij}\}}\D$, as desired.
 \end{proof}

\subsection{Some general facts about integral $p$-adic Hodge theory}
Suppose in this section that the Hodge--Tate weights of $V$ lie in the interval
$[0,r]$.  Let $T \subset V$ be a $G_K$-stable $\O_E$-lattice, and let
$\M$ be the
Kisin module attached to $T$ by the theory of
\cite{KisinCrys,MR2388556}; see the statements of \cite[Def.~3.1, Thm.~3.2, Def.~3.3, Prop.~3.4]{GLS-BDJ} for a concise
summary of the definitions and properties that we will need.
 The object $\M$ is a  finite free
$\fS_{\O_E}$-module with rank $d = \dim_{E} V$, together with an $\O_E$-linear
$\varphi$-semilinear map $\varphi : \M \to \M$ such that the cokernel
of $1 \otimes \varphi : \fS \otimes_{\varphi,\fS} \M \to \M$ is killed
by $E(u)^r$.   As in \cite[\S3]{GLS-BDJ} we use the contravariant functor
 $T_{\fS}$ to associate Galois representations to Kisin modules.

Set $\M^* = \fS \otimes_{\varphi, \fS}\M$, which by
\cite[Cor.~3.2.3]{MR2388556} (or \cite[Thm.~3.2(4)]{GLS-BDJ}) we can
view as a subset of $\D$.
Define $$\Fil ^{\{m_{ij}\}}\M^* := \M ^* \cap \fijd,$$
and  $M_{K, ij} = f_{ij}(\M^*) \subset D_{K, ij}$.  Similarly we
define $\M^*_i = \M^* \otimes_{\fS_E} \fS_{E,i}$ and
$\Fil^{\{m_{i,0},\ldots,m_{i,e-1}\}} \M^*_i := \M^*_i \cap \Fil^{\{m_{i,0},\ldots,m_{i,e-1}\}} \D_i.$. 
Note that $1 \otimes \varphi: \M^* \to \M$ is an  $\mathfrak S$-linear
map. By \cite[Cor.~3.2.3]{MR2388556} one also has  that 
\numequation
\label{eq: fil-Kisin}
\Fil ^r \M^* = \{ x\in \M^* | (1 \otimes \varphi) (x) \in E(u)^r \M\}.
\end{equation}

\begin{lemma}\label{lem: some integral properties} Let $\{m_{ij}\}$ be non-negative integers indexed by $i = 0,\ldots,e-1$ and $j=0,\ldots,f-1$.  Fix a pair
  $(i',j')$ and define  $m'_{ij} = m_{ij}+1 $ if $(i, j)= (i',j')$ and
  $m'_{ij} = m_{ij}$ otherwise.  Then we have the following.

\begin{enumerate}
\item $M_{K, ij}$ is an $\O_E$-lattice inside $D_{K, ij}$.
\item $\Fil ^{\{m_{ij}\}}\M^*/ \Fil ^{\{m'_{ij}\}}\M^*$ is a finite free $\O_E$-module.
\item If $\Fil^{m_{i'j'}+1}  D_{K,i'j'} = 0$ then $\Fil ^{\{m'_{ij}\}}\M^* = E_{i'j'} (u) \Fil ^{\{m_{ij}\}}\M^*$.
 \end{enumerate}
Moreover the natural analogues of (2) and (3) hold  for the filtration on $\M^*_i$.
 \end{lemma}
\begin{proof}
(1) Note that $M_K := f_\pi (\M^*) \subset D_K$ is a full $\O_K$-lattice. Indeed, $M_K \simeq \M ^*/ E(u)\M^*$ is a finite free $\O_K \otimes_{\Z_p} \O_E$-module of rank $d$ because $\M$ is finite $\fS_{\O_E}$-free of rank $d$. Then (1) follows quickly from this fact.

As usual the proofs of (2) and (3) will be the same for $\M^*_i$ as for
$\M^*$, and so we concentrate on the latter case.

(2) It is clear from the definition of $\Fil ^{\{m_{ij}\}}\M^*$ that
$\Fil ^{\{m_{ij}\}}\M^*/ \Fil ^{\{m'_{ij}\}}\M^*$ injects into $\Fil
^{\{m_{ij}\}}\D/ \Fil ^{\{m'_{ij}\}}\D$. By Proposition~\ref{prop:
  properties of filmij}(6),  $\Fil ^{\{m_{ij}\}}\D/ \Fil
^{\{m'_{ij}\}}\D$ is an $S_E /\Fil^1_{i'j'} S_E = E_{i'j'}$-module. In
particular it is $p$-torsion free, so the same is true of $\Fil
^{\{m_{ij}\}}\M^*/ \Fil ^{\{m'_{ij}\}}\M^*$. On the other hand  $\Fil
^{\{m_{ij}\}}\M^*/ \Fil ^{\{m'_{ij}\}}\M^*$ is an $\fS_{\O_E}/\Fil ^1
_{i'j'} \fS_{\O_E}= \O_{E_{i'j'}}$-module. Since it is $p$-torsion
free and finitely generated (note that $\M^*$ is finite
$\fS_{\O_E}$-free and $\fS_{\O_E}$ is noetherian), we see that $\Fil ^{\{m_{ij}\}}\M^*/ \Fil ^{\{m'_{ij}\}}\M^*$ is a finite free $\O_E$-module.

(3) Let $\cS$ be the
set of tuples $\{n_{ij}\}$ with $n_{ij}\le m'_{ij}$ for all $i,j$.
Fix any $r >
  \max_{i,j}\{m'_{ij}\}$,
define $$\widetilde\Fil^{\{m'_{ij}\}} \D = (\Fil ^1_{i'j'} S_E )\Fil
  ^{\{m_{ij}\}}\D + (\Fil^r S) \D,$$
and define $\widetilde \Fil^{\{n_{ij}\}} \D =
  \Fil^{\{n_{ij}\}}\D$ for all other tuples $\{n_{ij}\} \in \cS$.  (It
  is easy to check using  (1) and (5) of Proposition \ref{prop: properties of filmij}
  that the above definition does not depend on the choice of $r$.)
Our first goal is to show that  $\Fil ^{\{m'_{ij}\}}\D = \widetilde
\Fil ^{\{m'_{ij}\}}\D$.  This will follow from Proposition~\ref{prop:uniqueness-of-fil} once we can show that
\begin{enumerate}[(i)] \item $N(\widetilde \Fil ^{\{m '_{ij}\}}\D)
  \subseteq \Fil ^{\{m'_{ij}-1 \}}\D$,

\item  $f_\pi ( \widetilde \Fil ^{\{m '_{ij}\}}\D)  =  \bigoplus
  _{i,j} \Fil ^{m'_{ij}} D_{K, ij}$,

\item $E(u) \cdot \Fil ^{\{m '_{ij}-1\}}\D \subseteq
  \widetilde \Fil^{\{m'_{ij}\}} \D$.
\end{enumerate}

To check (i), let $s \in \Fil ^1_{i'j'} S_E$ and $x \in \fijd$. Then
$N(sx)= N(s) x + s N(x)$. We have $N(s) x \in \fijd \subseteq \Fil
^{\{m '_{ij}-1\}}\D$.  Similarly $N(x) \in \Fil ^{\{m_{ij}-1\}}\D$; by Proposition~\ref{prop: properties of filmij}(6) we have $s
N(x)$ and so also $N(sx)$ in $\Fil^{\{m '_{ij}-1\}}\D$.  Since
$N((\Fil^rS)\D) \subseteq N(\Fil^{\{r\}} \D) \subseteq \Fil^{\{m '_{ij}-1\}}\D$ as well, this checks (i).

Now let us check property (ii).
Certainly
$f_{\pi}$ does map  $(\Fil ^1_{i'j'} S_E )\Fil
  ^{\{m_{ij}\}}\D$ into  $ \bigoplus
  _{i,j} \Fil ^{m'_{ij}} D_{K, ij}$, and we have to check that this
  map is surjective.  (Since $f_{\pi}(\Fil^rS)=0$, the
$(\Fil^r S)\D$ term does not affect the statement.)
 Choose an element
$(x_{ij}) \in\oplus_{i,j} \Fil^{\{m'_{ij}\}} D_{K,ij}$.  Note that
$x_{i'j'}=0$ by our hypothesis that $\Fil^{m_{i'j'}+1} D_{K,i'j'} = 0$.
By Proposition~\ref{prop: properties of filmij}(4), there exists $y
\in \fijd$ such that $f_{ij}(y) = x_{ij}$ if $i \not = i'$ and
$f_{i'j}(y)=
\frac {1}{\pi_{i'j}- \pi_{i'j'}} x_{i'j}$ if $j \neq j'$. We
see that $f_{ij} (E_{i'j'}(u) y) = x_{ij}$ for all $i, j$, and this checks
property (ii).

Finally, property (iii) is an immediate consequence of
of the factorization $E(u) = \prod_{i,j} E_{ij}(u)$ together with
Proposition~\ref{prop: properties of filmij}(6) applied to
$E(u) \cdot \Fil^{\{m'_{ij}-1\}}\D$ repeatedly at all pairs \emph{other} than our
  fixed $(i',j')$.

Now turn to the statement we wish to prove.  It is clear from the
definitions and Proposition~\ref{prop: properties of filmij}(6) that $E_{i'j'} (u) \Fil ^{\{m_{ij}\}}\M^* \subseteq \Fil
^{\{m'_{ij}\}}\M^*$. Conversely, take $x \in \Fil ^{\{m'_{ij}\}}\M^*$.  By
the first part of the proof we can write
$x = s y + s'y' $ with $s \in \Fil ^1 _{i'j'} S_E$, $s' \in \Fil^r S$, $y
\in \fijd$, and $y' \in \D$. Since $\M^*$ is finite $\fS_{\O_E}$-free,
we can choose $e_1 , \dots ,e_d \in \M^*$ that are an
$\fS_{\O_E}$-basis of $\M^*$ (hence also an $S_E$-basis of $\D$). We
write $x = \sum _n a_n  e_n $ in terms of this basis, with $a_n \in
\fS_{\O_E}$ for all $n$.  The expression
$x = sy + s'y'$ shows that $a_n \in \Fil^1_{i'j'} S_E + \Fil^r S
=\Fil^1_{i'j'} S_E$ as well.  Hence $a_n  \in \Fil ^1_{i'j'} S_E \cap
\fS_{\O_E} = E_{i'j'}(u) \fS_{\O_E}$ by Lemma~\ref{lem: fil1ij}(3).  Now
$x = E_{i'j'}(u) y''$ with $y'' \in \M^*$. To conclude, we need to show
that $y'' \in \fijd$, and this follows from Proposition~\ref{prop: properties of filmij}(8).
 \end{proof}

As an application, we describe the Kisin modules
of crystalline characters.

\begin{lem}
  \label{lem:kisin-characters}
  Suppose that $V$ is a crystalline character such that $\HT_{\kappa_{ij}}(V)
  = \{r_{ij}\}$ with $r_{ij} \ge 0$ for all $i,j$.  Choose basis elements $e_i \in \M^*_i$
  for all $i$.  Then we have $\varphi(e_{i-1}) = \alpha_i
  \prod_{j=0}^{e-1} E_{ij}(u)^{r_{ij}} \cdot e_i$ with  $\alpha_i
  \in \O_E\llb u \rrb^{\times}$ for all $i$.
\end{lem}

\begin{proof}
 Choose $r \ge \max_{i,j} r_{ij}$.  We claim that
 $\Fil^r \M^*_i = \prod_{j=0}^{e-1} E_{ij}(u)^{r-r_{ij}} \cdot \M^*_i$ for all~$i$.
 It is immediate from the definition that
 $\Fil^{\{r_{ij}\}} \D= \D$, so that
 $\Fil^{\{r_{ij}\}} \M^* = \M^*$ as well.  Now the claim follows by
 repeatedly applying Lemma~\ref{lem: some integral properties}(3).

 Recall that we have $\M^* = \fS \otimes_{\varphi,\fS} \M$, so that
 $e_{i-1}$ (and not $e_i$) is a generator of $\M^*_i$.
Also recall from \cite[Lem.~4.3(1)]{GLS-BDJ} that we have $\Fil^r \M^*_i =
 \{x \in \M^*_i \, : \, (1 \otimes \varphi)(x) \in E^{\kappa_i} (u)^r
 \, \M_i\}$.
 Writing $\varphi(e_{i-1}) = \alpha'_i e_i$, we see that $\Fil^r \M_i^*
 = \{ \beta e_{i-1} \, : \, E^{\kappa_i} (u)^r \mid \beta \alpha'
 \}$.  From the shape of $\Fil^r \M^*$ we deduce that $\alpha'_i =
 \alpha_i \prod_{j=0}^{e-1} E_{ij}(u)^{r_{ij}}$ for some unit
 $\alpha_i$, as desired.
\end{proof}

 \subsection{Pseudo-Barsotti--Tate representations}
\label{sec:pseudo-BT}

 \begin{defn}
   \label{defn:pseudo-bt}
  Fix  integers $r_i \in [1,p]$ for all $i$.
  We say that a two-dimensional crystalline $E$-representation $V$ of
  $G_K$  is \emph{pseudo-Barsotti--Tate} (or pseudo-BT) of weight $\{r_i\}$
 if for all $0 \le i \le f-1$
  we have $\HT_{\kappa _{i,0}} ( V) =\{ 0 , r_i\}$, and $\HT_{\kappa_{ij}} (V)=\{0, 1 \}$ if $j \neq 0$.  (Strictly speaking, we should say that $V$ is pseudo-BT with respect
to the labelling $\kappa_{ij}$ of the embeddings $\Hom_{\Qp}(K,E)$.)
 \end{defn}

Equivalently,  a two-dimensional crystalline representation $V$ is
pseudo-BT if and only if
\begin{itemize}
\item $\dim _E \gr^0D_{K, i j}=1 $ for all $i, j$,
\item $\dim _E \gr ^{r_i} D_{K, i0} =1$ with $r_i \in
  [1,p]$ for all $i$,  and
\item $\dim _E \gr ^{1} D_{K, ij} =1$ for all $j \not =0 .$
\end{itemize}

Note that a pseudo-BT representation is actually Barsotti--Tate if and
only if $r_i = 1$ for all $i$.  For the remainder of this section, we assume that $V$ is pseudo-BT.

\para As in the previous subsection, we let $T \subset V$ be an $G_K$-stable
$\O_E$-lattice, and $\M$ the Kisin module attached to $T$.  Write $M_{K, ij} = f_{ij}(\M^*) \subset D_{K, ij}$. Note that $\M^*$ is a rank-two finite free $\fS_{\O_E}$-module.
 Set $\Fil ^{m} M_{K, ij} : = M_{K, ij} \cap \Fil ^m D_{K,ij}.$ By
 definition,
 \begin{itemize}
 \item $\Fil ^1 M_{K, ij}$ is a saturated rank-one free $\O_E$-submodule
   of $M_{K, ij}$,
 \item  $\Fil ^{r_i} M_{K, i0} = \Fil ^1 M_{K, i0}$ and $\Fil^{r_i +1}
   M_{K, i0}= \{0\}$, and
\item $\Fil ^2 M_{K, ij} = \{0\}$ if $j \not =0$.
 \end{itemize}

In the next few pages, we will establish several results about the
structure of the submodules $\Fil^{\{n,0,\ldots,0\}} \M_i^*$ of
$\M^*_i$, following roughly the same strategy as in
\cite[\S4]{GLS-BDJ}.  We remark that until Corollary~\ref{cor: shape of filp}, none
of these results will actually use the pseudo-BT hypothesis, only that
the crystalline representation $V$ has non-negative Hodge--Tate weights and $\HT_{\kappa_{i,0}}(V) =
\{0,r_i\}$.

We begin with following proposition, which should be compared with \cite[Prop.~4.5]{GLS-BDJ}.

\begin{prop}\label{prop: shape of 1st level} Suppose that there exists an $\fS_{\O_{E},i}$-basis $\{\e _i , \f_i\}$ of $\M^*_i$ such that
$\f_i \in \Fil ^{\{r_i , 0 \dots, 0\}}\M_i^*$ and
$f_{i,0} (\f_i )$ generates $\Fil ^{r_i} M_{K, i0}$. Then for any
$n\ge r_i$ we have
 $$\Fil^{\{n,0,\dots, 0 \}}\M^*_i = \fS_{\O_E,i} (u - \pi_{i,0}) ^{n} \e_i \oplus  \fS_{\O_E,i} (u - \pi_{i,0}) ^{n-r_i}\f _i .$$
\end{prop}

\begin{proof}
We first prove by induction that for $0 \le n \le r_i$ we have
\numequation \label{eq: shape of fil(1)} \Fil ^{\{n, 0, \dots , 0\}}\D_i = S_{E,{i}}(u - \pi_{i,0})^n \e_i \oplus S_{E,i} \f_i + (\Fil^n S_{E,i}) \D_i . \end{equation}
The case that $n = 0$ is trivial as $\{\e_i , \f_i\}$ is a basis of
$\M^*_i$. Suppose that the statement is valid for $n-1$, and consider
the statement for $n$.  Define $\widetilde \Fil ^{\{n, 0, \dots ,
  0\}}\D_i$ to be equal to the right-hand side of \eqref{eq: shape of
  fil(1)}, and set $\widetilde \Fil ^{\{m, 0, \dots ,
  0\}}\D_i = \Fil ^{\{m, 0, \dots ,  0\}}\D_i$ for $m < n$.  The
conditions of Proposition~\ref{prop:uniqueness-of-fil} for this
filtration are straightforward to check:\ use $\f_i \in \Fil^{\{r_i,0,\ldots,0\}} \D_i$ to see that $N(\f_i)
\in \widetilde \Fil^{\{n-1,0,\ldots,0\}} \D_i =
\Fil^{\{n-1,0,\ldots,0\}} \D_i $, so property (2) of \ref{prop:
  properties of filmij} holds; then to verify property (4) of \ref{prop:
  properties of filmij},  note that for  $1\le n \le
r_i$, $\Fil^nD_{K, i0}=\Fil^{r_i}D_{K, i0}$ is generated by
$f_{i,0}(\f_i)$ by hypothesis.
By
Proposition~\ref{prop:uniqueness-of-fil}  we deduce that $\widetilde\Fil ^{\{n, 0, \dots ,
  0\}}\D_i = \Fil ^{\{n, 0, \dots , 0\}}\D_i$, as desired.

Now by Lemma~\ref{lem: some integral properties} for $\D_i$, it suffices to show that
$$\Fil ^{\{r_i, 0, \dots , 0\}}\M^*_i = \fS_{\O_{E},i}(u - \pi_{i,0})^{r_i} \e_i \oplus \fS_{\O_E,i} \f_i . $$
Evidently $ \fS_{\O_{E},i}(u - \pi_{i,0})^{r_i} \e_i \oplus
\fS_{\O_E,i} \f_i  \subseteq \Fil ^{\{r_i, 0, \dots , 0\}}\M^*_i$. For
the converse, let $x \in \Fil ^{\{r_i, 0, \dots , 0\}}\M^*_i$. Then by
\eqref{eq: shape of fil(1)} we have $x = (s (u - \pi_{i,0})^{r_i}
+ t) \e_i + (s'+t')\f_i$ with $s, s' \in S_{E,i}$ and $t,t'\in \Fil ^{r_i} S_{E,i}$. Since $x \in \M_i ^*$, we see that
$s (u - \pi_{i,0})^{r_i} +t\in \fS_{\O_E}$ and $s' + t' \in
\fS_{\O_E}$. Note that $s (u - \pi_{i,0})^{r_i} +t$ is in $\Fil ^{r_i}
S_{E,i}\cap \fS _{\O_E}$, and since this is just $(u -
\pi_{i,0})^{r_i}\fS_{\O_E}$, we are done.
\end{proof}

 Via the identity $\M^* = \fS \otimes_{\varphi,\fS} \M$, we
can regard $\M_{i-1}$ as a $\varphi(\fS)$-submodule of $\M^*_i$ such that $\fS
  \otimes_{\varphi(\fS)} \M_{i-1} = \M^*_i$. (Note that
  the analogous statement immediately
  preceding the published version of \cite[Thm.~4.22]{GLS-BDJ} contains a mistake.  See
  Appendix~\ref{sec:corrigendum} for a correction.)
The following proposition is our key result on the structure of $\Fil
^{\{r_i , 0 , \dots, 0\}} \M _i^*$, and may be compared with \cite[Prop.~4.16]{GLS-BDJ}.

\begin{prop}\label{prop: Key shape}  \hfill

\begin{enumerate}
\item There exists a basis $\{\e_{i-1}, \f_{i-1} \}$ of $\M_{i-1} $ such that $f_{i,0} (\f_{i-1})$ generates $\Fil ^{r_i}M_{K, i0}$.
\item Suppose that $p \ge 3$.  There exists a basis  $\{\e'_{i}, \f'_{i} \}$  of $\M^*_i$ such
  that $\e'_{i} -\e_{i-1}$ and $\f'_i - \f_{i-1}$ are in $\m_E \M^*_i$ and
$\Fil^{\{ p, 0, \dots, 0\}}\M^*_i = \fS_{\O_E,i} (u - \pi_{i,0}) ^{p} \e'_i \oplus  \fS_{\O_E,i} (u - \pi _{i,0})^{p-r_i}\f' _i. $
\end{enumerate}
\end{prop}

\begin{proof}
 (1) We know by Lemma~\ref{lem: some integral properties}(1) that
 $f_{i,0}(\M^*_i)$ is an $\O_{E_i}$-lattice inside $D_{K, i0}$, so there exists an $\O_{E_i}$-basis $\{\bar \e , \bar \f\}$ of $M_{K, i0}$ such that $\bar \f$ is a basis of $\Fil ^{r_i} D_{K, i0}$. Pick any $\varphi(\fS_{\O_E,i})$-basis $\{\widetilde \e_{i-1}, \widetilde \f_{i-1}\}$ of $\M_{i-1}$. Then $\{f_{i,0}(\widetilde \e_{i-1}), f_{i,0} (\widetilde \f_{i-1})\}$ forms a basis of $f_{i,0} (\M^*_i)$. Let $A \in \GL_2(\O_E)$ be the matrix such that
$(f_{i,0}(\widetilde \e_{i-1}), f_{i,0} (\widetilde \f_{i-1})) A = (\bar \e , \bar \f)$. Then $(\e_{i-1}, \f_{i-1}) = (\widetilde \e_{i-1}, \widetilde \f_{i-1}) A $ is the desired basis.

(2) The proof is a minor variation of the proof of
\cite[Prop.~4.16]{GLS-BDJ}.
Let $\e := \e_{i-1}$ and $\f:=
\f_{i-1}$ be a basis of $\M_{i-1}$ as in part (1).  In this proof
only, we will write $\pi:= \pi_{i,0}$ and $r := r_i$, to avoid too many
subscripts and  to make the discussion easier to compare with the
proof of \cite[Prop.~4.16]{GLS-BDJ}.   We consider the following assertion:

($\star$) For each $n = 1, \ldots, r$ there exists $\f ^{(n)} \in \Fil ^{\{n , 0 , \dots, 0\}}\M^*_i$ such that
$$\f ^{(n)} = \f + \sum _{s= 1 }^{n-1} \pi^{p-s}(u -\pi)^s(a_s^{(n)}
\e + \widetilde a _s ^{(n)} \f)  $$ with $a^{(n)}_s , \widetilde
a^{(n)}_s \in \O_{E_i }.$   Once the assertion ($\star$) is
established, the
proposition follows from Proposition~\ref{prop: shape of 1st level}, taking $\e'_{i} = \e$ and $\f' _{i} = \f ^{(r)}$.

We prove $(\star)$ by
induction on $n$.  From the definition of $\Fil ^{\{1, 0 \dots
  , 0\}}\D_i$ and the defining property of $\f$ from (1), we see that $\f  \in \Fil ^{\{1
  , 0 , \dots, 0\}}\M^*_i$, so that for the base case $n=~1$ we can take $\f ^{(1)} = \f$.

Now assume that $(\star)$ is valid for some $1 \le n < r$, and let us consider the case $n + 1$.
Set $H(u) = \frac{u-\pi}{\pi} $ and $$\widetilde \f^{(n+1)}: = \sum _{\ell  =0} ^ n \frac{H(u)^\ell N^\ell (\f^{(n)})}{\ell !}. $$
As in the proof of \cite[Prop.~4.16]{GLS-BDJ} (and the proof of Proposition~\ref{prop: properties of
    filmij}(4) of this paper) one computes that
$$ N(\f^{(n+1)}) = \frac{H(u)^n}{n!} N^{n+1}(\f^{(n)}) + \sum_{\ell=1}^{n}
\frac{(1+N(H(u)))H(u)^{\ell-1} N^{\ell}(\f^{(n)})}{(\ell-1)!}.$$
Then, using $1+N(H(u)) \in \Fil^1_{i,0} S_{E,i}$, one deduces that $N(\widetilde \f^{(n+1)}) \in
  \Fil^{\{n,0,\ldots,0\}} \D$ and so $\widetilde \f^{(n+1)} \in \Fil ^{\{n+1 , 0, \dots , 0\}}\D$. Now by induction, we have
\nummultline\label{eq:bigeq}
\widetilde \f^{(n+1)} -  \f^{(n)} =   \sum _{\ell  =1} ^ n \frac{ (u - \pi)^\ell} { \pi ^\ell \ell!} N^\ell\left(\f + \sum _{s= 1 }^{n-1} \pi^{p-s}(u -\pi)^s(a_s^{(n)}  \e + \widetilde a _s ^{(n)} \f) \right )\\
= \sum _{\ell  =1} ^ n \frac{ (u - \pi)^\ell} { \pi ^\ell \ell!}
N^\ell (\f) + \hfill \\
  \sum_{\ell  =1}^n \sum_{s=1}^{n-1} \sum_{t=0}^{\ell} \frac{ \pi ^{p-s} (u -
   \pi)^\ell}{\pi ^\ell \ell!} \binom{\ell}{t} N^{\ell -t} ((u -\pi)^s) \left  (a _s^{(n)}  N^{t}(\e) + \widetilde a_s^{(n)}  N^{t}(\f) \right ).
\end{multline}
Write
\numequation\label{eq:subst}
N^t(\e) = \sum \limits_{m= 0}^\infty  (u-\pi)^m (c^t_m \e + d^t_m \f)
\quad \text{and} \quad N^t(\f) = \sum \limits_{m= 0}^\infty  (u-\pi)^m
(\widetilde c^t_m \e + \widetilde d^t_m \f)
\end{equation} with $c^t_m , d^t_m , \widetilde c^t_m , \widetilde d^t_m \in E$.
After substituting \eqref{eq:subst} into \eqref{eq:bigeq} and
expanding the terms $N^{\ell-t}((u-\pi)^s)$ using
\cite[Lem.~4.13]{GLS-BDJ}, we can collect terms and write
\numequation\label{eq:collected} \widetilde \f^{(n+1)} =   \f^{(n)} + \sum_{m = 1}^\infty (u -\pi)^m(b_m \e + \widetilde b _m \f )\end{equation}
with each $b_m, \widetilde b_m \in E$.  Now we delete all terms of $(u- \pi )$-degree at least
$n+ 1$ from this expression, and define
$$\f^{(n+1)}:  =   \f^{(n)} + \sum_{m = 1}^n  (u -\pi)^m(b_m \e + \widetilde b _m \f ). $$
It remains to show that $\pi ^{p-m} \mid b_m,\widetilde b_m$, or in
other words that if $m \le n$ then every occurrence of $(u-\pi)^m$ in
the terms collected to form \eqref{eq:collected}  has coefficient divisible by $\pi^{p-m}$.
A direct examination of these terms (just as in the last part of the
proof of \cite[Prop.~4.16]{GLS-BDJ}) shows that this comes down to the
claim that $\pi^{p-m} \mid c^\ell_m , d^\ell_m , \widetilde c^\ell_m ,
\widetilde d^\ell_m$ for all $1 \le \ell < p$ and $0 \le m < p$.  By
\cite[Cor.~4.11]{GLS-BDJ}, this follows from
Lemma~\ref{lem: newtechnical2} below (applied at $a_i = c^\ell_i$,
$d^\ell_i$, etc.), which generalizes
\cite[Lem.~4.12]{GLS-BDJ}.  (We remind the reader  that except for \cite[Lem.~4.12]{GLS-BDJ}, in the
results of \cite[\S4.2]{GLS-BDJ}
our ground field was an arbitrary finite extension of $\Qp$.)  The
application of \cite[Cor.~4.11]{GLS-BDJ} is where we use the
hypothesis that $p\ge 3$.
\end{proof}

Define $S' = W(k)\llb u^p,\frac{u^{ep}}{p}\rrb\left[\frac{1}{p}\right]
\cap S$ and set $\cI_\ell = \sum_{m=1}^{\ell} p^{\ell-m} u^{pm} S'
\subset S'$.

\begin{lemma}\label{lem: newtechnical2} Suppose that $y \in \mathcal I_\ell$ for
  some $1 \le \ell \le p$.
Write $y   = \sum\limits_{i=0} ^\infty a_{i} (u -\pi)^i $ with $a_{i}\in K_0$.
Then we
have
$\pi ^{p+ (\ell-1)\min(p,e) } \mid a_{0}$ and $\pi ^{p+e-i+ (\ell-1)\min(p,e) } \mid  a_{i}$ in
$\O_K$ for $1 \leq
i \leq p-1$.
\end{lemma}
\begin{proof} For any non-negative integer $n$, let $e(n) = \lfloor\frac n e\rfloor$. Note that any $x \in S$ can be written uniquely as $x = \sum \limits_{i =0}^\infty a_i \frac{u ^i}{e(i)!}$ with $a_i \in W(k)$.

 By hypothesis we have $y =  \sum\limits_{m =1}^\ell p
  ^{\ell -m} u ^{pm}z_{m} $ with $z_{m } \in S'$. We can write $z_m  =
  \sum\limits_{j =0}^\infty \frac{b_{j, m}u ^{p j}}{e(pj)!}$ with
  $b_{j, m }\in W(k)$.  Then
\begin{eqnarray*}
 y &= & \sum_{m=1}^\ell p ^{\ell -m} \left(\sum\limits_{j =0}^\infty b_{j, m  } \frac{u ^{p (j+m) }}{e(pj)!}\right) \\
&= & \sum_{m=1}^\ell\sum\limits_{j =0}^\infty p ^{\ell -m} b _{j, m }\frac{(u-\pi +\pi)  ^{p (j+m)}}{e(pj)!} \\
&= & \sum_{m=1}^\ell\sum_{j =0}^\infty p ^{\ell -m}\frac{b_{j, m }}{e(pj)!} \left ( \sum_{i=0}^{p(j+m) }  {\binom{p(j+m)}{i}} (u-\pi)^i \pi ^{p (j+m)-i}  \right )\\
&= & \sum _{i =0}^\infty  \left (\sum _{m=1}^\ell \sum_{j \geq s_{i , p ,
      m }}\frac{b_{j, m } \pi^{p(j+m )-i} p^{\ell  -m }}{e(pj)!} {\binom{p(j+m)}{i}} \right )  (u-\pi)^i,
\end{eqnarray*}
where $s _{i, p , m}= \max \{0, i/p-m\}$. Since we only consider $a_i$
for $0 \leq i \leq  p$, we have $s_{i, p, m }= 0 $ in all cases.
Note that $\pi^{pj}/e(pj)! \in \O_K$ for all $j \ge 0$. Let $v_{\pi}$
denote the valuation on $\O_K$ such that $v_{\pi}(\pi)=1$.
We first observe that  $v_\pi (a_0) \geq \min_{1 \le m \le \ell}
(pm+ e(\ell-m)) = p + (\ell-1)\min(p,e)$. If $1 \leq i \leq p-1$ then
$p$ divides  ${\binom{p(j+m)}{i}}$, so we get $$v_\pi (a_i)\geq
\min_{1 \le m \le \ell} (pm-i + e(\ell -m) +e) = p + e - i +
(\ell-1)\min(p,e)$$ instead.
\end{proof}

In what follows, when we write a product of matrices
as $\prod_{j=1}^{n} A_i$, we mean $A_1 A_2 \cdots A_n$.

\begin{cor} \label{cor: shape of filp}
Suppose that $p \ge 3$
and let $\M$ be the Kisin module corresponding
to a lattice in a pseudo-BT
representation $V$ of weight $\{r_i\}$. There exist matrices $Z'_{ij} \in \GL_2 (\O_E)$ for
$j = 1, \dots,  e-1$ such that $\Fil ^{p, p \dots, p}\M^*_i =
\fS_{\cO_E,i} \alpha_{i,e-1} \oplus \fS_{\cO_E,i} \beta_{i,e-1}$ with $$
(\alpha_{i,e-1}, \beta_{i,e-1})= (\e'_i,  \f'_i) \Lambda'_{i,0}
\left(\prod_{j =1}^{e-1} Z'_{ij} \Lambda'_{ij} \right)$$ where
$\e'_i,\f'_i$ are as in Proposition~\ref{prop: Key shape}(2), $\Lambda'_{i,0} = \begin {pmatrix} (u - \pi _{i,0})^p & 0 \\ 0 & (u - \pi_{i,0})^{p-r_i} \end{pmatrix}$ and $\Lambda'_ {ij} = \begin {pmatrix} (u - \pi _{ij})^p & 0 \\ 0 & (u - \pi_{i j})^{p-1} \end{pmatrix}$ for $j =1 , \dots, e-1$.
\end{cor}

\begin{proof} Define $\mf p_m := \{p, \dots , p, 0, \dots, 0\}$ where
  the tuple contains exactly $m+1$ copies of $p$.  We prove by
  induction on $m$ that there exist matrices $Z'_{ij} \in \GL_2 (\O_E)$
  for $j = 1,\ldots,m$ such that
$$\Fil ^{\mf p_m }\M^*_i = \fS_{\O_E,i} \alpha_{i,m}  \oplus
\fS_{\O_E,i} \beta_{i,m}$$ with
$(\alpha_{i,m}, \beta_{i,m}): = (\e'_i , \f'_i) \Lambda'_{i,0}
\left(\prod_{j =1}^{m} Z'_{ij} \Lambda'_{ij} \right)$.
If $m =0$ then this is Proposition~\ref{prop: Key shape}(2).

Suppose the statement holds for $m-1$,  and let us consider the
statement for $m$. We first show that $ M' _{K, im}:= f_{im}(\Fil
^{\mf p_{m-1}}\M^*_i)$ is an $\O_{E_{im}}$-lattice inside $ D_{K ,
  im}$, or equivalently that $\{ f_{im} (\alpha_{i,m-1}), f_{im} (\beta_{i ,m-1})\}$ is an $E_{im}$-basis of $D_{K, im}$. Since $\{\e'_i , \f'_i\}$ is a basis of $\M^*_i$,  it suffices to check that $f_{im} \left(\Lambda'_{i,0}
\left(\prod_{j =1}^{m-1} Z'_{ij} \Lambda'_{ij} \right)\right)$ is an
invertible  matrix in $\GL_2 (E) $. This holds because $Z'_{ij} \in
\GL_2(\O_E)$ and $ f_{im} (\Lambda'_{ij}) \in \GL_2 (E)$ for $j <m$.

 The fact that $ M' _{K, im}$ is an $\O_{E_{im}}$-lattice inside $
 D_{K , im}$ implies that  there exists an $\O_{E_i}$-basis $\bar
 \gamma' , \bar \delta'$ of $M'_{K, im}$ such that $\bar\delta'$
 generates $\Fil ^1 D_{K, im}$. Write  $\bar \alpha_{m-1} : = f_{im}
 (\alpha_{i ,m-1})$ and $\bar \beta_{m-1} := f_{im}(\beta_{i
   ,m-1})$. Let $Z'_{im} \in \GL_2 (\O_E)$ be the matrix such that
 $ (\bar \gamma ',\bar
 \delta' )   = (\bar \alpha_{i,m-1}, \bar \beta_{i,m-1}) Z'_{im} $ and define  $(\gamma _m,\delta_m):  = (\alpha_{i,m-1},
 \beta_{i,m-1}) Z'_{im}$. Let $\mf q_{m} := \{p, \dots , p, 1 , 0
 \dots, 0 \}$ where the tuple contains exactly $m$ copies of $p$. We claim
\numequation\label{eq:filqm} \Fil ^{\mf q_m} \M^*_i   =   \fS_{\O_E,i} (u - \pi_{im})\gamma_m \oplus  \fS_{\O_{E,{i}}}\delta_m . \end{equation}

We first show that $(u
-\pi_{im})\gamma_m, \delta_m$ are in $\Fil ^{\mf q_m}\M_i^*$. Note that
$\gamma_m, \delta _m$ generate $\Fil^{{\mf p}_{m-1}} \M^*$ by
construction, so by Proposition~\ref{prop: properties of filmij}(6)
for $\D_i$ it suffices to show that $\delta_m \in \Fil ^{\mf
  q_m}\D_i$. Note that $f_{im}(\delta_m)= \bar \delta' \in \Fil ^1
D_{K, im}$, so we just need to check that $N(\delta_m) \in \Fil
^{\{p-1 , \dots, p-1 , 0 , 0, \dots, 0\}}\D_i$, where there are
$m$ copies of $p-1$ in the superscript
(\emph{cf.}~Remark~\ref{rem:belowzero}). But this follows from the
fact that $\alpha_{i,m-1}$ and $\beta_{i,m -1}$ are in $\Fil ^{\mf
  p_{m-1} }\D_i$. Therefore, $\fS_{\O_{E,{i}}}(u -\pi_{im})\gamma_m
\oplus \fS_{\O_E,i} \delta_m \subseteq \Fil ^{\mf q_m}\M^*_i.$

 Now pick $x \in \Fil ^{\mf  q_m}\M^*_i\subseteq \Fil ^{{\mf
     p}_{m-1}}\M^*_i$. We have $x = a \gamma_m + b \delta_m$ with $a,
 b \in \fS_{\O_E,i}$. It suffices to show that $(u -\pi_{im}) \mid
 a$. Note that $f_{im}(x) = f_{im}(a) f_{im} (\gamma_m) + f_{im}(b)
 f_{im} (\delta_m )\in \Fil ^1 D_{K, im}$. But $(f_{im}(\gamma_m),
 f_{im}(\delta _m))$ is just $(\bar \gamma' , \bar \delta')$,
 which is  a basis of $D_{K, im}$, and $\bar \delta'$ generates
 $\Fil ^1 D_{K, im}$. This forces $f_{im}(a) = 0$ and then $(u -
 \pi_{im}) \mid a$ by Lemma~\ref{lem: fil1ij}.  This completes the proof of \eqref{eq:filqm}.

Finally, recall that $\Fil ^2 D_{K, im}= \{0\}$ since $V$ is pseudo-BT, so
the equality \eqref{eq:filqm}
together with Lemma~\ref{lem: some integral properties}(2) implies that  $$\Fil ^{\mf p_m }\M_i^*  =  \fS_{\O_E,i} (u - \pi_{im})^p \gamma_m \oplus (u - \pi_{im}) ^{p-1}\fS_{\O_E,i} \delta_m. $$
That is, $\Fil ^{\mf p_m}\M_i^*$ is generated by  $(\e'_i , \f'_i) \Lambda'_{i,0}
\left(\prod_{j =1}^{m} Z'_{ij} \Lambda'_{ij} \right)$. This completes the induction on $m$ and proves the proposition.
\end{proof}

\subsection{The structure theorem for pseudo-Barsotti--Tate Kisin
  modules}

In this subsection we prove our main result about
Kisin modules associated to pseudo-BT representations.  We retain the
notation from the previous subsection.

\begin{thm}\label{thm: structure}
Suppose that $p \ge 3$ and let $\M$ be the Kisin module corresponding
to a lattice in a pseudo-BT
representation $V$ of weight $\{r_i\}$.  Then there exists an $\O_E \llb u\rrb $-basis
$\{\e_i, \f_i \}$ of $\M_i $ for all $0 \le i \le f-1$ such
that
$$\varphi(\e_{i-1}, \f_{i-1}) = (\e_{i}, \f_{i})X_i \left(\prod_{j=1}^{e-1} \Lambda_{i,e-j} Z_{i,e-j} \right)
  \Lambda_{i,0} Y_i ,  $$
for $X_{i}, Y_i \in \GL_2(\O_E\llb u\rrb )$
with $Y_i \equiv I_2 \pmod{\m_E}$,  matrices $Z_{ij} \in \GL_2(\O_E)$ for all~$j$,
and $\Lambda_{i,0} = \begin {pmatrix} 1 & 0 \\ 0 & (u - \pi_{i,0})^{r_i} \end{pmatrix}$ and $\Lambda _ {ij} = \begin {pmatrix} 1 & 0 \\ 0 & u - \pi_{ij} \end{pmatrix}$ for $j =1, \dots, e-1$.

\end{thm}
\begin{proof} For all $i$ we let $\{\e_i, \f_i\}$  be the $\O_{E}\llb
  u\rrb $-basis of $\M_i$ in Proposition~\ref{prop: Key shape}(1) and
  write $\varphi (\e_{i-1}, \f_{i-1}) = (\e_{i}, \f_{i}) A_i$, where
  $A_i$ is a matrix with coefficients in $\O_{E}\llb u\rrb $. From \eqref{eq: fil-Kisin} we see that 
  $\Fil ^p \M^*_i$ is generated by $(\e_{i-1}, \f_{i-1}) B_i$, where
  $B_i$ is the matrix satisfying $A_i B_i = (E^{\kappa_i}(u))^p I
  _2$. On the other hand, Corollary~\ref{cor: shape of filp} shows
  that $\Fil ^p \M_i^*$ (which by definition is equal to  $\Fil ^{\{p, \dots ,
    p\}}\M^*_i$) is generated by $$(\alpha_{i,e-1},
  \beta_{i,e-1} )= (\e'_i , \f'_i ) \Lambda'_{i,0} \left (
    \prod_{j=1}^{e-1}  Z'_{ij} \Lambda'_{ij} \right) = (\e_{i-1} , \f_{i-1}) Y^{-1}_i  \Lambda'_{i,0} \left (
    \prod_{j=1}^{e-1}  Z'_{ij} \Lambda'_{ij} \right).$$ Here $Y_i$ is
  the matrix such that $(\e_i' , \f'_i)Y_i = (\e_{i-1} , \f_{i-1})$,
  which by
  Proposition~\ref{prop: Key shape}(2) is congruent to the identity
  modulo $\m_E$. Therefore there exists an invertible matrix $X_i \in \GL_2( \O_E\llb u\rrb )$ such that $(\e_{i-1}, \f_{i-1}) B_i  = (\alpha_{i,e-1}, \beta_{i,e-1})X_i^{-1}$. Hence we have
$$Y^{-1}_i  \Lambda'_{i,0} \left (
    \prod_{j=1}^{e-1}  Z'_{ij} \Lambda'_{ij} \right) X_i^{-1} = B_i.$$
  Then the relation $A_i B_i = (E^{\kappa_i}(u))^p I _2$ proves
  that $$A_i = X_i \left(\prod_{j=1}^{e-1} \Lambda_{i,e-j} Z_{i,e-j} \right)
  \Lambda_{i,0} Y_i $$ with $Z_{ij} = (Z'_{ij})^{-1}$ and
  $\Lambda_{ij} = E_{ij}(u)^p (\Lambda'_{ij})^{-1}$.
\end{proof}

\section{Semisimple reductions mod $p$ of pseudo-BT representations}
\label{sec:reductions-mod-p}

\subsection{\hskip 0 cm} \hskip -0.2 cm
In this  section we use Theorem~\ref{thm: structure} to
study the semisimple representations that can be obtained as the
reduction modulo~$p$ of pseudo-BT representations.
We begin with the
following notation.

\begin{defn}
  \label{defn:rank-one}
   Suppose $s_0,\ldots,s_{f-1}$ are non-negative integers and $a \in
   k_E^{\times}$.  Let $\barM(s_0,\ldots,s_{f-1}; a)$ be the Kisin
   module with natural $k_E$-action (in the sense of~\cite[\S 3]{GLS-BDJ}) that has rank one
   over $\fS_{\O_E} \otimes_{\O_E}
  k_E$ and satisfies
  \begin{itemize}
  \item $\barM(s_0,\ldots,s_{f-1}; a)_i$ is generated by $e_i$, and

\item $\varphi(e_{i-1}) = (a)_{i} u^{s_{i}} e_{i}$.
  \end{itemize}
Here $(a)_i = a$ if $i \equiv 0 \pmod{f}$ and $(a)_i = 1$ otherwise.
All Kisin modules of rank one have this form (see e.g.\ \cite[Lem.~6.2]{GLS-BDJ}).

Write $\kappabar_i$ for the embedding $k \to \Fpbar$ induced by
$\kappa_{ij}$ (this is independent of $j$).
For brevity we will sometimes write $\omega_i$ for the fundamental character $\omega_{\kappabar_i}$.

We refer the reader to \cite[\S3]{GLS-BDJ} for the definition of the
contravariant functors $T_{\fS}$ that associate a representation of
$G_{K_{\infty}}$ to each torsion or finite free Kisin module.
\end{defn}

\begin{lem}
  \label{lem:functor-on-characters}  We have $
  T_{\fS}(\barM(s_0,\ldots,s_{f-1}; a)) \simeq
  \chibar|_{G_{K_{\infty}}}$ for a unique character $\chibar : G_K \to k_E^{\times}$, and
  $\chibar$ satisfies
$ \chibar|_{I_K} \simeq \prod_{i=0}^{f-1} \omega_i^{s_i}.$
\end{lem}

\begin{proof}
Choose any integers $r_{ij} \ge 0$ such that $\sum_j r_{ij} = s_i$.
By Lemma~\ref{lem:kisin-characters} (together with
\cite[Lem.~6.4]{GLS-BDJ} and an analysis of the Kisin modules
associated to unramified characters as in  the proof of
\cite[Lem.~6.3]{GLS-BDJ}) we see that $\barM(s_0,\ldots,s_{f-1}; a)$
is isomorphic to $\M\otimes_{\O_E} k_E$ for a Kisin module $\M$ corresponding
to a lattice $T$ in a crystalline character $V$  with Hodge--Tate weights
$\HT_{\kappa_{ij}}(V) = \{r_{ij}\}$.  Then $\chibar = T \otimes_{\O_E}
k_E$.  The character $\chibar$ is unique since
$K_{\infty}/K$ is totally wildly ramified, so that
restriction to $G_{K_{\infty}}$ is faithful on characters of $G_K$.

For the last part of the statement it suffices to check that
$\overline{\psi}_{ij} |_{I_K} = \omega_i$, where $\psi_{ij}$ is a
crystalline character whose $\kappa_{i'j'}$-labeled Hodge--Tate is $1$
if $(i',j')=(i,j)$ and is $0$ otherwise.  For this see
\cite[Prop.~B.3]{conradlifting} and the proof of
\cite[Prop.~6.7]{GLS-BDJ}(1).
\end{proof}

We write $\Delta(\lambda_1,\ldots,\lambda_d)$ for the diagonal matrix with diagonal entries
$\lambda_1,\ldots,\lambda_d$.

\begin{prop}
  \label{prop:rank-one-subs}
Assume that $p \ge 3$, let $\M$ be the Kisin module corresponding
to a lattice in a pseudo-BT
representation $V$ of weight $\{r_i\}$, and write $\barM = \M \otimes_{\O_E} k_E.$

Suppose that $\barN \subset \barM$ is a sub-$\varphi$-module such that
$\barM/\barN$ is free of rank one as an  $\fS_{\O_E} \otimes_{\O_E} k_E$-module.  Then
$\barN \simeq \barM(s_0,\ldots,s_{f-1};a)$ for some $a \in
k_E^{\times}$, with $s_i = r_i + x_i$ or
$s_i = e-1-x_i$ for some $x_i \in [0,e-1]$ for all $i$.
\end{prop}

\begin{proof}
Choose a basis $\{\e_{i},\f_i\}$ for $\M_i$ as in Theorem~\ref{thm: structure}.
 Since we will work in $\barM$ for the remainder of the proof, no
 confusion will arise if we write $\{\e_{i},\f_i\}$ also for the image of
   that basis in $\barM$.

 A generator
$e_{i-1}$ of $\barN_{i-1}$ has the form $(\e_{i-1},\f_{i-1})
 \cdot (v,w)^T$ for some $v,w
 \in k_E\llb u \rrb$, by hypothesis at least one of which is a unit.  We know from
 Theorem~\ref{thm: structure} that
$$ \varphi(e_{i-1}) = (\e_i,\f_i) \overline{X}_{i} \left(\prod_{j=1}^{e-1}
\overline{\Lambda}_{i,e-j} \overline{Z}_{i,e-j}\right)
\overline{\Lambda}_{i,0} \cdot  (\varphi(v),\varphi(w))^{T}$$
where $\overline{X}_{i}$ and $\overline{Z}_{ij}$ are the reductions
mod $\m_E$ of $X_{i}$ and $Z_{ij}$, where $\overline{\Lambda}_{i,0} =
\Delta(1,u^{r_i})$, and where $\overline{\Lambda}_{ij} = \Delta(1,u)$ for $1 \le j
\le e-1$.

Observe that each entry of
$(\varphi(v),\varphi(w))^T$ is either a unit or
divisible by $u^p$, and at least one is a unit.   Since we have
$r_i \le p$ for all $i$, it follows that the largest power of
$u$ dividing the column vector $
\overline{\Lambda}_{i,0} \cdot  (\varphi(v),\varphi(w))^{T}$ is either
$u^{r_i}$ or $u^0$.

For any $s \ge 0$,  if $y$ is a column vector of length 2 that is exactly divisible  by $u^s$,
it is easy to see that
that $\Delta(1,u) \cdot y$ is exactly divisible by either $u^{s}$ or
$u^{s+1}$.  On the other hand if $Z$ is invertible, then $Z \cdot y$ is
still exactly divisible by $u^{s}$.  Applying these observations
iteratively to the invertible matrices $\overline{X}_{i}$ and
$\overline{Z}_{ij}$, and to the matrices $\overline{\Lambda}_{ij} = \Delta(1,u)$ for $1 \le j
\le e-1$, we see that $\varphi(e_{i-1})$ is divisible exactly by
$u^{s_i}$ where $s_i = r_i + x'_i$ or $s_i = x'_i$ and $0 \le x'_i \le
e-1$ is the
number of times that we took $u^{s+1}$ rather than $u^s$ when
considering the effect of the matrix $\overline{\Lambda}_{ij}$.
Setting $x_i = x'_i$ in the first case and $x_i = e-1-x'_i$ in the
latter case, the proposition follows.
\end{proof}

If $V$ is a pseudo-BT representation of weight $\{r_i\}$ and $\lambda =
\kappabar_i \in \Hom(k,\Fpbar)$, we write $r_{\lambda} := r_i$.

\begin{thm}
  \label{thm:reducible-case}
  Assume that $p \ge 3$.  Let $T$ be a lattice in a pseudo-BT
  representation $V$ of weight $\{r_i\}$, and assume that $\overline{T} = T \otimes_{\O_E}
  k_E$ is reducible.  Then there is a subset $J \subseteq \Hom(k,\Fpbar)$ and integers $x_{\lambda} \in [0,e-1]$ such that
 $$ \overline{T}|_{I_K} \simeq
\begin{pmatrix} \prod_{ {\lambda} \in
    J}\omega_{{\lambda}}^{r_{\lambda} + x_{\lambda}} \prod_{{\lambda} \not\in J} \omega_{\lambda}^{e-1-x_{\lambda}} &*\\ 0& \prod_{ {\lambda}\not\in
    J}\omega_{{\lambda}}^{r_{\lambda} + x_{\lambda}} \prod_{{\lambda} \in J} \omega_{\lambda}^{e-1-x_{\lambda}}  \end{pmatrix} . $$
\end{thm}

\begin{proof}
Let $\fM$ be the Kisin module associated to the lattice $T$.  We have
$\overline{T}|_{G_{K_{\infty}}} \simeq T_{\fS}(\barM)$.
 From
(the proof of) \cite[Lem.~5.5]{GLS-BDJ}, we see that $\barM$ is
reducible and has a submodule $\barN$ as in
Proposition~\ref{prop:rank-one-subs} such that
$\overline{T}$ has a quotient character ${\chibar}$ with
${\chibar}|_{G_{K_{\infty}}} \simeq T_{\fS}(\barN)$.
Take $J = \{{\kappabar}_i \, : \, s_i = e - 1 - x_i\}$. (In particular,
if it happens that $e-1-x_i = r+x_i$ then we have $i \in J$.) The result follows from
Lemma~\ref{lem:functor-on-characters} (and a determinant argument to
compute the sub-character of $\overline{T}$). \end{proof}

We now give the  analogue of Theorem~\ref{thm:reducible-case} when $\overline{T}$ is
absolutely irreducible.  Let $k_2$ denote the unique quadratic extension of $k$ inside the
residue field of $\barK$.  We say that a subset $J \subseteq \Hom(k_2,\Fpbar)$ is
\emph{balanced} if for each ${\lambda} \in  \Hom(k_2,\Fpbar)$ exactly one
of ${\lambda}$ and ${\lambda}^{q}$ lies in $J$, with $q = p^f$.  If ${\lambda} \in
\Hom(k_2,\Fpbar)$, write $r_{{\lambda}}$ for $r_{{\lambda}|_k}$.   The
result is as follows.

\begin{thm}
  \label{thm:irreducible-case}
  Assume that $p \ge 3$.  Let $T$ be a lattice in a pseudo-BT
  representation $V$ of weight $\{r_i\}$, and assume that $\overline{T} = T \otimes_{\O_E}
  k_E$ is absolutely irreducible.  Then there is a balanced subset $J
  \subseteq \Hom(k_2,\Fpbar)$ and integers $x_{\lambda} \in [0,e-1]$ so
  that $x_{{\lambda}}$ depends only on ${\lambda}|_k$ and
 \numequation \label{eq:irredform} \overline{T}|_{I_K} \simeq
\prod_{ {\lambda} \in
    J}\omega_{{\lambda}}^{r_{\lambda} + x_{\lambda}} \prod_{{\lambda} \not\in J} \omega_{\lambda}^{e-1-x_{\lambda}} \bigoplus \prod_{ {\lambda}\not\in
    J}\omega_{{\lambda}}^{r_{\lambda} + x_{\lambda}} \prod_{{\lambda} \in J} \omega_{\lambda}^{e-1-x_{\lambda}} .\end{equation}
\end{thm}

\begin{proof}
  Note that $V$ restricted to the unramified quadratic extension $K_2$ of
  $K$ remains pseudo-BT, and for each embedding ${\kappa}' : K_2 \to
  \Qpbar$ extending ${\kappa} : K \to \Qpbar$ we have $\HT_{{\kappa}'}(V|_{G_{K_2}})
  = \HT_{{\kappa}}(V)$.  Applying Theorem~\ref{thm:reducible-case} to
  the lattice  $T|_{G_{K_2}}$ shows that $\overline{T}|_{I_K}$ has the
  form \eqref{eq:irredform}, except that $J$ need not be
  balanced, nor  must $x_{{\lambda}} = x_{{\lambda}^q}$. It remains to
  be seen that these additional conditions may be taken to hold.
Assume that $\chibar : I_K \to
  \Fpbar^{\times}$ is a character of niveau $2f$ such that $\chibar \oplus
  \chibar^q$ is equal to a representation as in the right-hand side of
  \eqref{eq:irredform}; to complete the proof, we wish to show that $\chibar \oplus
  \chibar^q$ is also equal to a representation as in
  \eqref{eq:irredform} with $J$
  balanced and $x_{{\lambda}} = x_{{\lambda}^q}$.  We apologise to the
  reader for the argument that follows, which is entirely elementary but long and unenlightening.

It follows from \cite[Cor.~4.1.20]{blggu2} that if $e \ge p$, then as $J$ varies over all
balanced sets and the $x_{{\lambda}} \in [0,e-1]$ vary over all
possibilities with $x_{{\lambda}} = x_{{\lambda}^q}$, the right-hand side of
\eqref{eq:irredform} exhausts all representations
$\chibar\oplus \chibar^q$ of niveau $2f$ with determinant
$\prod_{{\lambda}\in \Hom(k,\Fpbar)}
\omega_{\lambda}^{r_{\lambda}+e-1}$.  The theorem is then automatic in this case, since
$\overline{T}|_{I_K}$ must have this determinant as well.  For the
remainder of the argument, then, we assume that $e \le p-1$.

Fix a character ${\kappabar}'_0$ extending ${\kappabar}_0$, and define
${\kappabar}'_i$ for $i \in \Z$ by $({\kappabar}'_{i+1})^p = {\kappabar}'_i$.  Write
$x_i$ for $x_{{\kappabar}'_i}$ for $i \in \Z$.
Similarly define $r_i$ for all $i \in \Z$ so that $r_{i+f} = r_i$.
   Define
$$\begin{cases} J_1  = J \cap (f+J) \\
J_2  = J^c \cap (f+J)^c \\
J_3  = J \cap (f+J)^c \\
J_4  = J^c \cap (f+J)\end{cases}$$
so that $J$ is balanced if and only if $J_1 = J_2 = \varnothing$.  To
simplify notation we will for instance write $i \in J_1$ in lieu of ${\kappabar}'_i
\in J_1$, and we will write $\omega_i$ for
$\omega_{{\kappabar}'_i}$. (We stress that the symbol $\omega_i$ here
denotes a fundamental character of level $2f$, and that the symbols
$\omega_i$ and $x_i$ are indexed modulo $2f$.)
The condition that the two summands on the right-hand side of \eqref{eq:irredform} are $q$th powers of
one another is equivalent to
$$ \prod_{i \in J_1} \omega_i^{r_i + x_i + x_{i+f} - (e-1)}
 \prod_{i \in J_2} \omega_i^{(e-1) - r_i - x_i - x_{i+f}}
 \prod_{i \in J_3} \omega_i^{x_i - x_{i+f}}
 \prod_{i \in J_4} \omega_i^{x_{i+f} - x_{i}} = 1.$$
Write $y_i$ for the exponent of $\omega_i$ in the above expression.
Since $r_i \in [1,p]$ for all $i$ and $e \le p-1$, we have
$$\begin{cases} y_i \in [-p+3,2p-2] & \text{if } i \in J_1 \\
y_i \in [-2p+2,p-3] & \text{if } i \in J_2 \\
y_i \in [-p+2,p-2] & \text{if } i \in J_3 \\
y_i \in [-p+2,p-2] & \text{if } i \in J_4.\end{cases}$$
Note that $y_i = y_{i+f}$ for all~$i$, so we can consider the $y_i$'s
as being labeled cyclically with index taken modulo~$f$.
As in the proof of \cite[Lem.~7.1]{GLS-BDJ}, one checks that since $\prod_i \omega_i^{y_i} = 1$
with $y_i \in [-2p+2,2p-2]$ for all~$i$,
 the tuple
$(y_0,\ldots,y_{f-1})$ must have the shape
$$ a_0 (p,0,\ldots,0,-1) + a_1 (-1,p,0,\ldots,0) + \cdots + a_{f-1} (0,\ldots,0,-1,p) $$
with $|a_i| \in \{0,1,2\}$ for all $i$, and in fact either $a_i = 2$
for all $i$, or $a_i = -2$ for all $i$, or else $a_i \in \{0,\pm 1\}$
for all $i$.
But if $a_i = 2$ for all $i$, so that $y_i = 2p-2$ for all $i$, we would have to have $J = \{0,\ldots,2f-1\}$ with $r_i=p$
and $e=x_i=p-1$ for all $i$.  But then $\chibar^q = \chibar$, i.e. $\chibar$ has niveau $f$ rather
than niveau $2f$, a contradiction.  The case where $a_i = -2$ for all $i$ is similarly
impossible.  So in fact we must have $a_i \in \{0,\pm 1\}$ for
all~$i$.
We now consider separately the case where some $a_i$ is equal to
$0$, and the case where $a_i = \pm 1$ for all $i$.

First let us suppose that at least one $a_i$ is equal to $0$.  The cyclic set of those $i$ with $y_i
\neq 0$ (with index $i$ taken modulo $f$) must break up as a disjoint union of sets of the form
$(i,i+1,\ldots,i+j)$ with $y_i = \pm 1$, $y_{i+j} = \pm p$, and $y_{i+1},\ldots,y_{i+j-1}
\in \{\pm p \pm 1\}$.  For every such interval $[i,i+j]$, choose a representative of
$i$ modulo $2f$ (that we also denote $i$), and perform the following
operation (noting that since $y_i\in [-p+2,p-2]$ for $i\in J_3\cup
J_4$, we have $i+1,\dots,i+j\in J_1\cup J_2$):
\begin{itemize}
\item Replace $J$ with $J \triangle \{i,\ldots,i+j\}$ if $i \in J_1
  \cup J_2$, or else with $J \triangle \{i+1,\ldots,i+j\}$ if $i \in
  J_3 \cup J_4$;
\item Replace $x_{\ell}$ with $x_{\ell+f}$ for each $\ell \in [i,i+j]$.
\end{itemize}
Here $\triangle$ denotes symmetric difference. It is easy to check that this operation does not change $\chibar$, and
that for the new choice of $J$ and $x_i$'s we have $y_i = 0$ for all
$i$.  We now have $r_i + x_i = e-1-x_{i+f}$ for each $i \in J_1 \cup J_2$,
so for each pair $\{i,i+f\} \subseteq J_1 \cup J_2$ we can again
replace $J$ with $J \triangle \{i\}$ and $x_i$ with $x_{i+f}$ without
changing $\chibar$.  When this operation is complete, our new set $J$ is
balanced.  Furthermore $y_i = 0$ for all $i$, and so $x_i = x_{i+f}$
for all $i$, and this case is complete.

Finally suppose that $a_i = \pm 1$ for all $i$.  Then $J_1 \cup J_2 =
\{0,\ldots,2f-1\}$, and in fact $i \in J_1$ if $a_i = 1$ while $i \in J_2$ if $a_i
= -1$.  Note
that if $i \in J_1$ and $i+1 \in J_2$ or vice-versa, then $r_i + x_i + x_{i+f}
- (e-1) = p+1$ (so in particular both $x_i,x_{i+f}$ are
nonzero), while if $i,\ i+1 \in J_1$ or $i,\ i+1 \in J_2$, then $r_i + x_i
+ x_{i+f} - (e-1) = p-1$.

By symmetry we can suppose without loss of generality that
$J_1 \neq \varnothing$.  If $J_2 = \varnothing$, then some $x_i$ is
not equal to $e-1$ (otherwise we have $x_i = e-1$ for all $i$, so that
$r_i = p-e$ for all $i$, and $\chibar$ has niveau $f$); changing our
choice of ${\kappabar}_0$ (if necessary) we suppose that $x_{f-1} \neq
e-1$.  On the other hand if $J_2 \neq \varnothing$, then by changing
our choice of ${\kappabar}_0$ (if necessary) we suppose that $0 \in J_1$
but $f-1 \in J_2$.

Take $J' = \{0,\ldots,f-1\}$, and for each
$0 \le i \le f-2$ we set
$$x'_i = x'_{i+f} = \begin{cases}
x_i & \text{ if } i \in J_1,\ i+1 \in J_1 \\
x_i-1 & \text{ if } i \in J_1,\ i+1 \in J_2 \\
x_{i+f}-1 & \text{ if } i \in J_2,\ i+1 \in J_1 \\
x_{i+f} & \text{ if } i \in J_2,\ i+1 \in J_2.\\
\end{cases}$$
We take $x'_{f-1}=x'_{2f-1} = x_{f-1}+1$ if $f-1\in J_1$, or
$x'_{f-1} = x'_{2f-1} = x_{2f-1}$ if $f-1 \in J_2$.  (In other words
if $i=f-1$ we take the common value of $x'_{f-1},\ x'_{2f-1}$ to be
$1$ more than
the value that would otherwise have been given by the above table.)  Note that we have $x'_i \in [0,e-1]$
for all $i$ by the observations and choices in the two preceding paragraphs.

We claim that $\chibar' := \prod_{i=0}^{f-1} \omega_i^{r_i +
  x'_i} \prod_{i=f}^{2f-1}\omega_i^{e-1-x'_i}$ is equal
to $\chibar$.  Since $J'$ is balanced and $x'_i = x'_{i+f}$ for all $i$, this will complete the proof.
Checking the claim is somewhat laborious, and we only give an
indication of the argument.  We wish to show that $\chibar' \chibar^{-1}$ is
trivial.  Write $\chibar' \chibar^{-1} = \prod_i \omega_i^{z_i}$ using the
defining formulas for $\chibar$ and $\chibar'$.  The values of the $z_i$ are
calculated by considering eight cases, depending on whether or not $i \in
\{0,\ldots,f-1\}$, whether or not $i \in J_1$, and whether or not $i+1
\in J_1$ (as well as making an adjustment by $1$ when $i=f-1,\ 2f-1$).  For instance if $i \in \{f,\ldots,2f-1\}$, $i \in J_1$, and
$i+1 \in J_2$ then we have
$$z_i = (e - 1 - (x_{i+f}-1)) - (r_i + x_i) = -p,$$
while if $i \in \{f,\ldots,2f-1\}$, $i \in J_1$, and
$i+1 \in J_1$ then $z_i = 1-p$ if $i \neq 2f-1$ and $z_i = -p$ if $i =
2f-1$.  (In all cases one finds that $|z_i| \in \{0, 1,p-1,p\}$ and $z_i$ depends only on $i$, $J_1$, and $J_2$, not on the $r_i$'s or
$x_i$'s.)  It is then straightforward to check that $\prod_i
\omega_i^{z_i} = 1$ (one just has to ``carry'' by replacing every $\omega_i^{\pm p}$ with~$\omega_{i-1}^{\pm 1}$).
\end{proof}

\section{The weight part of Serre's conjecture I: the semisimple case}
\label{sec:appl-weight-part}

For a detailed discussion of the weight part of Serre's conjecture for
$\GL(2)$, we refer the reader to~\cite[\S4]{blggu2}. In this
section, we will content ourselves with a brief explanation of the
consequences of the results of the previous sections for the weight
part of Serre's conjecture for (definite or indefinite) quaternion
algebras over totally real fields; we note that the analogous results
for compact at infinity unitary groups $\U(2)$ over CM fields follow
immediately from~\cite[Thm.~5.1.3]{blggu2}, together with the
discussion here.

\subsection{Local Serre weights}
 In order to make our various definitions associated to Serre weights,
it will be convenient to work in the local setting of the previous
sections.

\begin{defn}\label{defn:serre weight}
 A \emph{Serre weight} of $\GL_2(k)$ is by definition an
irreducible $\Fpbar$-representation of $\GL_2(k)$, which is
necessarily of the
form \[\sigma_{a,b} :=\otimes_{\lambda\in\Hom(k,\Fpbar)}\det{}^{b_\lambda}\otimes\Sym^{a_\lambda-b_\lambda}k^2\otimes_{k,\lambda}\Fpbar,\]for
some (uniquely determined) integers $a_\lambda,b_\lambda$ with
$b_\lambda,\ a_\lambda-b_\lambda \in [0,p-1]$ for all~$\lambda$,
and not all $b_{\lambda}$ equal to $p-1$.
  \end{defn}
Note that $\sigma_{a,b}$ has a natural model
$\otimes_{\lambda\in\Hom(k,k_E)}\det{}^{b_\lambda}\otimes\Sym^{a_\lambda-b_\lambda}k^2\otimes_{k,\lambda}
k_E$, and it will sometimes be convenient for us to think of
$\sigma_{a,b}$ as being defined over $k_E$ (or rather, it will be
convenient for us to identify $\Hom(k,k_E)$ with $\Hom(k,\Fpbar)$).

\para Suppose that $\rbar : G_K \to \GL_2(\Fpbar)$ is continuous.  In~\cite{blggu2}
and~\cite{GeeKisin} there are definitions of  several sets of
Serre weights $\Wexplicit(\rbar)$,
$\WBT(\rbar)$, and $\Wcris(\rbar)$. We now recall the
definitions of $\Wexplicit(\rbar)$ and
$\Wcris(\rbar)$; see~\cite[Def.~4.5.6]{GeeKisin} for $\WBT(\rbar)$.

Write $a_i,\ b_i$ in place of
$a_{\kappabar_i},\ b_{\kappabar_i}$.
We say that a de Rham lift $r$ of $\rbar$ has \emph{Hodge
  type $\sigma_{a,b}$} if for all $0 \le i \le f-1$
  we have $\HT_{\kappa _{i,0}} (r) =\{ b_i,a_i+1 \}$, and if
  $\HT_{\kappa_{ij}} (r)=\{0, 1 \}$ when $j \neq 0$.  Note that a crystalline representation of Hodge type
    $\sigma_{a,0}$ is pseudo-BT of weight $\{a_i+1\}$.
  \begin{defn}(\cite[Def.~4.1.4]{blggu2})
    \label{defn: Wcris}
    $\Wcris(\rbar)$ is the set of Serre weights
    $\sigma_{a,b}$ for which $\rbar$ has a crystalline lift
    of Hodge type $\sigma_{a,b}$.
  \end{defn}

As in the previous section let $k_2$ denote the unique quadratic extension of $k$ inside the
residue field of $\barK$.

  \begin{defn}(\cite[Def.~4.1.23]{blggu2})
  \label{defn: Wexplicit}If $\rbar$ is irreducible, then
  $\Wexplicit(\rbar)$ is the set of Serre weights $\sigma_{a,b}$ such
  that there is a balanced subset $J\subset\Hom(k_2,\Fpbar)$, and for each
  $\lambda \in\Hom(k,\Fpbar)$ an integer $0\le x_\lambda\le e-1$ such that
  if we write
  $x_{\lambda}$ for
  $x_{\lambda|_k}$ when $\lambda \in \Hom(k_2,\Fpbar)$, then \[\rbar|_{I_K}\cong
  \begin{pmatrix}\prod_{\lambda\in
      J}\omega_{\lambda}^{a_{\lambda}+1+x_{\lambda}}\prod_{\lambda\not\in
      J}\omega_\lambda^{b_{\lambda}+e-1-x_{\lambda}}&0\\ 0& \prod_{\lambda\not\in
      J}\omega_\lambda^{a_{\lambda}+1+x_{\lambda}}\prod_{\lambda\in
      J}\omega_\lambda^{b_{\lambda}+e-1-x_{\lambda}}

  \end{pmatrix}.\]

If $\rbar$ is reducible, then $\Wexplicit(\rbar)$ is
the set of weights $\sigma_{a,b}$ for which $\rbar$
   has a crystalline lift of type $\sigma_{a,b}$ of the
  form \[ \begin{pmatrix}\chi' &*\\ 0& \chi
  \end{pmatrix}.\]
 In particular (see the remark after~\cite[Def.~4.1.14]{blggu2}), if $\sigma_{a,b} \in \Wexplicit(\rbar)$ then it is
  necessarily the case that there is a subset
  $J \subseteq \Hom(k,\Fpbar)$ and for each $\lambda\in \Hom(k,\Fpbar)$ there is an integer
  $0\le x_\lambda\le e-1$ such that \[\rbar|_{I_K}\cong
  \begin{pmatrix}\prod_{\lambda\in
      J}\omega_\lambda^{a_{\lambda}+1+x_\lambda}\prod_{\lambda\not\in
      J}\omega_\lambda^{b_{\lambda}+x_\lambda}&*\\ 0& \prod_{\lambda\not\in
      J}\omega_\lambda^{a_{\lambda}+e-x_\lambda}\prod_{\lambda\in
      J}\omega_\lambda^{b_{\lambda}+e-1-x_\lambda}
 \end{pmatrix},\]
and when $\rbar$ is a sum of two characters this is necessary and sufficient.
  \end{defn}

\para The inclusion
  $\Wexplicit(\rbar)\subseteq\Wcris(\rbar)$ was proved in
  \cite[Prop.~4.1.25]{blggu2} and was conjectured there
  to be an equality (\cite[Conj.~4.1.26]{blggu2});
  this equality is the main local result of this paper.

  The definition of $\WBT(\rbar)$ is unimportant for us
in this paper; the only fact that we will need is that by~\cite[Cor.~4.5.7]{GeeKisin}, the inclusion
$\Wexplicit(\rbar)\subseteq\Wcris(\rbar)$
can be refined to
inclusions \[ \Wexplicit(\rbar)\subseteq\WBT(\rbar)\subseteq\Wcris(\rbar),\]
so that our main result will show that both these inclusions are
equalities.

  Observe
  that the definitions of the sets $\Wcris(\rbar)$ and $\Wexplicit(\rbar)$
  involve our fixed choice of embeddings $\kappa_{i,0}$.  We will prove
  that these sets in fact do not depend on this choice.
  Indeed the definition of $\WBT(\rbar)$ does not
  involve any choice of embeddings $\kappa_{i,0}$, so that once we have
  proved the equality $\Wexplicit(\rbar) =\Wcris(\rbar)$ it
  follows
  automatically that $\Wexplicit(\rbar)$ and $\Wcris(\rbar)$  do
  not depend on any choice of embeddings either. On the other hand this
  will also follow easily and directly from the arguments in this
  paper, and so we will give a direct proof as well.

In the above language, the main local result of this section is the
following.

\begin{thm}\label{thm:our main semisimple local result in Bexplicit
    language}Suppose that $p\ge 3$ and that $\rbar$ is
  semisimple. Then we have
  $\Wexplicit(\rbar)=\Wcris(\rbar)$.  Moreover,  these sets do not depend on
  our choice of embeddings~$\kappa_{i,0}$.
\end{thm}
\begin{proof} Suppose that $\sigma_{a,b}\in\Wcris(\rbar)$. Twisting, we
  may assume that $b_{\lambda} =0$ for all~$\lambda$.  If $r  : G_K \to
  \GL_2(\cO_E)$ is a pseudo-BT lift of $\rbar$ of Hodge type $\sigma_{a,0}$, we
  can freely enlarge the coefficient field $E$ so that it satisfies
  our usual hypotheses, and so that $r \otimes_{\cO_E} k_E$ is
  either reducible or absolutely irreducible.

 Now for the first part of the result it remains to show that
 $\Wcris(\rbar)\subseteq\Wexplicit(\rbar)$; this
  is immediate from Theorems~\ref{thm:reducible-case}
  and~\ref{thm:irreducible-case} and Definitions~\ref{defn: Wcris}
  and~\ref{defn: Wexplicit}, taking $r_i = a_i + 1$ for all $i$.
  The second part is automatic since the definition of $\Wexplicit(\rbar)$ when $\rbar$ is semisimple
does not depend on the choice of embeddings~$\kappa_{i,0}$.
\end{proof}

\subsection{Global Serre weights}  Let $F$ be a totally real field, and continue to assume that
$p>2$. Let $\rhobar:G_F\to\GL_2(\Fpbar)$ be continuous, absolutely irreducible,
and modular (in the sense that it is isomorphic to the reduction
modulo $p$ of a $p$-adic Galois representation associated to a Hilbert
modular eigenform of parallel weight two). Let $k_v$ be the residue
field of $F_v$ for each $v|p$.

 A \emph{global Serre weight} is by definition an irreducible
representation of the group $\prod_{v|p}\GL_2(k_v)$, which is necessarily of the
form $\sigma=\otimes_{v|p}\sigma_v$ with $\sigma_v$ a Serre weight of
$\GL_2(k_v)$ as above. Let $D$ be a quaternion algebra with centre
$F$, which is split at all places dividing $p$ and at zero or one
infinite places. Then~\cite[Def.~5.5.2]{GeeKisin} explains what
it means for $\rhobar$ to be modular for $D$ of weight
$\sigma$. (There is a possible local obstruction at the finite places
of $F$ at which $D$ is ramified to $\rhobar$ being
modular for $D$ for any weight at all; following~\cite[Def.~5.5.3]{GeeKisin}, we say that $\rhobar$ is \emph{compatible} with
$D$ if this obstruction vanishes. Any $\rhobar$ will be compatible
with some $D$; indeed, we could take $D$ to be split at all finite
places of $F$.) The following is a theorem of Gee--Kisin~\cite[Cor.~5.5.4]{GeeKisin}.
\begin{thm}\label{thm:GK main result}Assume that $p>2$, that $\rhobar$ is modular and compatible
  with $D$, that $\rhobar|_{G_{F(\zeta_p)}}$ is irreducible, and if
  $p=5$ assume further that the projective image
  of $\rhobar|_{G_{F(\zeta_p)}}$ is not isomorphic to $A_5$.

Then $\rhobar$ is modular for $D$ of weight $\sigma$ if and only if
$\sigma_v\in\WBT(\rhobar|_{G_{F_v}})$ for all $v|p$.
\end{thm}
The main global result of this section is the following (which is an immediate
consequence of Theorem~\ref{thm:our main semisimple local result in Bexplicit
  language} in combination with the above result of Gee--Kisin).

\begin{thm}\label{thm:our main semisimple global result}Assume that $p>2$, that $\rhobar$ is modular and compatible
  with $D$, that $\rhobar|_{G_{F(\zeta_p)}}$ is irreducible, and if
  $p=5$ assume further that the projective image
  of $\rhobar|_{G_{F(\zeta_p)}}$ is not isomorphic to $A_5$.

  Assume that $\rhobar|_{G_{F_v}}$ is semisimple for each place
  $v|p$. Then $\rhobar$ is modular for $D$ of weight $\sigma$ if and
  only if $\sigma_v\in\Wexplicit(\rhobar|_{G_{F_v}})$ for all $v|p$.
\end{thm}

Modulo the hypotheses on the image of $\rhobar$ (the usual hypotheses needed in order to
apply the Taylor--Wiles--Kisin method), Theorem~\ref{thm:our main semisimple
  global result} is the main conjecture of Schein
\cite{ScheinRamified}.

Note that we have thus far said nothing about the case where
$\rhobar|_{G_{F_v}}$ is reducible but non-split, where
$\Wexplicit(\rhobar|_{G_{F_v}})$ will depend on the extension
class of $\rhobar|_{G_{F_v}}$.  Our treatment of this more delicate case will occupy
the remainder of the paper.

\begin{remark}
  \label{rem:ds}
  If $\rhobar|_{G_{F_v}}$ is generic in a suitable sense, then
  \cite[Thm.~4.5]{DiamondSavitt} implies that $\Wexplicit(\rhobar|_{G_{F_v}}) = \WBT(\rhobar|_{G_{F_v}}) \cap
 \Wexplicit(\rhobar|^{\textrm{ss}}_{G_{F_v}})$.  On the other hand we
 have inclusions
$$ \WBT(\rhobar|_{G_{F_v}}) \subseteq  \Wcris (\rhobar|_{G_{F_v}})
\subseteq \Wcris(\rhobar|^\textrm{ss}_{G_{F_v}})  =
\Wexplicit(\rhobar|^\textrm{ss}_{G_{F_v}})  $$
where the equality is an application of Theorem~\ref{thm:our main semisimple local result in Bexplicit
    language}, and so we see in this case that
  $\Wexplicit(\rhobar|_{G_{F_v}}) = \WBT(\rhobar|_{G_{F_v}})$.   In
  particular we can already extend Theorem~\ref{thm:our main
    semisimple global
    result} to the case where $\rhobar|_{G_{F_v}}$ is either
  semisimple or generic for all $v|p$.

(We
  refer the reader
  \cite[Def.~3.5]{DiamondSavitt} for the definition of genericity that we
  use here, but note
  that genericity implies that $e \le (p-1)/2$.)
\end{remark}

\section{The weight part of Serre's conjecture II: the non-cyclotomic case}
\label{sec:non-cyclotomic}

In this section we make a detailed study of the extensions of rank-one
Kisin modules and use these results to prove that if $p \ge 3$ then
$\Wexplicit(\rbar) = \Wcris(\rbar)$ for reducible representations
$\rbar \simeq \begin{pmatrix} \chibar' & * \\ 0 &
  \chibar \end{pmatrix}$ with $\chibar^{-1} \chibar' \neq
\varepsilonbar$ (see Theorem~\ref{thm:non-cyclotomic} below).

\subsection{Extensions of rank-one Kisin modules}
\label{sec:extensions}

We begin with some basic results on extensions of rank-one Kisin modules, mildly generalising the results in
\cite[\S7]{GLS-BDJ}.  We begin with the following notation.

\begin{defn}
  \label{defn:alpha}
 If $\barN = \barM(s_0,\ldots,s_{f-1};a)$, define $\alpha_i(\barN) := \frac{1}{p^f-1}
 \sum_{j=1}^{f} p^{f-j} s_{j+i}$.
\end{defn}

It is immediate from the definition that these constants satisfy the
relations $\alpha_i(\barN) + s_i  = p
\alpha_{i-1}(\barN)$ for $i=0,\ldots,f-1$ and indeed are uniquely
defined by them.  We have the following easy lemma.

\begin{lem}
  \label{lem:map}  Write $\barN = \barM(s_0,\ldots,s_{f-1};a)$ and
  $\barN' = \barM(s_0',\ldots,s_{f-1}';a')$.  There exists a nonzero
  map $\barN \to \barN'$ if and only if $\alpha_i(\barN) -
  \alpha_i(\barN') \in \Z_{\ge 0}$ for all $i$, and $a=a'$.
\end{lem}

\begin{proof}
  Such a map, if it exists, must take the form $e_i \mapsto
  c u^{\alpha_i(\barN)-\alpha_i(\barN')} e'_i$ for all $i$, where $e'_i$ is the
  basis element for $\barN'_i$ as in Definition~\ref{defn:rank-one},
  and $c \in k_E^{\times}$.
\end{proof}

We can now check the following analogue of \cite[Prop.~7.4]{GLS-BDJ}.

\begin{prop}
  \label{prop:kisin-module-extensions}
  Let $\barN = \barM(s_0,\ldots,s_{f-1};a)$ and $\barP =
  \barM(t_0,\ldots,t_{f-1};b)$ be rank-one Kisin modules, and let
  $\barM$ be an extension of $\barN$ by $\barP$.  Then we can choose
  bases $e_i,f_i$ of the $\barM_i$ so that $\varphi$ has the form
  \begin{eqnarray*}
    \varphi(e_{i-1}) & = &(b)_i u^{t_i} e_i \\
    \varphi(f_{i-1}) & = &(a)_i u^{s_i} f_i + y_i e_i
  \end{eqnarray*}
with $y_i \in k_E\llb u \rrb$ a polynomial with $\deg(y_i) < s_i$,
except that when there is a nonzero map $\barN \to \barP$ we must
also allow $y_j$ to have a term of degree $s_j + \alpha_j(\barN) -
\alpha_j(\barP)$ for any one choice of $j$.
\end{prop}

\begin{proof}
  The proof is essentially identical to the proof of \cite[Prop.~7.4]{GLS-BDJ}
  except for the analysis of the exceptional terms.  Namely, it is
  possible to use a simultaneous change of basis of the form $f'_i = f_i + z_i e_i$
  for $i = 0,\ldots,f-1$ to eliminate all terms of degree at least
  $s_i$ in the $y_i$'s, except that if there is a sequence of
  integers $d_i \ge s_i$ satisfying
  \numequation \label{eq:d-recursion} d_i = p(d_{i-1} - s_{i-1}) + t_i \end{equation}
 for all $i$ (with subscripts taken modulo $f$), and if $a =
 b$,
 then for some $j$ the term of degree $d_j$ in $y_j$ may survive.
 (Note that taking $z_{i-1}$ to be a suitable monomial of degree $d_{i-1}-s_{i-1}$
 will eliminate the term of degree $d_{i-1}$ in $y_{i-1}$ but yield a term of
 degree $d_{i}$ in $y_{i}$ instead.  This is more-or-less what we called a
 \emph{loop} in the proof of  \cite[Prop.~7.4]{GLS-BDJ}.)

Comparing with the defining relations for $\alpha_i(\barN)$,
$\alpha_i(\barP)$, we must have
$$ d_i = p\alpha_{i-1}(\barN) -\alpha_i(\barP) =
\alpha_i(\barN)-\alpha_i(\barP) + s_i,$$
and these are integers greater than or equal to $s_i$ provided that
$\alpha_i(\barN)-\alpha_i(\barP) \in \Z_{\ge 0}$ for all $i$.  The
result now follows from Lemma~\ref{lem:map}.
\end{proof}

Next we combine Proposition~\ref{prop:kisin-module-extensions} with
Theorem~\ref{thm: structure} to study the extensions of Kisin modules
that can arise from the reduction mod $p$ of pseudo-BT
representations; this is our analogue of \cite[Thm.~7.9]{GLS-BDJ}.

\begin{thm}
  \label{thm:pseudo-BT-extensions}
  Suppose that $p \ge 3$.  Let $T$ be a $G_K$-stable $\cO_E$-lattice
  in a pseudo-BT representation of weight $\{r_i\}$.  Let $\M$ be the
  Kisin module associated to $T$, and let $\barM := \M \otimes_{\cO_E}
  k_E$.

  Assume that the $k_E[G_K]$-module $\overline{T} := T \otimes_{\cO_E}
  k_E$ is reducible, so that there exist rank-one Kisin modules $\barN =
  \barM(s_0,\ldots,s_{f-1};a)$ and $\barP =
  \barM(t_0,\ldots,t_{f-1};b)$ such that $\barM$ is an extension of
  $\barN$ by $\barP$.   Then for all $i$ there is an integer $x_i \in
  [0,e-1]$ such that $\{s_i,t_i\} = \{r_i + x_i, e-1-x_i\}$.

 We can choose
  bases $e_i,f_i$ of the $\barM_i$ so that $\varphi$ has the form
  \begin{eqnarray*}
    \varphi(e_{i-1}) & = &(b)_i u^{t_i} e_i \\
    \varphi(f_{i-1}) & = &(a)_i u^{s_i} f_i + y_i e_i
  \end{eqnarray*}
where
\begin{itemize}
\item $y_i \in k_E\llb u \rrb$ is a polynomial with $\deg(y_i) < s_i$,
\item if $t_i < r_i$ then the nonzero terms of $y_i$ have degrees
  in the set $\{t_i\} \cup [r_i,s_i-1]$,
\item except that when there is a nonzero map $\barN \to \barP$ we must
also allow $y_j$ to have a term of degree $s_j + \alpha_j(\barN) -
\alpha_j(\barP)$ for any one choice of $j$.
\end{itemize}
\end{thm}

\begin{proof}
  The fact that $\barM$ is an extension of two rank-one Kisin modules
  follows e.g.\ from \cite[Lem.~5.5]{GLS-BDJ}, and the fact that
  $\{s_i,t_i\} = \{r_i+x_i,e-1-x_i\}$ is an immediate consequence of
  Proposition~\ref{prop:rank-one-subs} and determinant considerations.  After applying
  Proposition~\ref{prop:kisin-module-extensions}, it remains to check
  that if  $t_i < r_i$ then the nonzero terms of $y_i$ have degrees
  in the set $\{t_i\} \cup [r_i,s_i-1]$ (except possibly for an exceptional term
  arising from the existence of a nonzero map $\barN \to \barP$).

 It is an immediate consequence of Theorem~\ref{thm: structure} that $\varphi(\barM_{i-1})$
 contains a saturated element (i.e.\ an element not divisible by $u$
 in $\varphi(\barM_{i-1})$) that is divisible by $u^{r_i}$ in
 $\barM_i$.  Such an element must be a saturated $\varphi(\fS)$-linear
 combination of $\varphi(e_{i-1})$ and $\varphi(f_{i-1})$, i.e.\ we
 must have $\gamma,\delta$, at least one of which is a unit, such that
 $u^{r_i}$ divides
 $$ \varphi(\gamma) \cdot ((a)_i u^{s_i}
 f_i + y_i e_i) +  \varphi(\delta) \cdot (b)_i u^{t_i} e_i.$$
The assumptions on $\{s_i,t_i\}$ imply
 that one of them is at least $r_i$.
 Suppose, then,
 that $t_i < r_i$.  It follows that $s_i \ge r_i$, so we have to have $$u^{r_i} \mid \varphi(\gamma) y_i +
 \varphi(\delta) (b)_i u^{t_i}.$$ If
 $\gamma$ were a non-unit, then we would have $u^p \mid
 \varphi(\gamma)$, so that $u^{r_i}$ divides both terms in the above
 sum separately.  But
 $t_i < r_i$, and so $\delta$ must also be a non-unit, a
 contradiction.  It follows that $\gamma$ is a unit.   Replacing
 $\delta$ with $\delta \gamma^{-1}$ we can suppose $\gamma = 1$, and
 we must be able to choose $\delta$ so that $u^{r_i}$ divides $y_i +
 \varphi(\delta)(b_i) u^{t_i}$.  That is, all terms of $y_i$ must have
 degree at least $r_i$, except for terms that could be canceled by the
 addition of  $\varphi(\delta)(b_i) u^{t_i}$.   These are the terms of
 degree $t_i + pn$ with $n \in \Z_{\ge 0}$.  But $t_i+p
 \ge r_i$, so the only extra term we may get this way is a term
 of degree~$t_i$.
\end{proof}

\begin{defn}
  \label{defn:pseudo-BT-classes}  Fix integers $r_i \in [1,p]$.
   Suppose that $\barN =
  \barM(s_0,\ldots,s_{f-1};a)$ and $\barP =
  \barM(t_0,\ldots,t_{f-1};b)$ are rank-one Kisin modules with
  $\{s_i,t_i\} = \{r_i + x_i, e-1-x_i\}$.  We let $\EpsBT(\barN,\barP)$ denote the subset of
 $\Ext^1(\barN,\barP)$ defined by the conditions of
 Theorem~\ref{thm:pseudo-BT-extensions}.  (The integers $\{r_i\}$ will
 be implicit in the notation.)  \end{defn}

From the proof of Theorem~\ref{thm:pseudo-BT-extensions} we see that
 $\EpsBT(\barN,\barP)$ can be characterised as the set of classes
 $\barM \in  \Ext^1(\barN,\barP)$ such that each $\varphi(\barM_{i-1})$
 contains a saturated element that is divisible by $u^{r_i}$ in
 $\barM_i$.

In the remainder of Section~\ref{sec:non-cyclotomic}, we prove that if $p \ge 3$ then
$\Wexplicit(\rbar) = \Wcris(\rbar)$ for reducible representations
$\rbar \simeq \begin{pmatrix} \chibar' & * \\ 0 &
  \chibar \end{pmatrix}$ with $\chibar^{-1} \chibar' \neq
\varepsilonbar$. We follow the same basic strategy as in
\cite[\S\S8-9]{GLS-BDJ}, which relies on (but, as we will explain, is necessarily somewhat
more complicated than)  a comparision between the
  the dimension of the space  $\EpsBT(\barN,\barP)$ and the
  dimension of an appropriate local Bloch--Kato group
  $H^1_f(G_K,-)$. Perhaps the main difference between the arguments
  here and the arguments of \cite[\S\S8-9]{GLS-BDJ} is that certain
  combinatorial issues which we were able to address in ad hoc ways
  in the unramified case (particularly those of \cite[\S8.2]{GLS-BDJ})
  are rather more involved in the ramified case, and so have had to be addressed
  systematically (\textit{cf.}\ Section~\ref{ss:maximalminimal}).

\subsection{Comparison of extension classes}
\label{sec:comparison}
We give a simple example that illustrates why the
  proof that $\Wexplicit(\rbar) = \Wcris(\rbar)$ in the indecomposable
  case is more complicated than one might initially
  expect, and in particular cannot follow immediately from comparing
  the dimension of the space  $\EpsBT(\barN,\barP)$ with the
  dimension of an appropriate local Bloch--Kato group $H^1_f(G_K,-)$.

\begin{example}
  \label{ex:structure-of-argumemt}

 Take $K = \Qp$.  The group $H^1_f(G_{\Qp},\cO_E(\varepsilon^{1-p}))$
 is torsion (since $H^1_f(G_{\Qp},E(\varepsilon^{1-p}))$ is trivial) and its
 $\varpi$-torsion has rank one, corresponding to the congruence
 $\varepsilon^{1-p} \equiv 1 \pmod{\varpi}$.  It follows that
 the subspace of $\Ext^1_{k_E[G_{\Qp}]}(\overline{1},\overline{1})$
 arising from crystalline extensions of $\varepsilon^{p-1}$ by $1$ is
 one-dimensional.  (This extension comes from the reduction mod
 $\varpi$ of a nontrivial extension $T$ of $\O_E(\varepsilon^{p-1})$
 by $\O_E$ inside the split representation
 $\varepsilon^{p-1} \oplus 1$.)  By Lemma~\ref{lem: injectivity of restriction on
   H^1} below, the subspace of $\Ext^1_{k_E[G_{K_{\infty}}]}(\overline{1},\overline{1})$
 arising from crystalline extensions of $\O_E(\varepsilon^{p-1})$ by $\O_E$ is
 also one-dimensional.

  On the other hand, by Proposition~\ref{prop:kisin-module-extensions}  there are no nontrivial
  extensions of $\barM(0;1)$ by $\barM(p-1;1)$ at all.  Why is this
  not a contradiction, given that the Kisin modules corresponding to
  $1$ and $\varepsilon^{p-1}$ reduce to $\barM(0;1)$ and
  $\barM(p-1;1)$ respectively?  The point is that although the functor
  $T_{\fS}$ is an  equivalence of categories, the inverse functor from
  lattices to Kisin modules need not be exact; and indeed the
  reduction mod $\varpi$ of the Kisin module
  corresponding to the lattice~$T$ of the previous paragraph turns out
  to be an extension of $\barM(p-1;1)$ by
  $\barM(0;1)$ rather than the reverse.
\end{example}

\para \label{para:maxmin} Suppose that $\rbar$ is
the reduction mod $\varpi$ of an $\O_E$-lattice in a pseudo-BT
representation $V$ of weight $\{r_i\}$, so that by
Theorem~\ref{thm:pseudo-BT-extensions} there exist $\barN$ and $\barP$
such that $\rbar|_{G_{K_{\infty}}}
  \simeq T_{\fS}(\barM)$ for some $\barM \in \EpsBT(\barN,\barP)$.  To
  prove that $\Wexplicit(\rbar)=\Wcris(\rbar)$, we wish to show that
  there exist crystalline characters $\chi',\chi : G_K \to
  \cO_E^{\times}$ and an extension $T$ of $\O_E(\chi)$ by
  $\O_E(\chi')$ such that
  $T[1/p]$ is pseudo-BT of weight $\{r_i\} $ and $\rbar \simeq T \otimes_{\O_E} k_E$.

 It is natural to try to argue by choosing $\chi'$ and $\chi$ so that their
 corresponding rank-one Kisin modules
 $\mathfrak{N},\mathfrak{P}$ are lifts of $\barN,\barP$ respectively,
  and then comparing the spaces $H^1_f(G_K,\cO_E(\chi^{-1} \chi'))$ and
  $\EpsBT(\barN,\barP)$ by a counting argument.  Unfortunately,
  Example~\ref{ex:structure-of-argumemt} shows that the Kisin module corresponding to an
  element of the first group may not reduce to an element of the
  latter set, and so an additional argument is needed.

  We consider instead all pairs of crystalline characters
  $\widetilde{\chi}',\widetilde{\chi}  : G_K \to \cO_E^{\times}$ with reductions
  $\chibar',\chibar$, and
  such that $\widetilde{\chi}'\oplus\widetilde{\chi}$ is pseudo-BT of
  weight $\{r_i\}$.  We will
  show that there is a preferred choice $\chi'_{\min},\chi_{\max}$
  (with Kisin modules $\mathfrak{N}_{\min},\mathfrak{P}_{\max}$ respectively) with the
  property that the reduction mod
  $\varpi$ of any element of any
  $H^1_f(G_K,\cO_E(\widetilde{\chi}^{-1}\widetilde{\chi}'))$ can be shown to occur as the image under
  $T_{\fS}$ of an
  element of $\EpsBT(\barN_{\min},\barP_{\max})$.  Then we \emph{can}
  proceed by comparing $H^1_f(G_K,\cO_E(\chi^{-1}_{\max} \chi'_{\min}))$
  with $\EpsBT(\barN_{\min},\barP_{\max})$.

 The construction of $\chi'_{\min},\chi_{\max}$
can be found in the proof of Theorem~\ref{thm:non-cyclotomic}.
 In the rest of this subsection we begin to carry out the above
strategy by proving the following proposition, which will
allow us to
compare the spaces $\EpsBT(\barN,\barP)$ (or at least the generic
fibres of the Kisin modules in those spaces) for certain varying choices of
$\barN$ and $\barP$.

\begin{prop}
  \label{prop:compare-extensions}
  Suppose that we are given Kisin modules
 $$\barN :=
  \barM(s_0,\ldots,s_{f-1};a) \quad \text{and} \quad \barP :=
  \barM(t_0,\ldots,t_{f-1};b)$$
as well  as
$$\barN' :=
  \barM(s'_0,\ldots,s'_{f-1};a) \quad \text{and} \quad  \barP' := \barM(t'_0,\ldots,t'_{f-1};b) $$
such that
\begin{itemize}
\item  there exist nonzero maps $\barP \to \barP'$ and $\barN' \to \barN$, and
\item $s_i + t_i = s'_i + t'_i = r_i + e-1$ and
  $\max\{s_i,t_i\},\max\{s'_i,t'_i\} \ge r_i$  for all $i$.
\end{itemize}
 For each
$\barM \in \EpsBT(\barN,\barP)$ there exists $\barM' \in
\EpsBT(\barN',\barP')$ such that $T_{\fS}(\barM) \cong
T_{\fS}(\barM')$.
\end{prop}

\begin{proof}
 The proof has two steps.  First we construct $\barM' \in
 \Ext^1(\barN',\barP')$ such that $T_{\fS}(\barM) \cong
T_{\fS}(\barM')$, and then we check that in fact $\barM' \in
\EpsBT(\barN',\barP')$.

We can choose
  bases $e_i,f_i$ of the $\barM_i$ so that $\varphi$ has the form
  \begin{eqnarray*}
    \varphi(e_{i-1}) & = &(b)_i u^{t_i} e_i \\
    \varphi(f_{i-1}) & = &(a)_i u^{s_i} f_i + y_i e_i
  \end{eqnarray*}
with the elements $y_i$ as in Theorem~\ref{thm:pseudo-BT-extensions}.
First we define the Kisin module $\barM'' \in \Ext^1(\barN',\barP)$ by
the formulas
 \begin{eqnarray*}
    \varphi(e''_{i-1}) & = &(b)_i u^{t_i} e''_i \\
    \varphi(f''_{i-1}) & = &(a)_i u^{s'_i} f''_i + y''_i e''_i
  \end{eqnarray*}
with $y''_i := u^{p(\alpha_{i-1}(\barN') - \alpha_{i-1}(\barN))} y_i$.
It is easy to check that there is a morphism $g : \barM'' \to \barM$
sending $e''_i \mapsto e_i$ and $f''_i \mapsto u^{\alpha_{i}(\barN')
  - \alpha_{i}(\barN)} f_i$ for all $i$, and since $g$ induces an
isomorphism $\barM''[1/u] \toisom \barM[1/u]$ of \'etale $\varphi$-modules, it induces an
isomorphism  $T_{\fS}(\barM) \toisom T_{\fS}(\barM'')$.

Next define the Kisin module $\barM' \in \Ext^1(\barN',\barP')$ by the
formulas
 \begin{eqnarray*}
    \varphi(e'_{i-1}) & = &(b)_i u^{t'_i} e'_i \\
    \varphi(f'_{i-1}) & = &(a)_i u^{s'_i} f'_i + y'_i e'_i
  \end{eqnarray*}
with $y'_i := u^{\alpha_i(\barP)-\alpha_i(\barP')} y''_i$.  Again, it is
easy to check that there is a morphism $g' : \barM'' \to \barM'$
sending $e''_i \mapsto u^{\alpha_i(\barP)-\alpha_i(\barP')}  e'_i$
and $f''_i \mapsto f'_i$, and that $g'$ induces an isomorphism
$T_{\fS}(\barM') \toisom T_{\fS}(\barM'')$.  Combining these
calculations shows that indeed $T_{\fS}(\barM) \cong T_{\fS}(\barM'')$.

It remains to check that $\barM' \in
\EpsBT(\barN',\barP')$.  Using the characterisation following
Definition~\ref{defn:pseudo-BT-classes}, we wish to show for each $i$
that there is a saturated element of $\varphi(\barM'_{i-1})$ that is
divisible by $u^{r_i}$ in $\barM'_i$.  When $t'_i \ge r_i$ this is
obvious (since $\varphi(e'_{i-1})$ will do), so we can suppose that $t'_i
< r_i$ and so $s'_i \ge r_i$.

Recall that $y'_i = u^{\alpha_i(\barP)-\alpha_i(\barP') +
  p(\alpha_{i-1}(\barN') - \alpha_{i-1}(\barN))} y_i$.  If
$\alpha_{i-1}(\barN')  > \alpha_{i-1}(\barN)$ then we are done, since
$u^p$ (hence also $u^{r_i}$) divides $y'_i$.  We may therefore assume
that  $\alpha_{i-1}(\barN')  = \alpha_{i-1}(\barN)$.  Since $s_i + t_i
= s'_i + t'_i$ for all $i$, it follows that  $\alpha_{i-1}(\barP')  =
\alpha_{i-1}(\barP)$ as well.    Note that whenever there is a map
$\barP\to \barP'$ with $\alpha_{i-1}(\barP')  =
\alpha_{i-1}(\barP)$, we must have $t_i \le t'_i$.  Indeed we have $t_i=p\alpha_{i-1}(\barP)-\alpha_i(\barP)$ and
  $t_i'=p\alpha_{i-1}(\barP')-\alpha_{i}(\barP')$, so that
  $t'_i-t_i=\alpha_i(\barP)-\alpha_i(\barP')\ge 0$ by Lemma~\ref{lem:map}.

In particular we also have $t_i < r_i$, and the assumption that $\barM
\in \EpsBT(\barN,\barP)$ implies that every term of $y_i$ has degree
at least $r_i$ except possibly for a term of degree $t_i$. (Note that
in the case that there is an extra term of degree $d_i$, we have
$d_i\ge s_i$, which is at least $r_i$ since $t_i<r_i$.)  Since $y_i$
divides $y'_i$, we see that every term of $y'_i$ has degree at least
$r_i$ except possibly for a term of degree
$$ t_i + \alpha_i(\barP)-\alpha_i(\barP') +
  p(\alpha_{i-1}(\barN') - \alpha_{i-1}(\barN)) = t_i + \alpha_i(\barP)-\alpha_i(\barP').$$
This quantity is easily seen to be congruent to $t'_i \pmod{p}$, so is
either equal to $t'_i$ or is at least $r_i$, and we are done, because
we can subtract a constant multiple of $\varphi(e'_{i-1})$ from
$\varphi(f'_{i-1})$ to obtain an element divisible by $u^{r_i}$.
\end{proof}

\subsection{Maximal and minimal models}
\label{ss:maximalminimal}

We now construct the maximal and minimal Kisin modules to which we alluded
in \ref{para:maxmin}.

\begin{lem}
  \label{lem:maximal-model}
  Fix integers $r_i \in [1,p]$.  Suppose that $\chibar : G_K \to k_E^{\times}$
  is a character and let $\cS$ be the space of rank-one
  Kisin modules $\barP = \barM(t_0,\ldots,t_{f-1};b)$ such that
  \begin{itemize}
  \item $T_{\fS}(\barP) = \chibar|_{G_{K_\infty}}$, and
  \item $t_i \in [0,e-1] \cup [r_i,r_i+e-1]$ for all $i$.
  \end{itemize}
  If $\cS \neq \varnothing$, then $\cS$ contains a maximal model.
  That is, there exists $\barP_{\max} \in \cS$ such that there is a
  nontrivial map $\barP \to \barP_{\max}$ for all $\barP \in \cS$.
\end{lem}

\begin{proof}  Assume that $\cS$ is non-empty.  Then there exists some
  $\barM(t'_0,\ldots,t'_{f-1};b) \in \cS$, and every other element of
  $\cS$ has the form $\barM(t_0,\ldots,t_{f-1};b)$ (i.e. the
  $b$ is the same).

 Write $\chibar|_{I_K} = \prod_i \omega_i^{m_i}$ with $m_i \in [0,p-1]$ and not
 all equal to $p-1$. For $0\le i\le f-2$, let $v_i$ be the $f$-tuple $(0,\ldots,-1,p,\ldots,0)$ with the $-1$ in
 the $i$th position (where the leftmost position is the zeroth), and similarly let $v_{f-1}$ be the $f$-tuple $(p,0,\ldots,0,-1)$. It is straightforward to see from Lemma~\ref{lem:functor-on-characters} that if
 $\barM(t_0,\ldots,t_{f-1};b) \in \cS$, then
\numequation \label{eq:sequation} (t_0,\ldots,t_{f-1}) = (m_0,\ldots,m_{f-1}) + \sum_{i=0}^{f-1} c_i
v_i \end{equation}
with $c_i \in \Z_{\ge 0}$ for all $i$.

 If we have $m_i \in \cI_i := [0,e-1] \cup [r_i,r_i+e-1]$
 for all $i$, then
 it is clear that $\barP_{\max} := \barM(m_0,\ldots,m_{f-1}; b)$ is a
 maximal model, e.g.\ because $\alpha_i(\barP_{\max}) \le
 \alpha_i(\barP)$ for all $\barP \in \cS$. This always holds for instance if $e \ge p$, or if
 $\chibar$ is unramified, so let us assume for the rest of the proof that $e \le p-1$ and that
 $\chibar$ is ramified, so that not every $m_i$ is equal to $0$.

With these additional hypothesis, we have $r_i + e-1 \le 2p-2$ for all $i$, and
also the integers $c_i$ in \eqref{eq:sequation} must all be $0$ or
$1$.   If $\barP = \barM(t_0,\ldots,t_{f-1};b)$ with integers
$c_i$ as in \eqref{eq:sequation}, write $J(\barP) := \{ i \, : \, c_i
\neq 0\}$.  To complete the proof, it suffices to show that there
exists a subset $J \subseteq \{0,\ldots,f-1\}$ such that
\begin{itemize}
\item if $\barP \in \cS$ then $J \subseteq J(\barP)$, and
\item there exists $\barP' \in \cS$ such that $J = J(\barP')$,
\end{itemize}
for then $\barP'$ is $\barP_{\max}$. We
construct the set $J$ as follows.  Note that if $x$ is a non-negative
integer with $x \not\in \cI_i$, then $x \ge e$, so that also $x+p
\not\in \cI_i$.

Let $K \subseteq \{0,\ldots,f-1\}$ be any set with the property that
if $\barP \in \cS$ then $K \subseteq J(\barP)$. Write
$(m'_0,\ldots,m'_{f-1}) = (m_0,\ldots,m_{f-1}) + \sum_{i\in K}
v_i$. Let  $\barP = \barM(t_0,\ldots,t_{f-1};b)$ be any element of
$\cS$ (here we use our assumption that $\cS$ is non-empty). Suppose
first that $i \in K$. Observe that $t_i \in
\{m'_i,m'_i+p\}$, and since $t_i \in \cI_i$ it follows by the last
sentence in the previous paragraph that $m'_i \in \cI_i$. Now suppose
instead that $i \not\in K$. Observe that $t_i \in
\{m'_i,m'_i+p,m'_i-1,m'_i+p-1\} \cap \cI_i$. If in fact $m'_i \not\in \cI_i$, it
follows by the last sentence in the previous paragraph that $t_i \neq
m'_i,m'_i+p$, and therefore $i \in J(\barP)$.  Write $\Delta(K) =\{
i \not\in K : m'_i \not\in \cI_i \}$; it follows that $K \cup \Delta(K)$ is
another  set that is contained in $J(\barP)$ for all $\barP \in \cS$.

Set $J_0 = \varnothing$, and iteratively define $J_i = J_{i-1} \cup
\Delta(J_{i-1})$ for $i \ge 1$.  Eventually this process stabilizes at some
$J_n$. By  construction $J_n \subseteq J(\barP)$ for all $\barP \in
S$. Again write $(m'_0,\ldots,m'_{f-1}) = (m_0,\ldots,m_{f-1}) + \sum_{i\in J_n}
v_i.$  We claim that $m'_i \in \cI_i$ for all~$i$, so that $J_n$ is
the desired set $J$:\ this is automatic if $i \in
J_n$ (by the first observation in the previous paragraph), and if $i
\not\in J_n$ follows from the fact that $\Delta(J_n) = \varnothing$.
\end{proof}

Evidently we must also have the following lemma, which essentially is dual to
the previous one.

\begin{lem}
  \label{lem:minimal-model}
  Fix integers $r_i \in [1,p]$.  Suppose that $\chibar' : G_K \to k_E^{\times}$
  is a character and let $\cS$ be the space of rank-one
  Kisin modules $\barN = \barM(s_0,\ldots,s_{f-1};a)$ such that
  \begin{itemize}
  \item $T_{\fS}(\barN) = \chibar'|_{G_{K_\infty}}$, and
  \item $s_i \in [0,e-1] \cup [r_i,r_i+e-1]$ for all $i$.
  \end{itemize}
  If $\cS \neq \varnothing$, then $\cS$ contains a minimal model.
  That is, there exists $\barN_{\min} \in \cS$ such that there is a
  nontrivial map $\barN_{\min} \to \barN$ for all $\barN \in \cS$.
\end{lem}

\begin{proof}
  Indeed, if $\chibar : G_K \to k_E^{\times}$ is a character such that $\chibar
  \chibar'|_{I_K} = \prod_i \omega_i^{r_i+e-1}$ and
  $\barM(t_0,\ldots,t_{f-1};b)$ is the maximal model for $\chibar$ given
  by the previous lemma, then the desired minimal model is given by $\barM(s_0,\ldots,s_{f-1};a)$ with $s_i = (r_i+e-1)-t_i$.
\end{proof}

Combining Proposition~\ref{prop:compare-extensions} with the above
lemmas, we obtain the following result, which was promised in~\ref{para:maxmin}.

\begin{prop}
  \label{prop:comparison}
  Fix characters $\chibar',\chibar : G_K \to k_E^{\times}$ and
  integers $r_i \in [1,p]$.  There exist rank-one Kisin modules
  $\barN_{\min}$ and $\barP_{\max}$ with the following property: if
  $\rbar \in \Ext^1_{k_E[G_K]}(\chibar,\chibar')$
  is the reduction modulo $\varpi$ of an $\cO_E$-lattice in a pseudo-BT
  representation of weight $\{r_i\}$, then $\rbar|_{G_{K_{\infty}}}
  \simeq T_{\fS}(\barM)$ for some $\barM \in
  \EpsBT(\barN_{\min},\barP_{\max})$.

 Moreover, if we write $\barN_{\min} =
  \barM(s_0,\ldots,s_{f-1};a)$ and $\barP_{\max} =
  \barM(t_0,\ldots,t_{f-1};b)$, then $\barN_{\min}$ and $\barP_{\max}$
  can be chosen so that for all $i$ we have $s_i + t_i = r_i + e-1$
  and $s_i,t_i \in
  [0,e-1] \cup [r_i,r_i+e-1]$.
\end{prop}

\begin{proof}
  If no such extensions $\rbar$ exist, then the Proposition is
  vacuously true, so we may suppose that some $\rbar$ exists as in the
  statement of the Proposition.  As in the proof of
  Theorem~\ref{thm:reducible-case}, there exist rank-one Kisin modules
  $\barN,\barP$ such that  $\rbar|_{G_{K_{\infty}}}
  \simeq T_{\fS}(\barM)$ for some $\barM \in \EpsBT(\barN,\barP)$.  In
  particular the set $\cS$ of Lemma~\ref{lem:maximal-model} is
  nonempty (it contains $\barP$, or possibly $\barN$ if
  $\rbar|_{G_{K_{\infty}}}$ is split), and similarly the set $\cS$ of
  Lemma~\ref{lem:minimal-model} is nonempty.
  Let $\barP_{\max}$ and $\barN_{\min}$ be the maximal and minimal
  models given, respectively, by those lemmas. (Note that these
  depend only on $\chibar$, $\chibar'$, and the integers $\{r_i\}$,
  and not on $\rbar$.)

 The fact that  $s_i,t_i \in
  [0,e-1] \cup [r_i,r_i+e-1]$ is given to us by
  Lemmas~\ref{lem:maximal-model} and~\ref{lem:minimal-model}, and the equality $s_i
  + t_i = r_i+e-1$ comes from the proof of Lemma~\ref{lem:minimal-model}. Now an application of Proposition~\ref{prop:compare-extensions}
  gives the claim in the first paragraph of the Proposition when $\rbar|_{G_{K_{\infty}}}$
    is non-split, while if $\rbar|_{G_{K_{\infty}}}$ is split we can
      take $\barM = \barP_{\max} \oplus \barN_{\min}$.
\end{proof}

\subsection{The non-cyclotomic case of the weight part of Serre's conjecture}
\label{sec:weight-part-serres}

We are now ready to prove the following result.

\begin{thm}
  \label{thm:non-cyclotomic}
  Suppose that $p\ge 3$ and that $\rbar : G_K \to \GL_2(k_E)$ is reducible, and write
  $\rbar \simeq \begin{pmatrix} \chibar' & * \\ 0 & \chibar \end{pmatrix}$.
  If $\chibar^{-1} \chibar' \neq \varepsilonbar$, then
  $\Wexplicit(\rbar) = \Wcris(\rbar)$.   Moreover,  these sets do not depend on
  our choice of embeddings~$\kappa_{i,0}$.
\end{thm}

\begin{proof}
  Suppose that $\sigma \in \Wcris(\rbar)$.  We may freely enlarge $E$
  so that $\rbar$ has a crystalline lift of Hodge type $\sigma$
  defined over $E$.  After twisting, we may assume
  that $\sigma$ has the shape $\otimes_{i} \Sym^{r_i-1} k^2 \otimes_{k,\kappabar_i}
  \Fpbar$ with integers $r_i \in [1,p]$, and the hypothesis that $\sigma \in \Wcris(\rbar)$ means that
  $\rbar$ is the reduction modulo $p$ of a lattice in a pseudo-BT
  representation of weight~$\{r_i\}$.
  Let  $\barN_{\min} =
  \barM(s_0,\ldots,s_{f-1};a)$ and $\barP_{\max} =
  \barM(t_0,\ldots,t_{f-1};b)$ be the rank-one Kisin modules given
  to us by Proposition~\ref{prop:comparison} applied to $\chibar'$,
  $\chibar$, and $\{r_i\}$.

    We construct a pair of crystalline characters
  $\chi'_{\min},\chi_{\max} : G_K \to \O_E^{\times}$ as follows.  If $t_i \ge r_i$, we take
  the ordered pair
  $(\HT_{\kappa_{i,0}}(\chi_{\max}),\HT_{\kappa_{i,0}}(\chi'_{\min}) ) = (r_i,0)$;
  and for $j > 0$ we take the pair
  $(\HT_{\kappa_{ij}}(\chi_{\max}),\HT_{\kappa_{ij}}(\chi'_{\min}) )$
  to be
  $(1,0)$ for exactly $t_i - r_i$ values of $j$ and to be $(0,1)$ for
  exactly $s_i$ values of $j$.  On the other hand if $t_i < r_i$, we take
  the pair
  $(\HT_{\kappa_{i,0}}(\chi_{\max}),\HT_{\kappa_{i,0}}(\chi'_{\min}) ) = (0,r_i)$;
  and for $j > 0$ we take the pair
  $(\HT_{\kappa_{ij}}(\chi_{\max}),\HT_{\kappa_{ij}}(\chi'_{\min}) )$
  to be
  $(1,0)$ for exactly $t_i$ values of $j$ and to be $(0,1)$ for
  exactly $s_i-r_i$ values of $j$.   It is then possible to choose the
  unramified parts of $\chi'_{\min},\chi_{\max}$ so that they reduce
  to $\chibar'$ and $\chibar$ respectively.

  Let us consider the extensions $T$ of $\O_E(\chi_{\max})$ by
  $\O_E(\chi'_{\min})$ such that $T[1/p]$ is crystalline.  Each of these is pseudo-BT of weight $\{r_i\}$, and
  so (by the defining property of $\barN_{\min},\barP_{\max}$ from Proposition~\ref{prop:comparison}) we must have
  $\overline{T}|_{G_{K_{\infty}}} \simeq T_{\fS}(\barM)$ for some $\barM \in
  \EpsBT(\barN_{\min},\barP_{\max})$.

  The space of crystalline extensions of
  $\chi_{\max}$ by $\chi'_{\min}$, which we identify with
  $H^1_f(G_K,E(\chi_{\max}^{-1} \chi'_{\min}))$, has dimension equal to
    the number of labeled Hodge--Tate weights of $\chi'_{\min}$ that exceed
    the corresponding weight of $\chi_{\max}$.  This is precisely
 $$  d = \sum_{i \, :\, t_i \ge r_i} s_i +
\sum_{i \, : \, t_i < r_i} (1 + s_i - r_i).$$
It follows as in \cite[Lem.~9.3]{GLS-BDJ} that the image of $H^1_f(G_K,\O_E(\chi_{\max}^{-1}
\chi'_{\min}))$ in $H^1(G_K,\chibar^{-1} \chibar')$ has dimension $d$ if
$\chibar \neq \chibar'$ and dimension $d+1$ if $\chibar=\chibar'$.

By Lemma~\ref{lem: injectivity of restriction on H^1} below, the
restriction map  $H^1(G_K,\chibar^{-1} \chibar') \to
H^1(G_{K_{\infty}},\chibar^{-1} \chibar')$  is injective.  (The
application of Lemma~\ref{lem: injectivity of restriction on H^1} is the
only place in the argument that we need our assumption that
$\chibar^{-1} \chibar' \neq \varepsilonbar$.)  It follows that the
number of elements of $H^1(G_{K_{\infty}},\chibar^{-1} \chibar')$ that
come from a crystalline extension of
  $\chi_{\max}$ by $\chi'_{\min}$ is exactly $|k_E|^{d+\delta}$, where
  $\delta=0$ if $\chibar
  \neq \chibar'$ and $\delta=1$ otherwise.

On the other hand, we know from Proposition~\ref{prop:comparison} that
the number of elements of
$H^1(G_{K_{\infty}},\chibar^{-1} \chibar')$ that come from the reduction
mod $p$ of some pseudo-BT representation of weight $\{r_i\}$ is
\emph{at most} $\#\EpsBT(\barN_{\min},\barP_{\max})$.  One easily checks
by counting, in the explicit description of the extensions in
Theorem~\ref{thm:pseudo-BT-extensions}, the number of terms in each $y_i$ that are permitted to be non-zero
(and noting that if
$\chibar=\chibar'$ there must exist a map
$\barN_{\min}\to\barP_{\max}$, by the maximality of $\barP_{\max}$) that
$\#\EpsBT(\barN_{\min},\barP_{\max}) = |k_E|^{d+\delta}$ as well.

In particular every  element of
$H^1(G_{K_{\infty}},\chibar^{-1} \chibar')$ that comes from the reduction
mod $p$ of some pseudo-BT representation of weight $\{r_i\}$ must in
fact come from a lattice in a crystalline extension of $\chi_{\max}$
by $\chi'_{\min}$.  Applying Lemma~\ref{lem: injectivity of
  restriction on H^1} again, the same must be true of every  element of
$H^1(G_{K},\chibar^{-1} \chibar')$ that comes from the reduction
mod $p$ of some pseudo-BT representation of weight $\{r_i\}$.  Since
in particular $\rbar$ is such an element, we deduce that $\sigma \in
\Wexplicit(\rbar)$, as desired.

Finally, note that the Kisin modules $\barN_{\min}$ and $\barP_{\max}$
do not depend on the choice of embeddings $\kappa_{i,0}$.  Then the
independence of $\Wexplicit(\rbar)$ from the choice of embeddings
$\kappa_{i,0}$ follows from the characterisation that $\sigma \in
\Wexplicit(\rbar)$ if and only if $\rbar|_{G_{K_\infty}} =
T_{\fS}(\barM)$ for some $\barM \in \EpsBT(\barN_{\min},\barP_{\max})$.
\end{proof}

\begin{lem}
  \label{lem: injectivity of restriction on H^1}Let $\chibar:G_K\to
  k_E^\times$ be a continuous character. If
  $\chibar\ne\varepsilonbar$, then the restriction
  map \[H^1(G_K,\chibar)\to H^1(G_{K_\infty},\chibar)\]is
  injective. If $\chibar=\varepsilonbar$, then the kernel of the
  restriction map is the tr\`es ramifi\'ee line determined by the
  fixed uniformiser $\pi$.
\end{lem}
\begin{proof}
  If $\chibar\ne 1,\varepsilonbar$, this is a special case
  of~\cite[Lem.~7.4.3]{emertongeesavitt}. If $\chibar=1$, then
  $H^1(G_K,\chibar)=\Hom(G_K,k_E)$ and
  $H^1(G_{K_\infty},\chibar)=\Hom(G_{K_\infty},k_E)$, so if the kernel
  of the restriction map is nonzero, there must be a Galois extension
  of $K$ of degree $p$ contained in $K_\infty$. This can only happen if
  $K$ contains a primitive $p$th root of unity, in which case
  $\varepsilonbar=1$, so $\chibar=\varepsilonbar$.

Finally, suppose that $\chibar=\varepsilonbar$. Kummer theory
identifies the restriction map with the natural
map \[K^\times/(K^\times)^p\to
K_\infty^\times/(K_\infty^\times)^p,\]and the kernel of this map is
evidently generated by $\pi$.
\end{proof}

\section{The weight part of Serre's conjecture III: the general case}
\label{sec:appl-weight-part-general-case}

\subsection{$(\varphi,\Ghat)$-modules}In order to complete our arguments in the remaining case, we will need
to make use of the second author's theory of
$(\varphi,\Ghat)$-modules. We refer the reader
to~\cite[\S5.1]{GLS-BDJ} (specifically, from the start of that section up to the
statement of Theorem~5.2)
for the definitions and notation that we will use, as well as to
\cite[(4.8)]{GLS-BDJ} for the definition of the operator $\tau$.

\para Consider a $(\varphi,\Ghat)$-module with natural $k_E$-action
$\hat\barM$, sitting in an extension of $(\varphi,\Ghat)$-modules with
natural $k_E$-action \[0 \to \hat\barP \to\hat \barM \to\hat \barN \to
0,\] where the underlying Kisin modules $\barN,\barP$ are given by
$\barN =\barM(s_0 , \dots, s_{f-1}; a) $ and $\barP = \barM(t_0 ,
\dots, t_{f-1}; b ) $ for some $a,b$, with
$\{s_i,t_i\}=\{r_i+x_i,e-1-x_i\}$, and the underlying extension of
Kisin modules is in $\EpsBT(\barN,\barP)$. We say that such a
$(\varphi,\Ghat)$-module is of \emph{reducible pseudo-BT type and
  weight $\{r_i\}$} if for all $x\in\barM$ there exist $\alpha\in R$
and $y\in R\otimes_{\varphi,\cS}\barM$ such that $\tau(x)-x=\alpha y$
and $v_R(\alpha)\ge \frac{p}{p-1}+\frac{p}{e}$.

\begin{lem}
  \label{lem: reduction of a pot BT rep gives phi Ghat of required
    type} Suppose that $p\ge 3$, and that $\rbar:G_K\to\GL_2(k_E)$ is reducible and arises as
  the reduction mod $p$ of a pseudo-BT representation $r$ of weight
  $\{r_i\}$. Then there is a $(\varphi,\Ghat)$-module with natural
  $k_E$-action $\hat\barM$ such that $\hat\barM$ is of reducible
  pseudo-BT type and weight $\{r_i\}$, and $\hat{T}(\hat\barM)\cong\rbar$.
\end{lem}
\begin{proof}
  Let $\hat\barM$ be the $(\varphi,\Ghat)$-module arising as the
  reduction mod $p$ of the $(\varphi,\Ghat)$-module corresponding to
  $r$ by \cite[Thm.~5.2(2)]{GLS-BDJ}. Then the underlying Kisin
  module $\barM$ is of the required kind by
  Theorem~\ref{thm:pseudo-BT-extensions}, and this extension
  can be extended to an extension of $(\varphi,\Ghat)$-modules
  by~\cite[Lem.~5.5]{GLS-BDJ}. Finally, the claim about the action of
  $\tau$ is immediate from~\cite[Cor.~5.10]{GLS-BDJ}.
\end{proof}

We have the following generalisation of~\cite[Lem.~8.1]{GLS-BDJ}.

\begin{lemma}\label{lem: hat G is determined by phi} Suppose that
  $p\ge 3$ and take
  $\barM\in\EpsBT(\barN,\barP)$. Except possibly for the case that $r_i =  p$ and $t_i=0$ for all $i
  = 0 , \dots, f-1$, there is at most one $(\varphi,\Ghat)$-module $\hat\barM$  of reducible
  pseudo-BT type and weight $\{r_i\}$ with underlying Kisin module
  $\barM$.
\end{lemma}
\begin{proof}
  We follow the proof of~\cite[Lem.~8.1]{GLS-BDJ}. Since by
  definition $\overline \M$ is contained in the $H_K$-invariants of
  $\hat{\overline{\M}}$, it suffices to show that the $\tau$-action on
  $ \hR \otimes_{\varphi, \fS} \barM$ is uniquely determined.  As
  usual we write
  $e_i,f_i$ for a basis of $\barM_i$ as given by
  Theorem~\ref{thm:pseudo-BT-extensions}.  We can
  write $$\tau (e_{i-1}, f_{i-1}) = (e_{ i-1 },
  f_{i-1}) \begin{pmatrix} \delta_{i} & \beta_{i} \\ 0 &
    \gamma_{i} \end{pmatrix}$$ with $\delta_i, \beta _i, \gamma_i \in
  (\hR/ p \hR) \otimes_{\Fp} k_E \subset R \otimes_{\Fp} k_E $.  If
  $\zeta \in R \otimes_{\Fp} k_E$ is written $\zeta = \sum_{i=1}^{n}
  y_i \otimes z_i$ with $z_1,\ldots,z_n \in k_E$ linearly independent
  over $\Fp$, write $v_R(\zeta) = \min_i \{ v_R(y_i)\}.$ By
  assumption, we have $v_R(\delta_i-1),
  v_R(\gamma_i-1),v_R(\beta_i) \ge \frac{p}{p-1}+\frac{p}{e}$ for all
  $i$.

Recalling
 that $\barM$
  is regarded as a $\varphi(\fS)$-submodule of $\hR \otimes_{\varphi,
    \fS}\barM$,
we may write $\varphi (e_{ i-1}, f_{ i-1 }) = (e_{ i }, f_{i}) \varphi (A_i)$ with $A_i = \begin{pmatrix} (b)_i u ^{t_i} &  x_i \\ 0 & (a)_i u ^{s_i} \end{pmatrix} $. Since $\varphi$ and $\tau$ commute, we have
$$ \varphi(A_i) \begin{pmatrix} \varphi(\delta_{i})  & \varphi(\beta_{i})  \\ 0 & \varphi(\gamma_{i}) \end{pmatrix} =  \begin{pmatrix} \delta_{i+1} & \beta_{i+1} \\ 0 & \gamma_{i+1} \end{pmatrix} \tau (\varphi(A_i)) . $$
 We obtain the following formulas:
\numequation \label{eq: unique tau 1}
   u ^{pt_i } \varphi (\delta_{i})= \delta_{i+1} (\ue u)
  ^{pt_i},  \ \   u ^{ps_i } \varphi (\gamma_{i})= (\ue u)
  ^{ps_i} \gamma_{i+1}
\end{equation}
 and
 \numequation\label{eq: unique tau 2}  (b)_i u ^{pt_i }
   \varphi (\beta_{i}) + \varphi(x_i)  \varphi(\gamma_{i}) =
   \delta_{i+1}\tau( \varphi (x _i)) + (a)_i (\ue u )^{ps_i} \beta_{i+1}
 \end{equation}
where for succinctness we have written $(a)_i$, $(b)_i$ in lieu of
$1\otimes (a)_i$, $1\otimes (b)_i$ in the preceding equation.

Let $\eta \in R$ be the element
defined in~\cite[Lem.~6.6(2)]{GLS-BDJ}, so that $\varphi^f(\eta) =
\ue \eta$. ($K/\Qp$ is assumed to be unramified throughout~\cite[\S6]{GLS-BDJ}, but it is easily checked that~\cite[Lem.~6.6(2)]{GLS-BDJ} remains valid with the same proof in our
setting.) From~\eqref{eq: unique tau 1} we see that $\varphi^f(\delta_i) =
\delta_i \ue^{\sum_{j=0}^{f-1} p^{f-j} t_{i+j}}$, and now \cite[Lem.~6.6(2)]{GLS-BDJ}
together with the requirement that
$v_R(\delta_i-1) > 0$
implies that
$$ \delta_i = \eta^{\sum_{j=0}^{f-1} p^{f-j} t_{i+j}} \otimes 1$$ for
all $i$.  Similarly we must have $\gamma_i = \eta^{\sum_{j=0}^{f-1}
  p^{f-j} s_{i+j}} \otimes 1$ for all $i$.  So at least the $\delta_i,\gamma_i$
are uniquely determined.

 Now suppose that there exists some other extension of $\barM$ to
a $(\varphi, \hat
 G)$-module $\hat \barMp$.
Then the $\tau$-action on $\hat \barMp$ is given by some $\delta'_i$,
 $\beta'_{i}$ and $\gamma'_i$ that also satisfy~\eqref{eq: unique tau
   1} and~\eqref{eq: unique tau 2}, and indeed we have already seen that
 $\delta'_i = \delta_i$ and $\gamma'_i = \gamma_i$.

Let $\tilde \beta_i = \beta_i - \beta'_i$.  Taking the difference
between~\eqref{eq: unique tau 2} for $\hat \barM$ and $\hat \barMp$
gives
$$ (b)_i u^{pt_i} \varphi(\tilde \beta_{i}) = (a)_i (\ue u)^{p
  s_i} \tilde \beta_{i+1},$$
which implies that
$$b u^{\sum_{j = 0}^{f-1} u^{p^{f-j} t_{i+j}}} \varphi^f (\tilde
\beta_i) = a (\ue u)^{\sum_{j = 0}^{f-1} u^{p^{f-j} s_{i+j}}} \tilde
\beta_i.$$
Considering the valuations of both sides, and recalling that $v_R(\underline \pi)= 1/e$, we see that if $\betat_i\ne
0$ then \numequation \label{eq:beta-valuation}  v_R(\tilde
\beta_i) = \frac{1}{e(p^f-1)} \sum_{j=0}^{f-1} p^{f-j} (s_{i+j} - t_{i+j}). \end{equation}
But since $s_i - t_i\le r_i+e-1$ is at most $p+e-1$ with equality
if and only if $t_i =0$ and $ r_i = p$,  the right-hand side
of~\eqref{eq:beta-valuation} is at most $\frac{p}{e(p-1)}(p+e-1)=\frac{p}{p-1}+\frac{p}{e}$ with
equality if and only if  $t_i =0$ and $ r_i = p$ for all $i$.  In particular,
since $v_R(\beta_i),v_R(\beta'_i) \ge \frac{p}{p-1}+\frac{p}{e}$, either
$\beta_i = \beta'_i$ for all $i$, or else $t_i =0$ and $ r_i = p$  for all $i$,
as required.
\end{proof}

\begin{prop}
  \label{prop: extending Ghat action to maximal model} Suppose that $p
  \ge 3$ and let $\hat\barM$
  be a $(\varphi,\Ghat)$-module  of reducible
  pseudo-BT type and weight $\{r_i\}$ with underlying Kisin module
  $\barM$. Suppose we have $\barN',\barP'$ as in the statement of
  Proposition~\ref{prop:compare-extensions} and let $\barM'$ be the Kisin module provided by
 that proposition. Then there is a
  $(\varphi,\Ghat)$-module $\hat\barM'$  of reducible
  pseudo-BT type and weight $\{r_i\}$ with underlying Kisin module
  $\barM'$, such that $\hat{T}(\hat\barM)\cong\hat{T}(\hat\barM')$.
\end{prop}
\begin{proof}
  From the proof of Proposition~\ref{prop:compare-extensions}, we see
  that there is a Kisin module $\barM''$ and morphisms
  $g:\barM''\to\barM$, $g':\barM''\to\barM'$, both of which induce
  isomorphisms after inverting $u$. Using these isomorphisms, the
  $\Ghat$-action on $\hR\otimes_{\varphi,\fS}\barM$ induces a
  $\Ghat$-action on $\hR\otimes_{\varphi,\fS}\barM'[1/u]$, and it is
  enough to show that this preserves
  $\hR\otimes_{\varphi,\fS}\barM'$ and makes it a
  $(\varphi,\Ghat)$-module of reducible
  pseudo-BT type and weight $\{r_i\}$.

Since $H_K$ acts trivially on $u$ and $\barM$, and since
$\tau(u)=\underline{\epsilon}u$ and $(\epsilon-1)\in I_+$, we see that
we only need to check that
$\tau(\barM')\subset\hR\otimes_{\varphi,\fS}\barM'$, and that for all $x\in\barM'$ there exists $\alpha\in R$
and $y\in R\otimes_{\varphi,\fS}\barM'$ such that $\tau(x)-x=\alpha y$
and $v_R(\alpha)\ge \frac{p}{p-1}+\frac{p}{e}$.

Take bases $e_i,f_i$ of $\barM_i$ and $e'_i,f'_i$ of $\barM'_i$ as in
the proof of Proposition~\ref{prop:compare-extensions}. Writing \[\tau (e_{i-1}, f_{i-1}) = (e_{ i-1 },
  f_{i-1}) \begin{pmatrix} \delta_{i} & \beta_{i} \\ 0 &
    \gamma_{i} \end{pmatrix},\] an easy calculation shows that \[\tau (e'_{i-1}, f'_{i-1}) = (e'_{ i-1 },
  f'_{i-1})\begin{pmatrix} \delta'_{i} & \beta'_{i} \\ 0 &
    \gamma'_{i} \end{pmatrix},\] where
\[\delta'_i=\delta_{i}\underline{\epsilon}^{p(\alpha_{i-1}(\barP')-\alpha_{i-1}(\barP))},\]
\[\beta'_i=\beta_{i}(u^{p(\alpha_{i-1}(\barN')-\alpha_{i-1}(\barN)+\alpha_{i-1}(\barP)-\alpha_{i-1}(\barP'))}\otimes
1)\underline{\epsilon}^{p(\alpha_{i-1}(\barN')-\alpha_{i-1}(\barN))},\]\[\gamma'_i=\gamma_{i}\underline{\epsilon}^{p(\alpha_{i-1}(\barN')-\alpha_{i-1}(\barN))}.\]
(The factors of $p$ in the above exponents come from the fact that $\tau$ acts on the
left-hand side of the twisted tensor product $\hR\otimes_{\varphi,\fS} \barM'$.)
Now,
$\underline{\epsilon}$ is a unit, so it is enough to check that the
exponent of $u$ in the expression for $\beta'_i$ is non-negative; but
this is immediate from Lemma~\ref{lem:map}.
\end{proof}

We are now in a position to prove that
$\Wexplicit(\rbar)=\Wcris(\rbar)$ in the reducible cyclotomic case,
and so to deduce this equality in full generality (for $p \neq 2$).

\begin{thm}
  \label{thm:cyclotomic}
  Suppose that $p\ge 3$ and that $\rbar : G_K \to \GL_2(k_E)$ is
  a continuous representation. Then
  $\Wexplicit(\rbar) = \WBT(\rbar) = \Wcris(\rbar)$.  Moreover,  these sets do not depend on
  our choice of embeddings~$\kappa_{i,0}$.
\end{thm}
\begin{proof}
 Recall that we have inclusions
$\Wexplicit(\rbar) \subseteq \WBT(\rbar) \subseteq \Wcris(\rbar)$
by~\cite[Cor.~4.5.7]{GeeKisin}, so
it
is enough to check that $\Wcris(\rbar) \subseteq \Wexplicit(\rbar)$.
Suppose that $\sigma \in \Wcris(\rbar)$.
 As in the proof of Theorem~\ref{thm:non-cyclotomic}, we
  may assume that $\sigma$ has the shape $\otimes_{i} \Sym^{r_i-1} k^2
  \otimes_{k,\kappabar_i} \Fpbar$ with integers $r_i \in [1,p]$, and that
  $\rbar$ is the reduction modulo~$p$ of a lattice in a pseudo-BT
  representation of weight~$\{r_i\}$.  By Theorem~\ref{thm:our main
    semisimple local result in Bexplicit language} we may assume that $\rbar$ is an
  extension of $\chibar$ by $\chibar'$; by
  Theorem~\ref{thm:non-cyclotomic}, we may assume that
  $\chibar^{-1}\chibar'=\varepsilonbar$.

  Suppose first that not all of the $r_i$ are equal to $p$. We
  may replace the appeals to Lemma~\ref{lem: injectivity of restriction
    on H^1} in the proof of Theorem~\ref{thm:non-cyclotomic} with
  appeals to Lemma~\ref{lem: reduction of a pot BT rep gives phi Ghat
    of required type}, Lemma~\ref{lem: hat G is determined by phi} and
  Proposition~\ref{prop: extending Ghat action to maximal model}; then
 the count of Kisin modules remains valid, and the
  argument goes through as before. (Recall that the only place that
  the assumption that $\chibar^{-1}\chibar'\ne\varepsilonbar$ was used
  in the proof of Theorem~\ref{thm:non-cyclotomic} was in the appeals
  to Lemma~\ref{lem: injectivity of restriction on H^1}. In particular
  the construction of $\chi_{\max}$ and $\chi'_{\min}$ carries over to
  the case that $\chibar^{-1} \chibar' = \varepsilon$.)

  Finally, in the case that all of the $r_i$ are equal to $p$, the
  same argument applies unless $\rbar$ has a model where all of the
  $t_i=0$. In this case, we see that the character $\chi_{\max}$ in
  the proof of Theorem~\ref{thm:non-cyclotomic} is unramified, while
  $\HT_{\kappa_{i,0}}(\chi'_{\min})=p$ and
  $\HT_{\kappa_{i,j}}(\chi'_{\min})=1$ if $j>0$.

Then every extension of $\chi_{\max}$ by an unramified twist of
$\chi'_{\min}$ is automatically crystalline. So, it suffices to show
that \emph{any} extension of $\chibar$ by $\chibar'$ lifts to an
extension of $\chi_{\max}$ by a twist of $\chi'_{\min}$ by an unramified
character with trivial reduction. This may be proved by exactly the
same argument as~\cite[Lem.~9.4]{GLS-BDJ}
(\emph{cf.}~\cite[Prop.~5.2.9]{GLS11}, which proves the claim in
the case that $K/\Qp$ is totally ramified).
\end{proof}

\begin{remark}
  \label{rem:unramified-twists}
  Suppose that $\rbar$ is an extension of $\chibar$ by $\chibar'$, and
  let $\chi_{\max},\chi'_{\min}$ be the crystalline lifts of
  $\chibar,\chibar'$ constructed in the proof of
  Theorem~\ref{thm:non-cyclotomic}.  Recall that there is a choice in
  this construction, namely $\chi_{\max}$ and $\chi'_{\min}$ are only
  specified up to twist by unramified characters with trivial
  reduction.  It follows immediately from our arguments that this
  choice is immaterial, i.e.\ that the image of
  the map $H^1_f(G_K,\O_E(\chi_{\max}^{-1}\chi'_{\min})) \to
  H^1(G_K,\chibar^{-1}\chibar')$ does not depend on the particular
  choice of  $\chi_{\max},\chi'_{\min}$,
  \emph{except} for the case where
  $r_i = p$ for all $i$ and $\chi_{\max}$ is unramified (\textit{cf.}\ \cite[Rem.~3.10]{bdj}).
\end{remark}

\subsection{Conclusion of the proof of the weight part of Serre's conjecture.}
We now extend the results of
Section~\ref{sec:appl-weight-part}, using our analysis of the
extension classes of reducible lifts to complete the proof of the
weight part of Serre's conjecture.

Continue to assume that $p>2$, let $F$ be a totally real field, and
let $\rhobar:G_F\to\GL_2(\Fpbar)$ be continuous, irreducible, and
modular. Again, let $D$ be a quaternion algebra with centre $F$, which
is split at all places dividing $p$ and at at most one infinite
place. The main global result of this paper is the following.

\begin{thm}\label{thm:our main global result}Assume that $p>2$, that $\rhobar$ is modular and compatible
  with $D$, that $\rhobar|_{G_{F(\zeta_p)}}$ is irreducible, and if
  $p=5$ assume further that the projective image
  of $\rhobar|_{G_{F(\zeta_p)}}$ is not isomorphic to $A_5$.

  Then $\rhobar$ is modular for $D$ of weight $\sigma=\otimes_{v|p}\sigma_v$ if and
  only if $\sigma_v\in\Wexplicit(\rhobar|_{G_{F_v}})$ for all $v|p$.
\end{thm}
\begin{proof}
The result is immediate from Theorems~\ref{thm:GK main
  result} and~\ref{thm:cyclotomic}.
\end{proof}

\subsection{Dependence on the restriction to inertia}
\label{sec:dependence-inertia}

We conclude the paper by checking (when $p \ge 3$) that the
local weight set $\Wexplicit(\rbar) = \Wcris(\rbar)$ depends only on
the restriction to inertia $\rbar|_{I_K}$.  This statement is part of
the formulation of some versions of the weight part of Serre's
conjecture.  For instance when $K/\Qp$ is unramified, the definition
of the  weight set
$\WBDJ(\rbar)$ in \cite{bdj} is modified in certain cases to ensure that $\WBDJ(\rbar)$ depends
only on $\rbar|_{I_K}$ (\emph{cf.}~the definition immediately
preceding \cite[Rem.~3.10]{bdj}).  An immediate consequence of the following
result (and its proof) is that the local weight set defined in
\cite{bdj} is equal to the local weight set considered in this paper
(again, when $p \ge 3$).

\begin{prop}
  \label{prop:dependence-inertia}
  Suppose $p \ge 3$.  Let $\rbar,\rbar' : G_K \to \GL_2(k_E)$ be two
  continuous representations with $\rbar|_{I_K} \simeq
  \rbar'|_{I_K}$.  Then $\Wexplicit(\rbar) = \Wexplicit(\rbar')$.
\end{prop}

\begin{proof}
When $K/\Qp$ is unramified, the corresponding statement for the weight
set $\WBDJ(\rbar)$ is proved in \cite[Prop.~3.13]{bdj}.  The proof
carries over to this context nearly word for word, except that in the
special case where $\rbar^{\mathrm{ss}}\simeq  \chibar \oplus \chibar$ is
scalar and the class in $H^1(G_K,k_E)$ defining $\rbar$ as an
extension of $\chibar$ by $\chibar$ is unramified, it is built into the definitions that $\WBDJ(\rbar) =
\WBDJ(\chibar \oplus \chibar)$.   (In the notation of \cite{bdj}, this is the containment $L_{\mathrm{ur}} \subseteq
L_\alpha$ in this case.)  In our context, it remains to show that
$\Wexplicit(\rbar) = \Wexplicit(\chibar \oplus \chibar)$.

We must prove that if $\chibar \oplus \chibar$ has a pseudo-BT
lift of weight $\{r_i\}$, then so does $\rbar$.  Let $\barN_{\min}$
and $\barP_{\max}$ be the rank-one Kisin modules given to us by (the
proof of) Proposition~\ref{prop:comparison}, with corresponding
parameters $s_0,\ldots,s_{f-1}$ and $t_0,\ldots,t_{f-1}$ as in that
Proposition.  In the special case where $r_i = p$ and $t_i = 0$ for
all $i$ (which is only possible when $\varepsilonbar = 1$ on $G_K$), the proof of Theorem~\ref{thm:cyclotomic} already shows that
\emph{every} extension of $\chibar$ by $\chibar$ has a pseudo-BT lift
of weight $\{r_i\}$, so for the rest of the proof we assume that we are not in this case.

Write $\calE$ for the set $\EpsBT(\barN_{\min},\barP_{\max})$.
Note that by construction (since $\barP_{\max}$ and $\barN_{\min}$
have the same generic fibre $\chibar$) there exists a nonzero map
$\barN_{\min} \to \barP_{\max}$.  We can therefore define a nontrivial
element
$\barM \in \mathcal{E}$ by taking $y_i = 0$ for all
$i$  (in the notation of
Theorem~\ref{thm:pseudo-BT-extensions}), except that $y_0$ is taken to
have a nonzero term of degree $d := s_0 +
\alpha_0(\barN_{\min}) - \alpha_0(\barP_{\max})$.  Indeed we obtain a
line's worth of such elements.

Let $K_p$ be the unramified extension of $K$ of degree $p$, with
residue field $k_p$, and suppose
without loss of generality that $k_p$ embeds into $k_E$.  Define
$\barM' = k_p \otimes_k \barM$, with $\varphi$ extended
$k_p$-semilinearly to $\barM'$.  It is straightforward to check that the
Kisin module $\barM'$ is split.  We briefly indicate how to see this.
Write the Kisin module $\barM'$ as an
extension as in Proposition~\ref{prop:kisin-module-extensions}. The
field $K_p$ has inertial degree $pf$.  The extension parameters $y'_i$
defining $\barM'$ are all zero, except for terms of degree~$d$ in
$y'_0,\ldots,y'_{(p-1)f}$, all with the same coefficient.  Since
terms of this degree are part of a loop  (in the terminology of the
proof of \cite[Prop.~7.4]{GLS-BDJ}) there is a
change of variables which replaces $y'_{if}$ with $0$ for each $0 < i < p$ and
replaces $y'_0$ with $y'_0 + \ldots y'_{(p-1)f }= p y'_0 = 0$.

Let $T$ be a lattice in a pseudo-BT representation of $G_K$ of weight
$\{r_i\}$ such that $\overline{T} |_{G_{K_{\infty}}} \simeq
T_{\fS}(\barM)$.  Since $T_{\fS}$ is faithful on $\calE$, the representation
$\overline{T}$ is non-split.   By Lemmas~\ref{lem: reduction of a pot BT rep gives phi Ghat of required
    type} and~\ref{lem: hat G is determined by phi} there is a unique
  $(\varphi,\Ghat)$-module $\hat\barM$ of reducible pseudo-BT type and
  weight $\{r_i\}$ with $\hat{T}(\hat\barM) \simeq \overline{T}$.
  Evidently we have $\overline{T} |_{G_{(K_p)_{\infty}}}  \simeq
T_{\fS}(\barM')$.  Since $T|_{G_{K_p}}$ is still pseudo-BT but $\barM'$ is split,
another application of Lemma~\ref{lem: hat G is determined by phi} shows that
  $\overline{T}|_{G_{K_p}}$ itself must be split.  (We remark that if
  $\varepsilonbar \neq 1$ on $G_K$, the fact that
  $\overline{T}|_{G_{K_p}}$  is split can be deduced more easily from
  Lemma~\ref{lem: injectivity of restriction on H^1}.)   It follows
  that the extension class in $H^1(G_K,k_E)$ defining $\overline{T}$
  as an extension of $\chibar$ by $\chibar$ lies in the kernel of
  the restriction map $H^1(G_K,k_E) \to H^1(G_{K_p},k_E)$, i.e.\ it is
  unramified.  Since the unramified subspace of $H^1(G_K,k_E)$ is a
  line, we have $\overline{T} \simeq \rbar$, and $\rbar$ has the
  desired pseudo-BT lift (namely $T$).
\end{proof}

\appendix

\section{Corrigendum to \cite{GLS-BDJ}}
\label{sec:corrigendum}

We take this opportunity to correct a minor but unfortunate mistake in
\cite{GLS-BDJ}.  In the sentence preceding the published version of \cite[Thm.~4.22]{GLS-BDJ},
we write that when we regard $\fM$ as a $\varphi(\fS)$-submodule of
$\fM^*$, we are regarding $\fM_{s+1}$ as a submodule of $\fM_s^*$.
As noted in Section~\ref{sec:pseudo-BT} of the present paper, this
should be
$\fM_{s-1}$ rather than $\fM_{s+1}$.  Although the main results
of \cite{GLS-BDJ} about the weight part of Serre's conjecture are unaffected, the
statement of  \cite[Thm.~4.22]{GLS-BDJ} and many of
the ensuing Kisin module and $(\varphi,\Ghat)$-module formulas in
\cite[\S\S6--8]{GLS-BDJ}, have some indices that are off by one.
(This issue does not affect the structure or content of any of the arguments in the paper, only the statements.)

For instance, where we write in \cite[Thm.~4.22]{GLS-BDJ} that
$$\varphi(e_{1,s},\ldots,e_{d,s}) = (e_{1,s+1},\ldots,e_{d,s+1}) X_s
\Lambda_s Y_s $$
we should instead have
$$\varphi(e_{1,s-1},\ldots,e_{d,s-1}) = (e_{1,s},\ldots,e_{d,s}) X_s
\Lambda_s Y_s ,$$
and the other corrections are all of a similar nature.  Corrected
versions of the paper are available on our websites and on the arXiv.
(Alternately, many of the calculations in \cite[\S\S6--8]{GLS-BDJ} are generalised
by results in Sections~\ref{sec:reductions-mod-p},\
\ref{sec:non-cyclotomic}, and~\ref{sec:appl-weight-part-general-case}
of the present paper.)

\bibliographystyle{amsalpha}
\bibliography{savitt}

\end{document}